%% file: main.tex
\pdfoutput=1
\documentclass[a4paper,10pt,twocolumn]{article}
\usepackage[T1]{fontenc}
\usepackage{lmodern}
\usepackage[nopar]{lipsum}
\usepackage{caption}
\usepackage{authblk}
\usepackage[top=2cm, bottom=2cm, left=2cm, right=2cm]{geometry}
\usepackage{fancyhdr}
\usepackage{url}
\usepackage{afterpage,float}
\usepackage{graphicx}
\usepackage{amssymb, gensymb}
\usepackage{amsmath}
\usepackage[authormarkuptext=id]{changes}
\usepackage{soul}
\usepackage{tikz}
\usepackage{pgfplots}
\pgfplotsset{compat=1.8}
\usetikzlibrary{patterns}
\usepackage{pgfplotstable, booktabs}  

\usepackage{graphicx,subfig}

\usepackage{tikz}       
\usetikzlibrary{hobby}  
\usepackage{multirow}
\usepackage{booktabs}

\usepackage{algpseudocode}
\usepackage{algorithm}
\restylefloat{algorithm}
\algrenewcommand\textproc{}
\usepackage{xpatch}

\algnewcommand{\algorithmicand}{\textbf{ and }}
\algnewcommand{\algorithmicor}{\textbf{ or }}
\algnewcommand{\OR}{\algorithmicor}
\algnewcommand{\AND}{\algorithmicand}

\algrenewcommand\alglinenumber[1]{{\sffamily\footnotesize#1}}

\makeatletter
\xpatchcmd{\algorithmic}{\itemsep\z@}{\itemsep=1ex plus2pt}{}{}
\makeatother

\newtheorem{theorem}{Theorem}
\newtheorem{definition}{Definition}
\newtheorem{lemma}{Lemma}

\newtheorem{proposition}{Proposition}
\newtheorem{remark}{Remark}
\newcommand{\specialformattingstuff}{a}
\newcommand{\specialformattingstufff}{b}
\newcommand{\specialformattingstuffff}{c}

\DeclareCaptionLabelFormat{cancaptionlabel}{#1 #2\specialformattingstuff}
\DeclareCaptionLabelFormat{cancaptionlabel2}{#1 #2\specialformattingstufff}
\DeclareCaptionLabelFormat{cancaptionlabel3}{#1 #2\specialformattingstuffff}

\begin{document}
\title {The inverse conductivity problem via the calculus of functions of bounded variation}
\author[1]{Antonios Charalambopoulos\thanks{acharala@math.ntua.gr}}
\author[1]{Vanessa Markaki\thanks{vamarkaki@central.ntua.gr}} 
\author[2]{Drosos Kourounis\thanks{drosos.kourounis@usi.ch}} 
\affil[1]{ Department of Applied Mathematical \& Physical Sciences, National Technical University Of Athens, 15780, Zografu, Greece}
\affil[2]{ Institute of Computational Science, Faculty of Informatics, Universit\`{a} della Svizzera italiana, 6904, Lugano, Switcherland}

\maketitle

\begin{abstract}
In this work, a novel approach for the solution of the inverse conductivity problem from one and multiple boundary measurements has been developed on the basis of the implication of the framework of $BV-$functions. The space of the functions of bounded variation is here recom\-mended as the most appropriate functional space hosting the conductivity profile under reconstruction. For the numerical solution, we propose and implement a suitable minimization scheme of an enriched - con\-structed herein - functional, by exploiting the inner structure of $BV-$ space. Finally, we validate and illustrate our theoretical results with numerical experiments.
\end{abstract}

MSC class: 35J25, 35R30,  35R05, 26B30, 65N21 (Primary), 35B27 (Secondary)
\date{}

\section{Introduction}
The problem of Electrical Impedance Tomography is of great importance from the theoretical and application point of view. It consists in the inverse problem aiming at the determination of the conductivity of an electrically conductive region when sufficient data are given on the surface of this region. The data of this inverse problem includes the knowledge of the voltage on the surface of the region -- which is usually incorporated in the physical assumptions of the problem -- along with the knowledge of the current density on the same surface (the role of these two surface fields may be interchanged), which is the outcome of a measurement process. The connection between surface voltage and current is linear and constitutes a well known and important operator, the so called voltage-to-current map, which is  of the type of Dirichlet-to-Neumann operators, encountered in elliptic boundary value problems of second order.  

One of the most interesting corner stone question has been whether the knowledge of the Dirichlet-to-Neumann map on the known surface $\partial \Omega$ of the conductive region is sufficient to determine the conductivity throughout the region $\Omega$. In case where the conductivity is isotropic, it was proposed by A. Calder\'{o}n \cite{calderon} in 1980 that any bounded conductivity might be determined solely from the boundary measurements, i.e. from the Dirichlet-to-Neumann operator. In \cite{astala} this has been confirmed for the two dimensional case. When the isotropic conductivity is smoother than just a $L^1$-function, the same conclusion is known to hold also in higher dimensions. 

The first global uniqueness result was obtained for a $C^\infty-$smooth conductivity in dimension $n \geq 3$ by J. Sylvester and G. Uhlmann in 1987 \cite{sylvester}. For the two dimensional case, A. Nachman \cite{nachman} produced in 1996 a uniqueness result for two times differentiable conductivities. The developed therein algorithm has been successfully implemented and proven to work efficiently even with real data \cite{silt}, \cite{mue}. The reduction of regularity assumptions has been since subject of active investigation. In two dimensions, the optimal $L^{\infty}-$regularity was obtained in \cite{astala}. The advantage of the reduction of the smoothness assumptions up to $L^{\infty}$ does not lie solely on the fact that many conductivities have jump-type singularities, but it also allows us to consider much more complicated singular structures such as porous media. 

The stability of the inversion is also very important and has been extensively investigated. The main argumentation can be encountered in  \cite{alessandrini},\cite{barcelo}, \cite{clop} and \cite{caro}, which are very representative references. The majority of the approaches establishing stability require some uniform control on the oscillations of the conductivity function and so deal with some kind of conditional stability. This is expected, since extreme oscillations of a sequence of conductivities create an instability of the Calder\'{o}n problem. This kind of asymptotics is well described via the implication of $H$ or $G-$convergence analysis, where the homogenization theory assigns the suitable convergence regime \cite{Ruiz}, \cite{jikov}, \cite{allaire}. 

As already mentioned, the investigation method depends strongly on the assumptions made for the regularity of the sought conductivity profile. Clearly, it is essential to assure the potentially discontinuous character of the conductivity function. This priority introduces the space of functions of bounded variation ($BV-$functions), as the most appropriate functional setting for the conductivity. In fact, modeling and formulation of a large number of problems in physics and technology require the introduction of new functional spaces, permitting discontinuities of the solution. In phase transitions, image segmentation, plasticity theory and the study of cracks and fissures, in the study of the wake in fluid dynamics and the shock theory in mechanics, the solution of the problem presents discontinuities along one-codimensional manifolds. Its first distributional derivatives are now measures, which may charge zero Lebesgue measure sets and its solutions cannot be considered as an element of Sobolev spaces throughout the entire domain of the problem. 

In the conductivity regime, the solution of the direct problem is an ordinary Sobolev function but the conductivity profile is not. Therefore, it seems gainful to select the space $BV$, as the appropriate hosting space for the conductivity functions. A function (with just $L^1-$ integrability behavior) belongs to $BV(\Omega)$ iff its first distributional derivatives are bounded measures. In particular, the space $SBV(\Omega)$ is a subspace of $BV(\Omega)$, which is more adequate for our purposes, since it contains functions whose distributional derivatives are free of the peculiar Cantor part \cite{butazzo}. Nevertheless, the supplementary, reasonable assumption that the conductivity profiles belong also to $L^{\infty}(\Omega)$ alleviates this concern and practically makes the profiles be retained inside $SBV(\Omega)$.

In the present work, the principal aim is the investigation of the inverse conductivity problem on the basis of a methodology developed inside the framework of $BV-$space. Our approach is motivated by the fundamental properties of the $BV-$functions, as these are well introduced in \cite{aubert} and \cite{butazzo}. Every $BV-$function (more accurately, every member of the subspace $SBV(\Omega)$) is a function in $L^1(\Omega)$, throughout the domain $\Omega$, which potentially disposes discontinuity surfaces. These surfaces consist the ``jump'' set, which represents the interfaces of the conductivity profile and is the support of the singular part (Young measures) of the distributional derivatives. Consequently, in case of reconstructing discontinuous conductivities, it is preferable to work within the $BV-$regime, since it is not necessary to define a priori a partition of the domain, every sub-domain of which supports one of the  continuous components of the conductivity function but just work with specific $BV-$functions, whose one of the main intrinsic characteristics is the possibly  multi-connected ``jump'' set. Therefore, a function in $BV$ represents simultaneously its discontinuity surfaces in the formation of its domain of definition. In \ref{ap1} we present, for completeness, the necessary properties of $BV-$functions. 

The implication of the regime of functions of bounded variation in the regularization of ill--posed problems has been already attempted and one of the primitive approaches can be encountered in \cite{vogel}. In the present work, a specific minimization scheme has been designed in order to solve the inverse conductivity problem, exploiting the advantages stemmed by the $BV$ calculus. The constructed functional deals with an ``energy'' over the whole domain $\Omega$, but when a $BV-$function is involved, a part of the energy is concentrated on the interface (``jump'' set) of this function \cite{butazzo}. The designed herein functional incorporates efficiently both the characteristics of the direct and inverse conductivity problem. Our motivation has been to detour - if necessary - solving directly the full conductivity problem at every step of the minimization process, but just perform the optimization descent, implementing simultaneously minimization on the conductivity profile and the $H^1(\Omega)$ potential fields participating in the direct Dirichlet and Neumann conductivity problems. The concept has been inspired by the works \cite{rich} and \cite{12332} and is based on the simple idea that the Dirichlet and Neumann solutions, corresponding to the surface data, should have a degree of compatibility (theoretically, they are identical given exact and noiseless data). Nevertheless, the aforementioned methodology is not restrictive, in the sense that it is always possible to solve the intermediate, consecutive, direct problems and transform the initial functional to a new one depending only on the unknown conductivity function. Then, this functional is very reminiscent (in case that a specific intrinsic weight parameter increases sufficiently) of the commonly used optimization functional, forcing the Dirichlet and Neumann surface data to compromise on the region boundary $\partial \Omega$.
Besides the alternative to detour the solution of a sequence of direct problems, the main motive to construct such an enriched functional is the anticipation that strong advantages will emerge when the method will be generalized in problems where the involving physical fields are also discontinuous functions, and thus should by themselves be represented as $BV-$functions (cracks,  fissures in mechanics etc). In that case, it might be preferable to develop a concurrent $BV-$minimization over profiles and fields.

When suitable, simultaneous convergence of physical parameters and fields occurs, it is useful to employ the regime of $G-$ and $H-$convergence, encountered in the framework of homogenization theory \cite{allaire}. Section \ref{sechom} is the introductory of this work and involves the fundamental properties of this theory, which are necessary to settle the forthcoming $BV-$analysis. It mainly contains two fruitful propositions controlling the convergence behavior of the intermediate problems (i.e. the members of the minimization sequence of the optimization scheme), which are greatly useful in all the stages of the forthcoming analysis. This is actually the reason why the main implementation of the method is preceded by Section \ref{sechom}. The proofs of these propositions are placed in \ref{ap2}, for authors' convenience. In Section \ref{icp}, the inverse conductivity problem is defined. Section \ref{func} includes the construction and the properties of the optimization functional. In Section \ref{tow}, we establish the necessary framework for the numerical implementation of the minimization scheme. Section \ref{numerics} involves the numerical investigation of the reconstruction method applied to some indicative, characteristic conductivity profiles.

\section{The regime of $H$ and $G-$ convergence} \label{sechom}
\subsection{The general framework}\label{subsec:generalframe}
In this section, we introduce the fundamental aspects of the homogenization theory, which in the next sections will be melt with the ingredients of the framework of functions of bounded variation. The notions of $G-$ convergence and $H-$ convergence are very essential, when involved limit processes are considered, in which we investigate the convergence of boundary value problems as the physical parameters alter drastically. More precisely, we consider the family of problems
\begin{align}
-\nabla \cdot (A^\epsilon(x) \nabla u_\epsilon (x)) &=f(x) \qquad &&\mbox{in} \  \Omega, \\
u_\epsilon(x)&=0 \qquad &&\mbox{on} \ \partial \Omega,
 \end{align} 
 with $A^\epsilon \in L^\infty(\Omega, \mathcal{M}_{b,c})$ standing for a sequence of matrix functions \cite{allaire}, where $\mathcal{M}_{b,c}$ represents the linear space $\mathcal{M}_N$ of the $N \times N$ matrices, which along with their inverses are bounded below (i.e. $M \in \mathcal{M}_{b,c}\ \mbox{iff}\ M \in \mathcal{M}_N \ \mbox{with} \ M \xi \cdot \xi \geq b |\xi|^2 \ \mbox{and} \ M^{-1} \xi \cdot \xi \geq c |\xi|^2,\ b,c > 0$). In this work we are interested in scalar isotropic problems, where the coefficients $A^\epsilon(x)$ are just proportional to the unit matrix: $A^\epsilon(x)=a^\epsilon(x) \mathbb{I}_{N \times N}$. The condition $a^\epsilon \in L^\infty(\Omega, [b,c])$ is now imposed, while the parameter $\epsilon$ goes to zero currying along the sequence of the corresponding boundary value problems. \\
The more general convergence notion is that of the $H-$ convergence
\begin{definition}
	A sequence of matrices $A^\epsilon(x) \in L^\infty(\Omega, \mathcal{M}_{b,c})$ is said to H-converge to the matrix $A^\ast(x) \in L^\infty(\Omega, \mathcal{M}_{b,c})$, as $\epsilon \rightarrow 0$, if for any right hand side $f \in H^{-1}(\Omega)$, the sequence of solutions of 
	\begin{align}
    -\nabla \cdot (A^\epsilon(x) \nabla u_\epsilon (x))&=f(x) \qquad &&\mbox{in} \  \Omega, \\
    u_\epsilon(x)&=0 \qquad &&\mbox{on} \ \partial \Omega,
	\end{align}   
	satisfies 
	\begin{align}
u_\epsilon&\underset{\mbox{weakly }}{{-\!\!\!\rightharpoonup}} u \ &&\mbox{in}\ H^{1}_{0}(\Omega), \\
A^\epsilon(x) \nabla u_\epsilon (x)&\underset{\mbox{weakly}}{{-\!\!\!\rightharpoonup}}  A^{\star}(x) \nabla u (x) \ &&\mbox{in}\ {(L^{2}(\Omega))}^N,\label{supplcond}
	\end{align}
	where $u$ is the solution of the homogenized boundary value problem
	\begin{align}
-\nabla \cdot (A^\star(x) \nabla u(x)) &=f(x) \qquad && \mbox{in} \  \Omega, \\
u(x)&=0 \qquad && \mbox{on} \ \partial \Omega.
	\end{align}
	\end{definition}
A fundamental result \cite{allaire} is that for any sequence of matrices $A^\epsilon(x)$ in $L^\infty(\Omega, \mathcal{M}_{b,c})$, there exists a subsequence, still denoted by $A^\epsilon(x)$, and an homogenized matrix $A^\star(x) \in L^\infty(\Omega, \mathcal{M}_{b,c})$ such that $A^\epsilon$ H-converges to $A^\star$. \\
In case that the matrices $A^{\epsilon}(x)$ are symmetric, a primitive notion of operator convergence had been introduced, the so called $G-$convergence. The only difference from the H-convergence is that the weak convergence (\ref{supplcond}) of the flux $A^\epsilon(x) \nabla u_\epsilon (x)$ is not required any more. Consequently, G-convergence is a weaker notion than H-convergence, in the sense that if a symmetric matrix H-converges to a symmetric homogenized matrix, then it automatically G-converges to the same limit. The converse is not obvious, but it turns out to be true \cite{allaire} and then, a sequence of symmetric matrices in $L^\infty(\Omega, \mathcal{M}_{b,c})$ G-converges to a limit iff H-converges to that limit. In order to relate these operator convergences with the usual ones, we mention that if a sequence of matrices  $A^\epsilon(x)$ in $L^\infty(\Omega, \mathcal{M}_{b,c})$ converges strongly to a limit matrix $A^\star(x)$ in $L^1(\Omega, \mathcal{M}_N)$ or converges to $A^\star(x)$ almost everywhere in $\Omega$, then $A^\epsilon(x)$ also H-converges to $A^\star(x)$. 

\subsection{Investigation of limiting processes of conductivity problems via $H$ and $G-$ convergence}\label{HGherein}
For the purposes of the present work, we reformulate in some extent the above, well established terminology by considering a countable family of boundary value problems, each of them identified with the natural number $n \in \mathbb{N}$. The parameter $\epsilon$ is discretized and actually could be assigned the values $\epsilon_n=\frac{1}{n}$ or independently the convergence would be considered in the sense $n \rightarrow \infty$. As already stated, we are interested in conductivity profiles of the type $A^{\epsilon_n}(x)=\alpha_n(x) \mathbb{I}_{N \times N}$. In particular, the case 
$\alpha_n\underset{L^1(\Omega)}{\longrightarrow} \alpha$ will emerge in Section \ref{func}, in combination with the supplementary convergence $\alpha_n \rightharpoonup \alpha \ \mbox{weakly} ^\ast \ \mbox{in} \ L^\infty(\Omega)$, while all the scalar conductivities $\alpha_n,\alpha$, belong to $L^\infty(\Omega,[b,c])$. As a result, we have that $\alpha_n(x) \mathbb{I}_{N \times N}\underset{H}{\longrightarrow} \alpha(x) \mathbb{I}_{N \times N}$. This means that $\forall f \in H^{-1}(\Omega)$, the sequence $\hat{u}_n^{\alpha_n}(x)$ of the solutions of the problems
\begin{align}
-\nabla \cdot (\alpha_n(x) \nabla \hat{u}_n^{\alpha_n}(x))&=h(x) \qquad && \mbox{in} \  \Omega,  \label{ouf1}\\
\hat{u}_n^{\alpha_n}(x)&=0 \qquad && \mbox{on} \ \partial \Omega, \label{ouf2}
\end{align}
obeys to the convergence rules:
\begin{align}
\hat{u}_n^{\alpha_n}&\underset{\mbox{{\it weakly }}}{{-\!\!\!\rightharpoonup}} u^{\alpha} \ &&\mbox{in}\  H^{1}_{0}(\Omega), \\
\alpha_n(x) \nabla \hat{u}_n^{\alpha_n}(x)&\underset{\mbox{{\it weakly }}}{{-\!\!\!\rightharpoonup}}\alpha(x) \nabla u^{\alpha}(x)  \ &&\mbox{in}\  {(L^{2}(\Omega))}^N \label{supplcond2}
\end{align}
as $n \rightarrow \infty$, where $u^\alpha$ is the unique solution of the problem
\begin{align}
-\nabla \cdot (\alpha(x) \nabla u^{\alpha}(x))&=h(x) \qquad &&\mbox{in} \  \Omega, \label{oriako1} \\
 u^{\alpha}(x)&=0 \qquad &&\mbox{on} \ \partial \Omega. \label{oriako2}
\end{align}

The following two kinds of problems will emerge in Section \ref{func}. \\
\underline{Problem I}: 
\begin{align}
-\nabla \cdot (\alpha_n(x) \nabla u_n^{\alpha_n}(x))&=h_n(x) \qquad&& \mbox{in} \  \Omega, \label{pr1}\\
u_n^{\alpha_n}(x)&=0 \qquad&& \mbox{on} \ \partial \Omega.
\end{align} 
\underline{Problem II}:
\begin{align}
\nabla \cdot (\alpha_n(x) \nabla w_n^{\alpha_n}(x))&=0 \qquad&& \mbox{in} \  \Omega, \\
\alpha_n(x) \frac{\partial w_n^{\alpha_n}}{\partial n}(x)&=g(x) \qquad&& \mbox{on} \ \partial \Omega,
\end{align}
where $h_n \in H^{-1}(\Omega)$ and $g \in H^{-\frac{1}{2}}(\partial \Omega)$. \\
The first auxiliary lemma, whose proof is reminiscent of the arguments presented in Proposition 1.2.19 of \cite{allaire} - and placed in \ref{ap2} - establishes the convergence of the fluxes given the convergence of the fields.   
\begin{lemma} \label{prop1}
	In case that the solution of problem I obeys to the convergence $u_n^{\alpha_n} \underset{n \rightarrow \infty}{-\!\!\!\rightharpoonup} v$ ( weakly \ in $H^1_0(\Omega)$) and the sequence $h_n$ converges strongly in $H^{-1}(\Omega)$, then $\alpha_n(x) \nabla u_n^{\alpha_n} \underset{n \rightarrow \infty}{-\!\!\!\rightharpoonup} \alpha \nabla v$ (weakly \ in ${\left( L^2(\Omega)\right)}^N$).
\end{lemma}
The second auxiliary lemma deals with the energy convergence of the sequence of the problems of type I. 
\begin{lemma} \label{prop2}
	In case that the stimulus $h_n$ of the Problem I obeys to the convergence $h_n \underset{n \rightarrow \infty}{-\!\!\!\rightarrow} h$ (strongly in $H^{-1}(\Omega)$), then the sequence $u_n^{\alpha_n}$ obeys to the energy convergence
	\begin{align}
	\int_{\Omega}\alpha_n|\nabla u_n^{\alpha_n}{|}^2 \longrightarrow \int_{\Omega}\alpha |\nabla u^{\alpha}{|}^2. \end{align}
\end{lemma}
{\bf Proof.}
See also in \ref{ap2}.

We introduce here, the  trace operator $\gamma : H^{1}(\Omega) \rightarrow H^{\frac{1}{2}}(\partial \Omega)$, which is one to one, onto and disposes as right inverse the extension operator $\eta : H^{\frac{1}{2}}(\partial \Omega) \rightarrow H^{1}( \Omega)$, which extends surface functions into $H^1$ - functions defined in $\Omega$. This operator will be useful in the sequel and just mention here, that when we state, as example, that  $u_n^{\alpha_n}(x)=0 \ \mbox{on} \ \partial \Omega$, we mean that $\gamma u_n^{\alpha_n}=0$. \\
The interesting class of stimuli $h_n$ will be of the form
\begin{align}
h_n(x)= \nabla \cdot(\alpha_n(x) \nabla (\eta f)(x) ) \in H^{-1}(\Omega), \label{glykos}
\end{align}
$\mbox{with} \ f \in H^{\frac{1}{2}}(\partial \Omega)$.\newline
Next, we consider the Problem I furnished with this particular type of non-homogeneous terms (\ref{glykos}). In this framework, we state (and prove in \ref{ap2}) the following two Propositions.   

\begin{proposition} \label{prop3}
	The sequence $u_n^{\alpha_n},\ n \in \mathbb{N}$ converges to $u^\alpha$ weakly in $H^1_0(\Omega)$ as $n \rightarrow \infty$ and furthermore $\alpha_n(x) \nabla {u}_n^{\alpha_n}(x)  \rightharpoonup  \alpha(x) \nabla u^\alpha (x)\ \mbox{weakly \ in}\ {(L^{2}(\Omega))}^N$, where $u^{\alpha}$ is the unique solution of the problem (\ref{oriako1})-(\ref{oriako2}). In addition, it holds that 
	\begin{align}
	\int_{\Omega}\alpha_n|\nabla (u_n^{\alpha_n}+\eta f)|^2 \longrightarrow \int_{\Omega}\alpha |\nabla (u^{\alpha}+\eta f)|^2.
	\end{align} 
\end{proposition} 

\begin{proposition} \label{prop4}
	The sequence of solutions of the Problem II converges to $w^\alpha$ weakly in $H^1_0(\Omega)$ as $n \rightarrow \infty$, where $w^\alpha$ is the solution of the problem 
\begin{align}
-\nabla \cdot (\alpha(x) \nabla w^{\alpha}(x))&=0 \qquad&& \mbox{in} \  \Omega, \label{firstfinal} \\
\alpha(x) \frac{\partial w^{\alpha}}{\partial n}(x)&=g(x) \qquad&& \mbox{on} \ \partial \Omega.  \label{secondfinal}
\end{align}
In addition 
\begin{align}
\int_{\Omega}\alpha_n|\nabla w_n^{\alpha_n}|^2 \longrightarrow \int_{\Omega}\alpha |\nabla w^{\alpha}|^2.
\end{align} 
\end{proposition}

\section{The inverse conductivity problem} \label{icp}
Let $\Omega \subset \mathbb{R}^N$, $N=2,3$, denotes a conductive medium, whose conductivity $\alpha(x),\ x\in \Omega$, is the target of reconstruction. Usually, there exists a probably disconnected object $D$, which is an open subset of $\Omega$, being the target of detection and diversified from the remaining conductive material (i.e. the matrix) via the possibly discontinuous change of the conductivity parameter on the interface $\partial D$. 
\newline
The potential $u$ inside the structure $\Omega$ satisfies the direct boundary value problem 
\begin{align}
\nabla \cdot (\alpha(x) \nabla u(x))&=0 \ \ && x \in \Omega, \\
\alpha(x) \frac{\partial u}{\partial n} (x)&=g(x) \ \ && x \in \partial \Omega,
\end{align}
with the coefficient of conductivity having the form
\begin{align}
\alpha(x) = \widetilde{b}(x) \chi_{\Omega \setminus D}(x) + \widetilde{c}(x) \chi_{D}(x)
\end{align}
where we meet the characteristic function $\chi_A$ of a set $A$. The functions $\widetilde{b}(x), \widetilde{c}(x)$ characterize the background (matrix) and the inclusion region respectively. In several applications, these functions are constant and constitute representations of two-phase materials. Here, we allow them to be variable functions with discontinuous behavior on the interfaces of the inclusions. The only quantitative restriction is that they take values in the interval $[b,c]$, where $b, c$ are real thres\-hold values with $0 < b < c < \infty$. So, it holds that $\alpha(x)\in L^{\infty}(\Omega,[b,c])$. The additional assumption in this work is that $\alpha(x)\in BV(\Omega)$. We mention here the simple pilot case of a two - valued conductivity profile. In that case, the inverse conductivity problem consists in the determination of the measurable function $\alpha(x) \in L^\infty(\Omega,\{b,c\}) \cap BV(\Omega)$. The boundary condition of the problem imposes a specific current $g \in H^{-\frac{1}{2}}(\partial \Omega)$ on the whole boundary $\partial \Omega$ of the body. It is well known, that such a boundary value problem is solvable only when the surface data obeys to compatibility condition $\left\langle g,1 \right\rangle_{H^{-\frac{1}{2}} (\partial \Omega) \times H^{\frac{1}{2}}(\partial \Omega)}=0$, and there emerges an equivalence class of solutions $u \in  H^{1}(\Omega)$, whose members differ each other by a real constant. To unify the members of the equivalence classes above, we introduce the quotient space $H^{1}(\Omega)/\mathbb{R}$ with the appropriate norm $\|u\|_{H^{1}(\Omega)/\mathbb{R}}=\inf_{\lambda \in \mathbb{R}}\|u+\lambda\|_{H^{1}(\Omega)}$. Notice that every solution $u$, along with the electric current flux $\alpha \frac{\partial u}{\partial n}$, is continuous across the interface $\partial D$.  

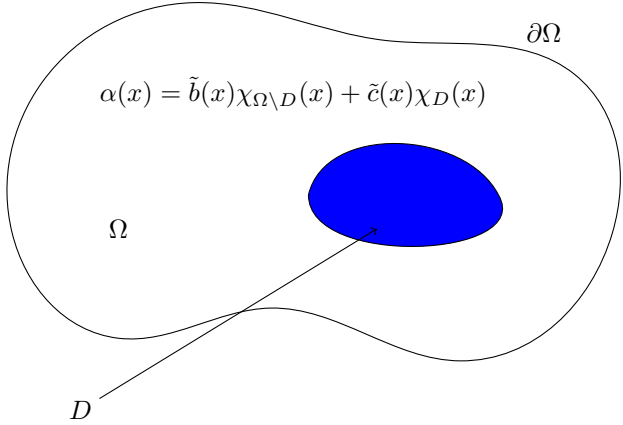
\begin{figure}[H] 
	\centering
	 \input{figure1smaller.tex}
	\caption{The inverse conductivity problem consists in the determination of $\alpha(x)$.}
	\label{fig1}
\end{figure} 

 The inverse conductivity problem consists in the determination of $\alpha(x)$, in case that we have at hand as supplementary data, in the realm of measurements, the voltage $u|_{\partial \Omega}$ on the boundary of the body. This is equivalent to say that we have the knowledge of the Neumann to Dirichlet map\footnote{or the current to voltage map} denoted by $\Lambda_{N.t.D}$ $:H^{-\frac{1}{2}}(\partial \Omega) \rightarrow H^{\frac{1}{2}}(\partial \Omega)$ or equivalently its inverse $\Lambda_{D.t.N}$ $:H^{\frac{1}{2}}(\partial \Omega) \rightarrow H^{-\frac{1}{2}}(\partial \Omega)$, fact corresponding to the situation, in which we initially have knowledge about the Dirichlet values of the field and we acquire a posteriori additional information about the current field on the surface $\partial \Omega$. We will next denote $\Lambda_{\alpha}=\Lambda_{D.t.N}$, in order to make clear the dependence of this surface operator on the conductivity parameter function $\alpha(x)$ and then, we have immediately that $\Lambda_{\alpha}^{-1}=\Lambda_{N.t.D}$.\newline
Summarizing, the following inverse problem is formulated
\begin{align}\label{eq:PDE}
\nabla \cdot (\alpha(x) \nabla u(x))&=0 \ \ && x \in \Omega, \\
u(x)&=f(x) \ \ && x \in \partial \Omega, \\
\alpha(x) \frac{\partial u}{\partial n} (x)&=g(x) \ \ && x \in \partial \Omega,
\end{align}   
with known data $f \in H^{\frac{1}{2}}(\partial \Omega)$ and $g=\Lambda_\alpha(f) \in H^{-\frac{1}{2}}(\partial \Omega)$. \\ A thorough analysis, incorporating the necessary compatibility condition of the Neumann data, leads to adopting the more appropriate functional setting $\Lambda_{\alpha}: \mathcal{B}\left(:=H^{\frac{1}{2}}(\partial \Omega)/\mathbb{R}\right) \rightarrow {\left(H^{\frac{1}{2}}(\partial \Omega)/\mathbb{R}\right)}^{\ast}$, where an isometric isomorphism of the dual space ${\mathcal{B}}^{\ast}={\left(H^{\frac{1}{2}}(\partial \Omega)/\mathbb{R}\right)}^{\ast}$ is the space of distributions $\mathcal{A}:=\left\{g \in H^{-\frac{1}{2}} (\partial \Omega) \right. $ $ \left. : \left\langle g,1 \right\rangle_{H^{-\frac{1}{2}} (\partial \Omega) \times H^{\frac{1}{2}}(\partial \Omega)}=0\right\}$, in which the data $g$ should belong. 
\section{On the investigation of the appropriate minimization functional for the inverse conductivity problem} \label{func}
Building the appropriate functional under minimization is based on connecting efficiently the conductivity function $\alpha$, along with the solution of the Dirichlet and the Neumann boundary value problem. More precisely, the functional under examination should involve the term
\begin{align}
J_1(\alpha,w)=\frac{1}{2}\int_{\Omega}\alpha(x) |\nabla w(x)|^2 d x-\left\langle g,\gamma w \right\rangle_{\mathcal{A} \times \mathcal{B}}. \label{first}
\end{align}
The minimization of this term  - even under the assumption of a specific function $\alpha$ - over all possible $w \in H^{1}( \Omega)/\mathbb{R}$, would lead to the solution of the Neumann boundary value problem $w^{\alpha}$, corresponding to the conductivity coefficient $\alpha(x)$. However, this is not the ultimate settlement, since the function $\alpha$ is not known at all - as a matter of fact it is the target of our investigation - and more data should be incorporated. 
The functional under construction should also contain the term
\begin{align}
J_2(\alpha,u)=\frac{1}{2}\int_{\Omega}\alpha(x) |\nabla (u(x)+(\eta f)(x))|^2 d x. \label{second}
\end{align}
Here we meet again the continuous right inverse $\eta$ of the trace operator $\gamma$, encountered already in Section \ref{sechom}. The term (\ref{second}) would independently be minimized - over all functions $u \in H_0^1(\Omega)$  - by the unique solution $u^\alpha(x)+(\eta f)(x)$ of the Dirichlet problem, in case that the conductivity function was considered as a fixed function. \newline The third term of the functional should be a forcing term, imposing a theoretical vanishing of the difference of the solutions $u^{\alpha}+\eta f$ and $w^{\alpha}$ in $\Omega$, in case that the data respected exactly the relation $g=\Lambda_{\alpha}f$ - without measurement errors - and the function $\alpha$ was indeed the ``correct" conductivity coefficient over $\Omega$. This term should then incorporate the integral $\frac{1}{2}\int_{\Omega}\alpha |\nabla (u+\eta f-w)|^2 d x$ forcing the fields to coincide \cite{rich}. Actually, we force, this way, the gradient of the fields to compensate each other. We avoid to add a term forcing the $L^2-$norm (weighted by $\alpha(x)$) of the difference $u+\eta f-w$ to diminish, since this $L^2-$distance does not obey to the good $H-$convergence results emerging in the minimization process. Furthermore, in order to manipulate this minimization descent to be monitored, so as to converge to the profiles and fields of the target conductivity problem, we need to incorporate, in the construction of $J_3$, the term $\int_{\Omega}\alpha \nabla (u+\eta f-w) \cdot \nabla w d x + \left\langle g,\gamma w-f \right\rangle_{\mathcal{A}\times \mathcal{B}}$. The specific manner this term leads the minimization tracing will be clarified in the proof of Theorem \ref{thfund}. Selecting a regularization parameter $\kappa$, we schedule then the penalty term $J_3$ as follows 
\begin{align}
 J_3(\alpha,u,w)=&\frac{\kappa}{2}\int_{\Omega}\alpha |\nabla (u+\eta f-w)|^2 d x \nonumber \\ &+ \kappa \int_{\Omega}\alpha \nabla (u+\eta f-w) \cdot \nabla w d x \nonumber \\ &+ \kappa \left\langle g,\gamma w-f \right\rangle_{\mathcal{A} \times \mathcal{B}}. \label{third}
\end{align}   
 It is clear that all terms of (\ref{third}) are identically zero in case of exact conductivity and data.\

The fourth term of the functional is a term properly referring to the selection over conductivity profiles obeying to restriction rules concerning their total masses and some reference compatibility. This last term of the functional involves two regularization parameters $\lambda, \mu$ and obtains the following form
\begin{align}
J_4=\frac{\lambda}{2}\|R\alpha-\alpha_0\|^2_{L^2(\Omega)}+\mu \int_{\Omega} \phi(|D \alpha|), \label{fourth}      
\end{align}
where $\phi$ is a convex function of its argument, $R$ a bounded operator $R: L^2(\Omega) \rightarrow L^2(\Omega)$ and $\alpha_0$ a reference value for the conductivity coefficient. Adding the previous terms (\ref{first})-(\ref{fourth}), we construct the total functional
\begin{align}
 E(u,w,\alpha)=\tilde{E}(u,w,\alpha)&+\frac{\lambda}{2}\|R\alpha-\alpha_0\|^2_{L^2(\Omega)}\nonumber\\ &+\mu \int_{\Omega} \phi(|D \alpha|), \label{total}
\end{align}  
where
\begin{align}
\tilde{E}(u,w,\alpha)=&\frac{\kappa+1}{2}\int_{\Omega}\alpha|\nabla (u+\eta f)|^2 d x \nonumber\\ & +(1-\kappa)\left[\frac{1}{2}\int_{\Omega}\alpha |\nabla w|^2 d x-\left\langle g,\gamma w \right\rangle_{\mathcal{A} \times \mathcal{B}}\right]\nonumber \\ &-\kappa\left\langle g,f \right\rangle_{\mathcal{A} \times \mathcal{B}} \label{total1}
\end{align}
with the parameter $\kappa$ being restricted in the interval $(-1,1)$, in order to preserve the coercivity of the first two terms of (\ref{total1}).\newline We will clarify now the admissible set for the conductivity functions $\alpha$, and as a byproduct, the special role of the operator $R$ is going to be revealed. As explained previously, the abrupt discontinuities in the physical properties of the domain $\Omega$ are very well represented via the functions of bounded variation. One additional, essential property of the $\alpha$ coefficient  stems from the requirement that there exists a very thin region away the surface $\partial \Omega$, in which the conductivity does not change its value. We introduce the space $\Omega_{\delta}=\left\{x \in \Omega : \mbox{dist}(x, \partial \Omega) < \delta \right\}$ and demand that $\alpha \in BV(\Omega) \cap L_\delta^\infty(\Omega,\left\{b,c \right\})$, in case of the pilot two constant valued profile, where $L_\delta^\infty(\Omega,\left\{b,c \right\})=\left\{ \beta \in L^\infty(\Omega,\left\{b,c \right\}) : \|\beta-b\|_{L^\infty(\Omega_\delta)}=0 \right\}$. Similarly, in case of a variable profile, we state that $\alpha\in BV(\Omega) \cap L_\delta^\infty(\Omega,[b,c])$, where now the values of the conductivity range inside the interval $[b,c]$.
We mention here, that the necessity of the implication of the parameter $\delta$ and the zone $\Omega_\delta$ is due to the desire to have the ability to take strong limits of the Calder\'{o}n's operator $\Lambda_{\alpha}$, when $\alpha$ follows specific limit processes \cite{Ruiz}. \\
As far as the reference operator $R$ is concerned, we could select as $R \alpha$ the restriction $\alpha\mid_{\Omega_\delta}$ and in consequence $\|R\alpha-\alpha_0\|_{L^2(\Omega)}$ is defined to be the norm $\|\alpha-b\|_{L^2(\Omega_\delta)}$, expressing the $L^2(\Omega_\delta)$-difference of the conductivity coefficient from the background level value $b$. This is a quantitative manner to force the conductivity to be close to the value $b$ in the vicinity of the boundary, throughout all the steps of the minimization process. Alternatively, the requirement $\alpha\mid_{\Omega_\delta}=b$ could be exactly imposed on the members of the minimization conductivity scheme. In any case, there exist certain realizations of $\|R\alpha-\alpha_0\|_{L^2(\Omega)}$, exploiting the function $\alpha_0$ as an initial guess for the unknown conductivity in the whole or some part of $\Omega$.  \\      
The fundamental result of this section is the following theorem
\begin{theorem} \label{thfund}
The minimization problem $\inf_{(u,w,\alpha^\prime)} E(u,w,\alpha^\prime)$ over elements $(u,w,\alpha^\prime)$ in $H_0^1(\Omega) \times H^{1}( \Omega)/\mathbb{R} \times BV(\Omega) \cap L_\delta^\infty(\Omega,\left\{ b,c \right\})$ admits a solution belonging to the space $H_0^1(\Omega) \times H^{1}( \Omega)/\mathbb{R} \times BV(\Omega) \cap L_\delta^\infty(\Omega,\left[ b,c \right])$. 
\end{theorem}
{\bf Proof.}
	The proof will be accomplished in several steps. Let us consider a minimizing sequence $(u_n,w_n,\alpha_n)$, $n \in \mathbb{N}$ such that the sequence $E(u_n,w_n,\alpha_n)$ converges to the finite real number $\inf_{(u,w,\alpha^\prime)} E(u,w,\alpha^\prime)$. Clearly, there exists a positive constant $C$ such that 
	\begin{align}
	E(u_n,w_n,\alpha_n) \leq C \ \mbox{for \ all} \ n \in \mathbb{N}. \label{v1}
	\end{align}
	On the other hand, given that $\alpha_n \in L_\delta^\infty(\Omega,\left\{ b,c \right\})$, we deduce that $\|\alpha_n\|_{L^\infty(\Omega)} \leq c$ and so there exists a subsequence, still denoted by $\alpha_n$, and an element $\alpha \in L^\infty(\Omega)$ such that 
	\begin{align}
	\alpha_n \rightharpoonup \alpha \ \mbox{weakly} ^\ast \ \mbox{in} \ L^\infty(\Omega) \ (\mbox{as} \ n \rightarrow \infty). \label{firstone}
	\end{align}
	Thanks to a very interesting result \cite{tartar}, the weak$^\ast$ closure of the space $L^\infty(\Omega, K)$, where $K$ is a subspace of $\mathbb{R}^p$, is $L^\infty(\Omega, \mathcal{K})$, where $\mathcal{K}$ is the closed convex hull of $K$. In our case, we deduce that $\alpha \in L^\infty(\Omega,\left[b,c \right])$. Furthermore, it is obvious that $\int_{{\Omega}_{\delta}} \alpha_n \varphi \rightarrow \int_{{\Omega}_{\delta}} \alpha \varphi$ for every $\varphi \in C_0^\infty(\Omega_\delta) \subset L^1(\Omega)$ and then $\int_{{\Omega}_{\delta}} (b-\alpha(x)) \varphi(x) dx = 0, \ \forall \varphi \in C_0^\infty(\Omega_\delta)$. Consequently, we deduce that $\alpha(x)=b$, almost everywhere in $\Omega_{\delta}$ \cite{brezis} and so finally $\alpha \in L_\delta^\infty(\Omega,\left[b,c \right])$.   
	Returning to (\ref{v1}), we start taking advantage from this uniform boundedness by claiming that there exists a uniform bound $M$ (independent of $n$) such that
	\begin{align}
	\mathcal{C}_{n}:=\left\langle g,\gamma w_n \right\rangle_{\mathcal{A} \times \mathcal{B}}-\frac{1}{2}\int_{\Omega}\alpha_n |\nabla w_n|^2 d x \leq M. \label{wproof}
	\end{align}
	Indeed, due to the boundedness of $g$ and the trace operator $\gamma$, we have that 
	\begin{align}
	&\left\langle g,\gamma w_n \right\rangle_{\mathcal{A} \times \mathcal{B}}-\frac{1}{2}\int_{\Omega}\alpha_n |\nabla w_n|^2 d x \nonumber\\ &\leq C^\prime \|g\|\|w_n\|_{H^{1}( \Omega)/\mathbb{R}}-\frac{b}{2}\|\nabla w_n\|^2_{L^2(\Omega)}  \nonumber \\
	&=C^\prime \|g\| \inf_{\nu \in \mathbb{R}} \|w_n+\nu\|_{H^{1}(\Omega)}-\frac{b}{2}\|\nabla w_n\|^2_{L^2(\Omega)}\nonumber \\
	& \leq C^\prime \|g\|  \|w_n -\frac{1}{|\Omega|} \int_\Omega w_n(x)dx \|_{H^{1}(\Omega)}-\frac{b}{2}\|\nabla w_n\|^2_{L^2(\Omega)} \nonumber \\
	& \leq C^{\prime \prime} \|g\| \|\nabla w_n\|_{L^2(\Omega)} - \frac{b}{2}\|\nabla w_n\|^2_{L^2(\Omega)}, \label{help}
	\end{align}
	where we have used the Poincar$\acute{e}$ - Wirtinger inequality for the field $w_n-\frac{1}{|\Omega|} \int_\Omega w_n(x)dx$, with zero mean value over the finite region $\Omega$. The last right term of the inequality (\ref{help}) is necessarily bounded, since it is the negative of a coercive function. So (\ref{wproof}) has been proved. \\
This equation, along with the special form of $E(u_n,w_n,\alpha_n)$, lead easily to the uniform boundedness of the following terms: 
	\begin{align}
	&\int_{\Omega} \phi(|D \alpha_n|), \ \|R \alpha_n-\alpha_0\|_{L^2(\Omega)}, \nonumber \\ & \mbox{and} \ \int_{\Omega}\alpha_n|\nabla (u_n+\eta f)|^2 d x \ < \ \infty. \nonumber\label{sbounds}
	\end{align} 
Using the convexity of $\phi$, we further prove that $|D \alpha_n|(\Omega) \leq C$. Due to the aforementioned coercivity of the terms $(-\mathcal{C}_{n})$, it is necessary that the integrals $\int_{\Omega}\alpha_n |\nabla w_n|^2 d x$ are also uniformly bounded. We deduce from these results and the fact that $\alpha_n$ is bounded below by the value $b$, that the norms $\|\nabla(u_n+\eta f)\|_{L^2(\Omega)}$ and  $\|\nabla w_n\|_{L^2(\Omega)}$ are again uniformly bounded. This is inherited of course to the norms $\|\nabla u_n\|_{L^2(\Omega)}$ and using the Poincar$\acute{e}$ inequality for the functions $u_n \in H_0^1(\Omega)$, we infer that 
	\begin{align}
	\|u_n\|_{H_0^1(\Omega)} \leq C, \ \ n \in \mathbb{N}. 
	\end{align} 
	In addition, thanks again to the Poincar$\acute{e}$ - Wirtinger inequality, we find that 
\begin{align}
\inf_{\nu \in \mathbb{R}} \|w_n+\nu\|_{H^{1}(\Omega)}\nonumber \\ & \leq \|w_n -\frac{1}{|\Omega|} \int_\Omega w_n(x)dx \|_{H^{1}(\Omega)} \nonumber \\ &\leq  \|\nabla w_n\|_{L^2(\Omega)} \leq C, \ \ n \in \mathbb{N}. 
	\end{align} 
	We keep from the last two bounds that there exist $u \in H_0^1(\Omega), w \in H^{1}(\Omega)/\mathbb{R}$, such that 
	\begin{align}
	& u_n \underset{n \rightarrow \infty}{-\!\!\!\rightharpoonup} u \ \in H_0^1(\Omega), \label{guer1} \\
	 &\nabla w_n  \underset{n \rightarrow \infty}{-\!\!\!\rightharpoonup} \nabla w \ \in {(L^2(\Omega))}^N. \label{guer2}
	\end{align}
	Moreover, $\alpha_n \in L_\delta^\infty(\Omega,\left\{ b,c \right\})$ and so $\alpha_n \in L^1(\Omega)$ with $\|\alpha_n\|_{L^1(\Omega)} \leq c |\Omega|$. We infer that the sequence $\alpha_n$ is uniformly bounded in $BV(\Omega)$ and so thanks to the compactness theorem (2.2.3 of \cite{aubert}), we extract a subsequence, with the same symbolism, and an element $\tilde{\alpha} \in BV(\Omega)$ such that 
	\begin{align}
	&\alpha_n \underset{BV-w^\ast}{-\!\!\!\rightharpoonup} \tilde{\alpha}\ \ \ (\mbox{i.e.} \ a_n \underset{L^1(\Omega)}{\longrightarrow} \tilde{\alpha} \ \mbox{and}\nonumber \\ &\int_\Omega \varphi D\alpha_n \rightarrow \int_\Omega \varphi D\tilde{\alpha} \ \ \forall \varphi \in {(C_0^1(\Omega))}^N).
	\end{align} 
	A fortiriori, we obtain clearly that $\int_\Omega \alpha_n \varphi \rightarrow \int_\Omega \tilde{\alpha} \varphi$, $\forall \varphi \in C_0^\infty(\Omega)$. On the other hand, on the basis of (\ref{firstone}), we obtain $\int_\Omega \alpha_n \varphi \rightarrow \int_\Omega {\alpha} \varphi$, $\forall \varphi \in C_0^\infty(\Omega) \subset L^1(\Omega)$. These two last results imply that $\int_\Omega (\tilde{\alpha}-\alpha) \varphi=0$, $\forall \varphi \in C_0^\infty(\Omega)$ and so we finally obtain that $\tilde{\alpha}=\alpha, \ \mbox{a.e} \ \in \Omega$ \cite{brezis}. So, the convergence of the sequence of the conductivities $\alpha_n$ obeys to the rule
	\begin{align}
	\alpha_n  \underset{L^1(\Omega)}{\longrightarrow} {\alpha} \in BV(\Omega) \cap L_\delta^\infty(\Omega,\left[b,c \right]),
	\end{align}
	plus the convergence (\ref{firstone}). Furthermore, given that the sequence $\alpha_n \in L_\delta^\infty(\Omega,\left\{b,c \right\})$ satisfies $\|\alpha_n\|_{L^2(\Omega)} \leq c^2 |\Omega|$, the sequence $\alpha_n$ - actually a subsequence of it with which we keep on working - converges weakly \ in $L^2(\Omega)$ and so 
	\begin{align}
	R \alpha_n \underset{L^2(\Omega)}{-\!\!\!\rightharpoonup} R \alpha. 
	\end{align}  
	The ultimate goal of the theorem is to take the limit in (\ref{v1}) as $n \rightarrow \infty$ and prove some kind of lower semi-continuity for the functional, in order to obtain the same inequality for the limiting triple $(u,w,\alpha)$. First, we remark that due to the convexity of $\phi$ and the lower semi-continuity of the $BV-w^\ast$ convergence (see Appendix A), it holds that 
	\begin{align}
	\mu \int_{\Omega} \phi(|D \alpha|) \leq \mu \underset{n \rightarrow \infty}{\underline{\lim}}\int_{\Omega} \phi(|D \alpha_n|). \label{paok1}
	\end{align}   
	In addition, from the lower semi-continuity of the $L^2-$norm in the weak topology, we obtain
	\begin{align}
	\frac{\lambda}{2} \|R\alpha-\alpha_0\|^2_{L^2(\Omega)}\leq \frac{\lambda}{2} \underset{n \rightarrow \infty}{\underline{\lim}}\|R\alpha_n-\alpha_0\|^2_{L^2(\Omega)}.  \label{paok2}
	\end{align}
	To step further, we consider the quadratic and the linear terms of the portion $\tilde{E}(u_n,w_n,\alpha_n)$ of the functional. We insert in the scene the fields $u_n^{\alpha_n}$ and $w_n^{\alpha_n}$, encountered in Section \ref{sechom} (Problem I, Problem II), which solve the exact Dirichlet and Neumann problems correspondingly, in case that the conductivity coefficient coincides with $\alpha_n(x)$. Then, we write in condensed form
\begin{align}
&\frac{\kappa+1}{2}\int_{\Omega}\alpha_n|\nabla (u_n+\eta f)|^2=\nonumber \\ &\frac{\kappa+1}{2}\int_{\Omega}\alpha_n|\nabla (u_n^{\alpha_n}+\eta f)|^2 + \frac{\kappa+1}{2}\int_{\Omega}\alpha_n|\nabla (u_n-u_n^{\alpha_n})|^2 \nonumber \\
& +(\kappa+1)\int_{\Omega}\alpha_n \nabla (u_n^{\alpha_n}+\eta f) \cdot  \nabla (u_n-u_n^{\alpha_n}) \nonumber \\
& =\frac{\kappa+1}{2}\int_{\Omega}\alpha_n|\nabla (u_n^{\alpha_n}+\eta f)|^2 +\frac{\kappa+1}{2}\int_{\Omega}\alpha_n|\nabla (u_n-u_n^{\alpha_n})|^2,\label{dec11}
\end{align}
where we have used that $u_n-u_n^{\alpha_n} \in H^1_0(\Omega)$ and the generalized Green's formula to prove that $$(\kappa+1)\int_{\Omega}\alpha_n \nabla (u_n^{\alpha_n}+\eta f) \cdot\nabla (u_n-u_n^{\alpha_n})=0.$$
In addition
	\begin{align}
	&\frac{1-\kappa}{2}\int_{\Omega}\alpha_n|\nabla w_n|^2= \nonumber\\ &\frac{1-\kappa}{2}\int_{\Omega}\alpha_n|\nabla w_n^{\alpha_n}|^2 + \frac{1-\kappa}{2}\int_{\Omega}\alpha_n|\nabla (w_n-w_n^{\alpha_n})|^2 \nonumber \\
	& +(1-\kappa)\int_{\Omega}\alpha_n \nabla w_n^{\alpha_n} \cdot  \nabla (w_n-w_n^{\alpha_n})\nonumber\\
	& =\frac{1-\kappa}{2}\int_{\Omega}\alpha_n|\nabla w_n^{\alpha_n}|^2+\frac{1-\kappa}{2}\int_{\Omega}\alpha_n|\nabla (w_n-w_n^{\alpha_n})|^2 \nonumber \\
	& +(1-\kappa)\left\langle g,\gamma w_n-\Lambda_{\alpha_n}^{-1} g \right\rangle_{\mathcal{A} \times \mathcal{B}},\label{dec2}
	\end{align} 
	where we have used again the Green's formula and the Calder\'{o}n operator to handle the surface term. Finally, the linear term of $\tilde{E}$ becomes
	\begin{align}
	-(1-\kappa)\left\langle g,\gamma w_n \right\rangle_{\mathcal{A} \times \mathcal{B}}&=\nonumber \\ &-(1-\kappa)\left\langle g,\gamma w_n^{\alpha_n} \right\rangle_{\mathcal{A} \times \mathcal{B}}\nonumber \\ &-(1-\kappa)\left\langle g, \gamma w_n- \Lambda_{\alpha_n}^{-1} g \right\rangle_{\mathcal{A} \times \mathcal{B}}.\label{dec3}
	\end{align} 
	The decomposition made in (\ref{dec11})-(\ref{dec3}) allows to write 
	\begin{align}
	\tilde{E}(u_n,w_n,\alpha_n)=&\tilde{E}(u_n^{\alpha_n},w_n^{\alpha_n},\alpha_n)\nonumber \\ &+\frac{\kappa+1}{2}\int_{\Omega}\alpha_n|\nabla (u_n-u_n^{\alpha_n})|^2 \nonumber \\
	& +\frac{1-\kappa}{2}\int_{\Omega}\alpha_n|\nabla (w_n-w_n^{\alpha_n})|^2. \label{bf}
	\end{align}
	The form of (\ref{bf}) is very helpful for monitoring the minimization descent. The structure of the sub-functional $J_3$ has been responsible for intermediate eliminations of linear terms. That's why it has been selected in this particular form. The outcome of these cancellations consists in expressing the difference of $\tilde{E}(u_n,w_n,\alpha_n)$ and $\tilde{E}(u_n^{\alpha_n},w_n^{\alpha_n},\alpha_n)$ as a positive flux - difference between the members of the whole minimizing sequence and those of the sequence solutions of the intermediate direct problems.
	\newline We exploit now the convergence results presented in Propositions \ref{prop3} and \ref{prop4}. Then, clearly we have \begin{align}
	\lim_{n \rightarrow \infty} \tilde{E}(u_n^{\alpha_n},w_n^{\alpha_n},\alpha_n)=\tilde{E}(u^{\alpha},w^{\alpha},\alpha). \label{stack}
	\end{align} 
	Taking the limit $n \rightarrow \infty$ in (\ref{bf}) and exploiting (\ref{paok1}),(\ref{paok2}),(\ref{stack}), we obtain that
\begin{align}
&\underset{(u,w,\alpha^\prime)}{\inf} E(u,w,\alpha^\prime)=\lim_{n \rightarrow \infty} E(u_n,w_n,\alpha_n) \nonumber\\ &\geq \underset{n \rightarrow \infty}{\lim}\tilde{E}(u_n^{\alpha_n},w_n^{\alpha_n},\alpha_n) +\frac{\lambda}{2} \underset{n \rightarrow \infty}{\underline{\lim}}\|R\alpha_n-\alpha_0\|^2_{L^2(\Omega)} \nonumber \\ & + \mu \underset{n \rightarrow \infty}{\underline{\lim}}\int_{\Omega} \phi(|D \alpha_n|) +\frac{\kappa+1}{2}\underset{n \rightarrow \infty}{\underline{\lim}}\int_{\Omega}\alpha_n|\nabla (u_n-u_n^{\alpha_n})|^2\nonumber \\ &+\frac{1-\kappa}{2}\underset{n \rightarrow \infty}{\underline{\lim}} \int_{\Omega}\alpha_n|\nabla (w_n-w_n^{\alpha_n})|^2 \nonumber \\ &\geq \tilde{E}(u^{\alpha},w^{\alpha},\alpha)+\frac{\lambda}{2} \|R\alpha-\alpha_0\|^2_{L^2(\Omega)}+ \mu \int_{\Omega} \phi(|D \alpha|)\nonumber \\ &+ \frac{\kappa+1}{2}\underset{n \rightarrow \infty}{\underline{\lim}}\int_{\Omega}\alpha_n|\nabla (u_n-u_n^{\alpha_n})|^2\nonumber \\ & +\frac{1-\kappa}{2}\underset{n \rightarrow \infty}{\underline{\lim}} \int_{\Omega}\alpha_n|\nabla (w_n-w_n^{\alpha_n})|^2 \Rightarrow \nonumber \\
 &\underset{(u,w,\alpha^\prime)}{\inf} E(u,w,\alpha^\prime) \geq E(u^{\alpha},w^{\alpha},\alpha)\nonumber \\ &+\frac{\kappa+1}{2}\underset{n \rightarrow \infty}{\underline{\lim}}\int_{\Omega}\alpha_n|\nabla (u_n-u_n^{\alpha_n})|^2 \nonumber \\
& + \frac{1-\kappa}{2}\underset{n \rightarrow \infty}{\underline{\lim}} \int_{\Omega}\alpha_n|\nabla (w_n-w_n^{\alpha_n})|^2. \label{good}
\end{align}
	On the basis of the convergences (\ref{guer1}),(\ref{guer2}), the limit processes declared by Propositions \ref{prop3},\ref{prop4} and the lower semi-continuity of the $L^2(\Omega)$-norm in the weak topology, we obtain
	\begin{align}
	&\underset{n \rightarrow \infty}{\underline{\lim}}\int_{\Omega}\alpha_n|\nabla (u_n-u_n^{\alpha_n})|^2 \nonumber \\ &\geq b \underset{n \rightarrow \infty}{\underline{\lim}}\int_{\Omega}|\nabla (u_n-u_n^{\alpha_n})|^2 \geq b \int_{\Omega}|\nabla (u-u^{\alpha})|^2, \\
	&\underset{n \rightarrow \infty}{\underline{\lim}} \int_{\Omega}\alpha_n|\nabla (w_n-w_n^{\alpha_n})|^2 \nonumber \\ &\geq b \underset{n \rightarrow \infty}{\underline{\lim}} \int_{\Omega}|\nabla (w_n-w_n^{\alpha_n})|^2 \geq b \int_{\Omega}|\nabla (w-w^{\alpha})|^2.
	\end{align}   
	Consequently the relation (\ref{good}) becomes
	\begin{align}
	\underset{(u,w,\alpha^\prime)}{\inf} E(u,w,\alpha^\prime) &\geq E(u^{\alpha},w^{\alpha},\alpha)\nonumber\\  &+b\frac{\kappa+1}{2}\int_{\Omega}|\nabla (u-u^{\alpha})|^2 \nonumber\\
	& +b\frac{1-\kappa}{2} \int_{\Omega}|\nabla (w-w^{\alpha})|^2. \label{excellent}
	\end{align}
	The last inequality can be satisfied only if $\nabla u= \nabla u^{\alpha}$, $\nabla w = \nabla w^\alpha$, a.e. in $\Omega$ or equivalently 
	\begin{align}
	u=u^{\alpha} \qquad &\mbox{in}\ H^{1}_{0}(\Omega), \\
	w=w^{\alpha} \qquad &\mbox{in}\ H^{1}(\Omega)/\mathbb{R}.  
	\end{align}
Necessarily, the inequality (\ref{excellent}) becomes the equation $\underset{(u,w,\alpha^\prime)}{\inf} E(u,w,\alpha^\prime)=E(u,w,\alpha)$, which has the interesting interpretation
	\begin{align}
	\underset{(u,w,\alpha^\prime)}{\inf} E(u,w,\alpha^\prime)&=\lim_{n \rightarrow \infty} E(u_n,w_n,\alpha_n)\nonumber\\ &=E(\lim_{n \rightarrow \infty} u_n, \lim_{n \rightarrow \infty} w_n, \lim_{n \rightarrow \infty} \alpha_n)\nonumber\\ &=E(u^\alpha,w^\alpha,\alpha).\nonumber \\
	\end{align}
	So the minimization problem admits a solution, which necessarily is the limit of the minimizing sequence.

\begin{remark} \label{rem1}
	It is important to notice that one of the fundamental outcomes of the Theorem \ref{thfund} is that it is not necessary to solve a direct conductivity problem in every step of the process (with a specific conductivity profile). This would of course accelerate the convergence of the minimization process but the descent toward the minimum could be implemented independently for the three variables $u,w,\alpha$ of the problem. 
\end{remark}
\begin{remark}
	Without loss of generality, the results of Theorem \ref{thfund} hold, even if the members of the minimizing sequence belong themselves to $BV(\Omega)\cap L^{\infty}_{\delta}(\Omega, [b,c])$ (This space is closed with respect to the particular kind of activity the sequence elements participate in). The specific form of the theorem has been adopted, since it is common to use step wise constant functions as approximation profiles, and so as to reveal that even starting with constant profiles the minimization process is well accomplished. This gives some extra degrees of freedom to the arsenal of the numerical implementation.
\end{remark}
Theorem \ref{thfund} states that the infima of the functional occur at triples $(u^{\alpha},w^{\alpha},\alpha)$, where the functions $u^{\alpha},w^{\alpha}$ solve the Dirichlet and Neumann boundary value problems corresponding to the conductivity profile $\alpha$ and the data $f,g$. Actually, this would emerge naturally if we solved repeatedly the direct problem for every member $\alpha_n$ of the convergent sequence of profiles, as stated in Remark \ref{rem1}. If this was the plan, then the functional part $\tilde{E}(u_n,w_n,\alpha_n)$ would be replaced immediately by the term $\tilde{E}(u_n^{\alpha_n},w_n^{\alpha_n},\alpha_n)$ and the total functional would depend only on the conductivity function $\alpha$,i.e   
\begin{align}
&E_1(\alpha)=\tilde{E}(u^{\alpha},w^{\alpha},\alpha)+\frac{\lambda}{2}\|R\alpha-\alpha_0\|^2_{L^2(\Omega)}\nonumber \\ &+ \mu \int_{\Omega} \phi(|D \alpha|)\nonumber \\ &=\frac{\kappa+1}{2}\left\langle \Lambda_{\alpha}f,f \right\rangle_{\mathcal{A} \times \mathcal{B}} -\frac{(1-\kappa)}{2}\left\langle g,\Lambda_{\alpha}^{-1}g  \right\rangle_{\mathcal{A} \times \mathcal{B}}\nonumber \\ & -\kappa\left\langle g,f \right\rangle_{\mathcal{A} \times \mathcal{B}}+\frac{\lambda}{2}\|R\alpha-\alpha_0\|^2_{L^2(\Omega)}+\mu \int_{\Omega} \phi(|D \alpha|) \nonumber \\
 &=\left\langle g,f-\Lambda_{\alpha}^{-1}g  \right\rangle_{\mathcal{A} \times \mathcal{B}}+\frac{\kappa+1}{2}\left\langle g-\Lambda_{\alpha} f,\Lambda_{\alpha}^{-1}g-f  \right\rangle_{\mathcal{A} \times \mathcal{B}}\nonumber \\ &+ \frac{\lambda}{2}\|R\alpha-\alpha_0\|^2_{L^2(\Omega)}+\mu \int_{\Omega} \phi(|D \alpha|).\nonumber\\
\label{totalpa}
\end{align}  
In Theorem \ref{thfund}, the parameter $\kappa$ has been taken to belong strictly in the interval $(-1,1)$, so that the concurrent minimization over fields and profiles shares coercive behavior. However, in the scheme (\ref{totalpa}), where only the conductivity profile appears, the parameter $\kappa$ can be released and detour the critical value $1$. Thus, taking a sufficiently large parameter $\kappa$ would render the term $\left\langle g-\Lambda_{\alpha} f,\Lambda_{\alpha}^{-1}g-f  \right\rangle_{\mathcal{A} \times \mathcal{B}}=\left\langle g-\Lambda_{\alpha} f,\Lambda_{\alpha}^{-1}(g-\Lambda_{\alpha}f)  \right\rangle_{\mathcal{A} \times \mathcal{B}}$ primitively significant. It is worthwhile to notice that if it was additionally known that $g,\Lambda_{\alpha}f \in L^2(\partial \Omega)$ (some slight additional regularity), then the second term of $E_1$ would be responsible for the suppression of the difference $\|g-\Lambda_{\alpha}f{\|}_{L^2(\partial \Omega)}^2$ between data, over the measurement surface $\partial \Omega$. This is the commonly used minimization functional forcing the surface data to be in accordance. 

Theorem \ref{thfund} establishes existence of a minimum, but does not offer any argument concerning some kind of uniqueness. In fact, it is not expected to have a unique solution working with just a pair of data $(f,g)$. As far as convexity is concerned, a more attentive examination of the terms participating in the formation of the functional $E_1(\alpha)$ leads to the result that there exist opposing terms. More precisely, the last two terms of (\ref{totalpa}), controlling purely the functional structure of the conductivity profiles, are convex \cite{aubert}. The first two terms, emerging from the treatment of the involved direct problems, contain necessarily a concave term. Indeed, working rather with the decomposition 
\begin{align}
	\tilde{E}(u^{\alpha},w^{\alpha},\alpha)&=\frac{\kappa+1}{2}\left\langle \Lambda_{\alpha}f,f \right\rangle_{\mathcal{A} \times \mathcal{B}}\nonumber\\ &-\frac{(1-\kappa)}{2}\left\langle g,\Lambda_{\alpha}^{-1}g  \right\rangle_{\mathcal{A} \times \mathcal{B}}\nonumber \\ & -\kappa\left\langle g,f \right\rangle_{\mathcal{A} \times \mathcal{B}},
\end{align}   
we see that the first term is concave. This can been shown by considering two discrete conductivity profiles $\alpha_i,\ i=1,2$, along with the accompanying fields $u_i,w_i, \ i=1,2$, which solve the corresponding direct problems as described in Theorem \ref{thfund}. Let us select an arbitrary convex combination $\alpha_s=s\alpha_1+(1-s)\alpha_2$ ($0<s<1$) of the profiles $\alpha_i, i=1,2$ (with corresponding fields $u_s,w_s$). It turns out obviously, that this function is an admissible conductivity function.  
We remark that
\begin{align}
\left\langle\Lambda_{\alpha_s}f,f \right\rangle_{\mathcal{A} \times \mathcal{B}}&=\int_{\Omega}\alpha_s|\nabla (u_s+\eta f)|^2\geq s \int_{\Omega}\alpha_1|\nabla(u_1+\eta f)|^2 \nonumber \\ &+(1-s) \int_{\Omega}\alpha_2|\nabla (u_2+\eta f)|^2  \nonumber \\
&=s\left\langle \Lambda_{\alpha_1}f,f \right\rangle_{\mathcal{A} \times \mathcal{B}}+(1-s)\left\langle \Lambda_{\alpha_2}f,f \right\rangle_{\mathcal{A} \times \mathcal{B}},
\end{align}
which proves that (when of course $\kappa+1>0$) the term $\frac{\kappa+1}{2}\left\langle \Lambda_{\alpha}f,f \right\rangle_{\mathcal{A} \times \mathcal{B}}$ is concave. In addition 
\begin{align}
&-\left\langle g,\Lambda_{\alpha_s}^{-1}g \right\rangle_{\mathcal{A} \times \mathcal{B}}=\int_{\Omega}\alpha_s|\nabla w_s|^2-2\left\langle g,\gamma w_s \right\rangle_{\mathcal{A} \times \mathcal{B}}\nonumber \\ &=s \int_{\Omega}\alpha_1|\nabla w_s|^2 +(1-s)\int_{\Omega}\alpha_2|\nabla w_s|^2 -2 \left\langle g, \gamma w_s \right\rangle_{\mathcal{A} \times \mathcal{B}}\nonumber \\ & \geq s \left[\int_{\Omega}\alpha_1|\nabla w_1|^2- 2\left\langle g, \gamma w_1 \right\rangle_{\mathcal{A} \times \mathcal{B}}\right] \nonumber \\ & +(1-s)\left[\int_{\Omega}\alpha_2|\nabla w_2|^2-2\left\langle g, \gamma w_2 \right\rangle_{\mathcal{A} \times \mathcal{B}}\right] \nonumber \\ &=s\left[-\left\langle g,\Lambda_{\alpha_1}^{-1}g \right\rangle_{\mathcal{A} \times \mathcal{B}}\right]+(1-s)\left[-\left\langle g,\Lambda_{\alpha_2}^{-1}g \right\rangle_{\mathcal{A} \times \mathcal{B}}\right]. \nonumber \\
\label{conv11}
\end{align}
Then, even the second term $-\frac{(1-\kappa)}{2}\left\langle g,\Lambda_{\alpha}^{-1}g  \right\rangle_{\mathcal{A} \times \mathcal{B}}$ of $\tilde{E}(u^{\alpha},w^{\alpha},\alpha)$ is concave when $\kappa<1$. In case that $\kappa>1$, this term becomes convex and acts in collaboration with the term $\frac{\lambda}{2}\|R\alpha-\alpha_0\|^2_{L^2(\Omega)}+\mu \int_{\Omega} \phi(|D \alpha|)$. Indeed, when we perform a mixed minimization over $(u,w,\alpha)$ (and then respect the coercivity condition $\kappa<1$) we have more local minima than in the case of minimizing $E_1$ over $\alpha$, when the relaxed selection $\kappa>1$ can be adopted. Nevertheless, even the commonly used functional $\left\langle g-\Lambda_{\alpha} f,\Lambda_{\alpha}^{-1}g-f  \right\rangle_{\mathcal{A} \times \mathcal{B}}$ written as $\left\langle g,\Lambda_{\alpha}^{-1}g  \right\rangle_{\mathcal{A} \times \mathcal{B}}+$$\left\langle \Lambda_{\alpha} f,f  \right\rangle_{\mathcal{A} \times \mathcal{B}}-$$2\left\langle g,f  \right\rangle_{\mathcal{A} \times \mathcal{B}}$, disposes adversarial terms. 

\begin{remark}
	A well known result (see for example \cite{butazzo}, p 409) holds, according to which when a $BV-$function $\alpha$ belongs additionally to $L^{\infty}(\Omega)$ and there exists a closed subset $K$ of $\Omega$ (with finite Hausdorff measure, i.e. $\mathcal{H}^{N-1}(K) < \infty$) such that $\alpha \in W^{1,1}(\Omega \backslash K)$, then $\alpha$ is necessarily a function of $SBV(\Omega)$, whose distributional derivative  disposes no Cantor part, while its jump part is a subset of $K$ (i.e. $S_\alpha \subset K$). This result assures that all the possible discontinuities lay in the set $K$ and the conductivities emerging as members of minimization sequences or limit functions, cannot dispose more intrinsic derivatives than the Young measures, delimiting the support of interfaces. 
\end{remark}

\section{Towards the numerical implementation of the inverse conductivity problem} \label{tow}
The following theorem is the first stage of the investigation stepping to the establishment of the suitable numerical scheme.
\begin{theorem}
	The functional $E(u,w,\alpha)$ expressed by (\ref{total}) and (\ref{total1}) is lower semicontinuous with respect to the strong $L^1(\Omega)-$convergence on $\alpha$ and the weak $L^1(\Omega) \times L^1(\Omega)-$convergence on $(\nabla u, \nabla w )$ (This combined convergence defines a topology, which will next be denoted $\tau-$topology). \label{th2} 
\end{theorem} 
{\bf Proof.}
	We evoke herein mainly the Theorems 13.1.1 and 16.4.1 of \cite{butazzo} to establish the sought semicontinuity. More precisely, we consider the function $h(\alpha,v)$, where $v=(v_1,v_2) \in \mathbb{R}^N \times \mathbb{R}^N$, defined as   
	\begin{align}
	h(\alpha,v)= \frac{\kappa+1}{2} \alpha | (v_1+\nabla (\eta f))|^2 +\frac{(1-\kappa)}{2}\alpha |v_2|^2 . 
	\end{align}
	Since the function $h(\alpha,v)$ is convex in $v$ and lower semicontinuous in $\alpha$, we infer, on the basis of the aforementioned theorems of \cite{aubert}, that the functional 
	\begin{align}
	H(\alpha,v)=\int_{\Omega}h(\alpha(x),v(x)) dx
	\end{align}
	is lower semicontinuous with respect to the strong $L^1-$convergence on $\alpha$ and the weak $L^1-$convergence on $v$. We take $v=(v_1,v_2)=(\nabla u, \nabla w)$ and consider sequences $(u_n,w_n)$ with the properties implied by Theorem \ref{thfund}, according to which $\nabla u_n \underset{L^2(\Omega)}{-\!\!\!\rightharpoonup} \nabla u$ and $\nabla w_n \underset{L^2(\Omega)}{-\!\!\!\rightharpoonup} \nabla w$. Then, the pair $(\nabla u_n,\nabla w_n)$ converges $L^1-$weakly to $(u,w)$, while $\alpha_n \rightarrow \alpha$ strongly. Incorporating the term $-\kappa\left\langle g,f \right\rangle_{\mathcal{A} \times \mathcal{B}}$, we construct the functional  $\tilde{E}(u,w,\alpha)$ given by (\ref{total1}), which obeys then to the semicontinuity result:
	\begin{align}
	\tilde{E}(u,w,\alpha) \leq \underset{n \rightarrow \infty}{\underline{\lim}}\tilde{E}(u_n,w_n,\alpha_n).
	\end{align}
Taking into consideration (\ref{paok1}) and (\ref{paok2}), it is finally proved that the total functional (\ref{total}) satisfies 
	\begin{align}
	{E}(u,w,\alpha) \leq \underset{n \rightarrow \infty}{\underline{\lim}}{E}(u_n,w_n,\alpha_n).
	\end{align}

It is more informative to give the following extended form of the lower semicontinuous functional ${E}(u,w,\alpha)$ \cite{butazzo}
\begin{align}
&E(u,w,\alpha)=\tilde{E}(u,w,\alpha)+\frac{\lambda}{2}\|R\alpha-\alpha_0\|^2_{L^2(\Omega)} \nonumber \\
&+\mu \int_{\Omega} \phi(|\nabla \alpha|) dx+\mu t \int_{S_\alpha}|\alpha^{+}-\alpha^{-}| d\mathcal{H}^{N-1}\nonumber \\ &+\mu t \int_{\Omega-S_\alpha}|C_\alpha|, \label{before}
\end{align}
where we recognize the decomposition of the total mass into three parts: the absolutely continuous part with respect to the Lebesgue measure, the jump part and the Cantor measure. The parameter $t$ stands for the limit $\lim_{s \rightarrow \infty}\frac{\phi(s)}{s}$. This makes sense, since the strictly convex function $\phi$ is selected to be non decreasing function from $\mathbb{R}^{+}$ to $\mathbb{R}^{+}$, with $\phi(0)=0$, while the condition $ts-\tau \leq \phi(s) \leq ts+\tau$ with specific constants $t>0$ and $\tau \geq 0$, is satisfied for every $s \geq 0$. \\
The functional (\ref{before}) could be exactly the target of the minimization process. However, it is common to use quadratic schemes in the optimization setting. We construct then, inspired by \cite{aubert}, the auxiliary function 
\begin{align}
\phi_\epsilon(s)=
\left\{
\begin{array} {llrr}
\frac{\phi^\prime(\epsilon)}{2 \epsilon}s^2+\phi(\epsilon)-\frac{\epsilon \phi^\prime(\epsilon)}{2} \ 
&\mbox{if}& \ 0 \leq s \leq\epsilon, \\ 
\phi(s) \ 
&\mbox{if}& \ \epsilon \leq s \leq \frac{1}{\epsilon}, \\ 
\frac{\epsilon \phi^\prime(\frac{1}{\epsilon})}{2}s^2+\phi(\frac{1}{\epsilon})-\frac{\phi^\prime(\frac{1}{\epsilon})}{2\epsilon} \ 
&\mbox{if}& \ s \geq \frac{1}{\epsilon} . \end{array} \right. \label{phie}
\end{align} 
Clearly, for every $\epsilon$, it holds that $\phi_\epsilon(s) \geq 0$, $\forall s$ and $\lim_{\epsilon \rightarrow 0} \phi_\epsilon(s) = \phi(s)$. We  set $\mathcal{E}=H_0^1(\Omega) \times H^{1}( \Omega)/\mathbb{R}$ and define now the functional $E_\epsilon$  as follows 
\begin{align}
&E_\epsilon(u,w,\alpha)=\nonumber \\ &\left\{\begin{array} {llll} \tilde{E}(u,w,\alpha)+\frac{\lambda}{2}\|R\alpha-\alpha_0\|^2_{L^2(\Omega)} +\mu \int_{\Omega} \phi_\epsilon(|\nabla \alpha|)dx,  \\ \ \ \mbox{if} \ (u,w) \in \mathcal{E} \ \mbox{and}\ \alpha \in H^1(\Omega) \cap L_\delta^\infty(\Omega,[b,c]), \\ \infty, \\ \ \ \mbox{if} \ (u,w) \in \mathcal{E} \ \mbox{and} \\ \ \ \alpha \in (BV(\Omega)-H^1(\Omega)) \cap L_\delta^\infty(\Omega,[b,c]). \end{array}\right. \nonumber \\
\end{align}
Clearly, evoking classical arguments, simpler but reminiscent of the steps followed in Theorem \ref{thfund}, we find that
\begin{proposition}
	For each $\epsilon>0$, the functional $E_\epsilon$  admits a minimum $(u_\epsilon,w_\epsilon,\alpha_\epsilon)$ in $\mathcal{E} \times (H^1(\Omega)\cap L_\delta^\infty(\Omega,[b,c]))$. \label{rega}
\end{proposition} 
In addition, we define
\begin{align}
\overline{E}(u,w,\alpha)=\left\{ \begin{array} {llll} E(u,w,\alpha) \nonumber \quad \mbox{if} \ (u,w) \in \mathcal{E} \ \mbox{and}\\ \hspace{1.5cm} \alpha \in H^1(\Omega)\cap L_\delta^\infty(\Omega,[b,c]),  \\ \infty  \hspace{1.5cm} \mbox{if} \ (u,w) \in \mathcal{E} \ \mbox{and} \\ \ \alpha \in (BV(\Omega)-H^1(\Omega)) \cap L_\delta^\infty(\Omega,[b,c]). \end{array}\right.\nonumber\\
\end{align}
\begin{proposition}
	The lower semicontinuous envelope $\mathcal{R}_\tau \overline{E}$ of $\overline{E}$, with respect to the strong $L^1(\Omega)-$convergence on $\alpha$ and the weak $L^1(\Omega) \times L^1(\Omega)-$convergence on $(\nabla u, \nabla w )$, coincides with $E$.   \label{aintee}
\end{proposition}
{\bf Proof.} It is proved in Theorem \ref{th2}, that $E$ is lower semicontinuous with respect to the mixed type $\tau-$topology stated there. Given that clearly $\overline{E} \geq E$, what remains to be proved in order that $\mathcal{R}_\tau \overline{E}=E$ is that there exists a sequence $(u_h,w_h,\alpha_h)$ in $\mathcal{E} \times (H^1(\Omega)\cap L_\delta^\infty(\Omega,[b,c]))$ converging in the $\tau$ topology to $(u,w,\alpha)$ and $E(u,w,\alpha)=\underset{h \rightarrow 0}{\underline{\lim}}\overline{E}(u_h,w_h,\alpha_h)$. Such a sequence can be constructed using classical approximation theorems \cite{evans, demengel}. We just pay attention on the exploitation of the intermediate convergence of measures, assuring the appropriate convergence of the involved measures. 

\begin{proposition}
	Every cluster point $(u,w,\alpha)$ - with respect to the strong $L^1(\Omega)-$conver\-gence on $\alpha$ and the weak $L^1(\Omega) \times L^1(\Omega)-$convergence on $(\nabla u, \nabla w )$ - of the sequence $\left\{(u_\epsilon,w_\epsilon,\alpha_\epsilon); \epsilon>0\right\}$, introduced in Proposition \ref{rega}, is a minimizer of $E$ and $E_\epsilon(u_\epsilon,w_\epsilon,\alpha_\epsilon)$ converges to $E(u,w,\alpha)$.
\end{proposition}
{\bf Proof.} By construction, it is obvious that $E_\epsilon(u,w,\alpha)$ is a decreasing sequence converging pointwise to $\overline{E}(u,w,\alpha)$. According to a well known result of the theory of $\Gamma-$convergence (Theorem 2.1.8 in \cite{aubert}), $E_\epsilon$ $\Gamma-$converges to the lower semicontinuous envelope $\mathcal{R}_\tau \overline{E}$ of $\overline{E}$, which according to Proposition \ref{aintee}, coincides with the functional $E$. In addition, thanks to the uniform bounds $b,c$ of the conductivity coefficients the functionals $E_\epsilon$ are equicoercive and so the set $\left\{(u_\epsilon,w_\epsilon,\alpha_\epsilon); \epsilon>0\right\}$ is relatively compact in the aforementioned topologies. According to Theorem 12.1.1 of \cite{butazzo}, every cluster point $(u,w,\alpha)$ of the set $\left\{(u_\epsilon,w_\epsilon,\alpha_\epsilon); \epsilon>0\right\}$ is a minimizer of $E$ and $E_\epsilon(u_\epsilon,w_\epsilon,\alpha_\epsilon)$ converges to $E(u,w,\alpha)$ as $\epsilon$ goes to zero.  

\section{The numerical implementation} \label{numerics}
Following the arguments presented in Section 3.2.4 of \cite{aubert}, we are interested in constructing a semiquadratic algorithm implementing numerically the optimization scheme described in the previous Section. More precisely, we restrict ourselves to the case that the auxiliary function $\phi$, introduced above and participating in the part of the functional controlling the masses of the conductivity profiles, is selected to satisfy the requirements of Proposition 3.2.4 (encountered in \cite{aubert}). The main condition is that the function $\phi(\sqrt{s})$ is concave in $(0,\infty)$, while we recall that $\phi(s)$ is a non-decreasing function. Consequently, the non negative parameters $L=\lim_{s \rightarrow \infty}\frac{\phi^\prime(s)}{2s}$ and $M=\lim_{s \rightarrow 0}\frac{\phi^\prime(s)}{2s}$ are well defined and obey to the ordering $L<M$. According to the main results of the Proposition 3.2.4, there exists a convex and decreasing function $\psi:]L,M] \rightarrow [\beta_1, \beta_2]$ such that
\begin{align}
\phi(s)=\underset{L \leq \omega \leq M}{\inf} (\omega s^2 +\psi(\omega)),
\end{align}
where $\beta_1=\lim_{s \rightarrow 0^+}\phi(s)$ and $\beta_2=\lim_{s \rightarrow +\infty}(\phi(s)-\frac{s \phi^\prime(s)}{2})$.
Moreover, for every $s \geq 0$, the value $\omega$ at which the minimum is reached is exactly $\omega=\frac{\phi^\prime(s)}{2s}$. This approach has the characteristics of a duality process and the variable $\omega$ is called the dual variable. \\ 
Applying this duality approach to the - quadratic near the end points - auxiliary function $\phi_\epsilon$, given by (\ref{phie}), we have the possibility \cite{aubert} to express the $\underset{(u,w,\alpha)}{\inf}E_\epsilon(u,w,\alpha)$ in the form
\begin{align}
\underset{(u,w,\alpha)}{\inf}E_\epsilon(u,w,\alpha)=\underset{\omega}{\inf}\underset{(u,w,\alpha)}{\inf} J_\epsilon(u,w,\alpha,\omega),\label{problemtoreconstruct}
\end{align}    
where
\begin{align}
J_\epsilon(u,w,\alpha,\omega)=&\tilde{E}(u,w,\alpha)+\frac{\lambda}{2}\|R\alpha-\alpha_0\|^2_{L^2(\Omega)}\nonumber \\ &+\mu \int_\Omega (\omega |\nabla \alpha{|}^2+\psi_\epsilon(\omega)) dx\label{funtoreconstruct}
\end{align}
with $(u,w) \in \mathcal{E}$ and $\alpha \in H^1(\Omega)$. 
 \begin{remark}\label{discandhom}
	The sequence of functions $\omega^{n}(x)$ can be seen as an indicator of contours.\\
	1. If $\omega^n(x) \approx 0$, then $x$ belongs to a discontinuity surface. \\
	2. If $\omega^n(x) \approx 1$, then $x$ belongs to a homogeneous region. 
\end{remark}

\subsection{$BV$ Regularized Inversion algorithm}\label{subsec:Algorithm}
The main property of the semiquadratic algorithm we apply on (\ref{funtoreconstruct}) is to detect progressively the discontinuities of the conductivity profile $\alpha$. Our computational process will be started from an initial estimation of $(u,w,\alpha,\omega)$. Following \cite{aubert}, we have further taken advantage of the role of the dual variable, as attributed in Remark \ref{discandhom}. The initial guess $\omega^0$ allows us to introduce into the formulation of the minimization problem (\ref{problemtoreconstruct}) prior information about the location of the inhomogeneity's boundary. 
 The $BV$ Reconstruction algorithm is then read as: given a step size control $h\in (0,1)$, an initial approximation $(u^0,w^0,\alpha^0,\omega^0)$, a smoothing parameter $\epsilon$, a maximum number of iterations ($\max$\_iters $=10$) and setting the regularization parameter $\lambda=0$:

 \begin{algorithm}[H]
 	 	\caption{$BV$ {\it Regularized Inversion}}
 	\begin{algorithmic}[1]
 		\Require $(u^0,w^0,\alpha^0)$
 		 \Comment{{\footnotesize initial guess for $(u^0,w^0,\alpha^0,\omega^0)$}}
 		\Ensure  $u_{\min},w_{\min},\alpha_{\min},\omega_{\min}$
 	
		\Function{BVReconstruction}{$u^0,w^0,\alpha^0,\omega^0$}
        \For{ $n \gets 0$ : max\_iters} 
 		\State $(u^{n+1}_\epsilon,w^{n+1}_\epsilon,\alpha_\epsilon^{n+1})
 		\gets {\footnotesize\underset{(u,w,\alpha)}{\mbox{argmin}}}J_\epsilon(u,w,\alpha,\omega^n)$
 		\label{opt_step}
 		\State $\omega^{n+1} \gets \frac{\phi_{\epsilon}^\prime(|\nabla \alpha_\epsilon^{n+1}|)}{2|\nabla \alpha_\epsilon^{n+1}|}$
 		\label{omega_step}
 		\EndFor
 		\State \Return {$u_{\min},w_{\min},\alpha_{\min},\omega_{\min}$}
 		\EndFunction
 	\end{algorithmic}
 	\label{alg:BV}
 \end{algorithm}
 
\noindent At the $n$th iteration of the algorithm, the variables $u,w,\alpha$ are computed for given $\omega$ as the $\underset{(u,w,\alpha)}{\mbox{argmin}}J_\epsilon(u,w,\alpha,\omega^n)$ (Step 3). Since we would like to exploit knowledge about lower and upper bounds on the conductivity, we adopt the interior point method implemented in the software package IPOPT~\cite{IPOPT2005,IPOPT2006}, which accepts lower and upper bounds on the control parameters, can exploit gradients, and approximates the Hessian using BFGS updates \cite{Fletcher}. Gradient-based optimizers, such as IPOPT, typically terminate as soon as the dual infeasibility has been decreased below a specified threshold.
Then, $\omega^{n+1}$ is computed as
  $\omega^{n+1} = \frac{\phi_{\epsilon}^\prime(|\nabla \alpha_\epsilon^{n+1}|)}{2|\nabla \alpha_\epsilon^{n+1}|} =  \underset{\omega}{\mbox{argmin}}J_\epsilon(u_\epsilon^{n+1},w_\epsilon^{n+1},\alpha_\epsilon^{n+1},\omega)$ (Step 4).
 
  \subsection{Numerical tests}\label{subsec:numtests}
 We now illustrate the theoretical results with numerical tests. For this purpose we consider the boundary value problem
 \begin{align}\nonumber
\nabla \cdot (\alpha(x) \nabla u(x))&=0 \ \ && x \in \Omega, \\ \nonumber
u(x)&=f(x) \ \ && x \in \partial \Omega, \\ \nonumber
\alpha(x) \frac{\partial u}{\partial n} (x)&=g(x) \ \ && x \in \partial \Omega,  \nonumber
\end{align}   
with known data $(f,g)$ and the conductive region $\Omega$ be the disc
\begin{equation}
\Omega:=\left\{(x,y) \in \mathbb{R}^2 :\rho =\sqrt{x^2+y^2} < 2 \right\}.
\end{equation} The threshold values $b$ and $c$ in the definition domain of $\alpha(x)$ can vary accordingly. \newline We examine the efficiency of the conductivity reconstruction method via a single measurement for one concentric and two different eccentric inhomogeneities, each represented by the inclusion $D$ defined as 
\begin{equation}
    D:=\left\{ (x,y)\in \mathbb{R}^2: \rho_D=\sqrt{(x-x_{D})^2+(y-y_D)^2}< 1\right\}
\end{equation}
with $(x_D,y_D)$ its center coordinates. 
We assume that the - known a priori though target of reconstruction - conductivity profile has the following simple profile
\begin{align}
\alpha =\alpha(x,y)=\begin{cases} 
  2, \qquad (x,y) \in D,  \\  1, \qquad (x,y) \in \Omega\setminus D.  
\end{cases} \label{june1}
\end{align} \newline Without loss of generality, the parameter $b$ has been selected - in the most part of our experiments - equal to $1$ in order to comply, in a simple manner, with the requirements on the protective region $\Omega_{\delta}$. The threshold parameter $c$ has been chosen to be not less than the critical value $2$. \newline
 We discretize the computational domain, defined by the outer circle ($\rho=2$), using Triangle~\cite{shewchuk96b, shewchuk2002a} with characteristic mesh size ($h$), reported in Table \ref{table:infomesh} below along with the number of elements and the number of nodes for each mesh. We adopt the piecewise linear, continuous family of finite elements (P1) for the discretization of the operators $J_\epsilon(u,w,\alpha,\omega)$. 

\begin{figure}[!ht]\onecolumn
  \centering
 	\hspace{-1.5cm}
	\begin{minipage}{.08\textwidth}
  	\vspace{0.2cm}
	\includegraphics[width=0.71\textwidth,height=0.50\textheight]{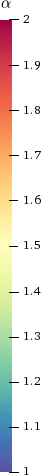}
 \end{minipage} \hspace{.3cm}
	\centering
\begin{minipage}{.96\textwidth}
\centering
		{\footnotesize{$\boldsymbol{\ell=0.2 \ }$ \; }} \subfloat[$\omega^0$]{\includegraphics[width=0.19\linewidth]{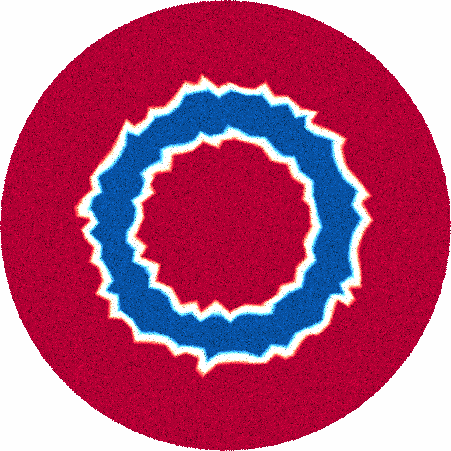}}
		\subfloat[$\alpha^0$]{\includegraphics[width=0.19\linewidth]{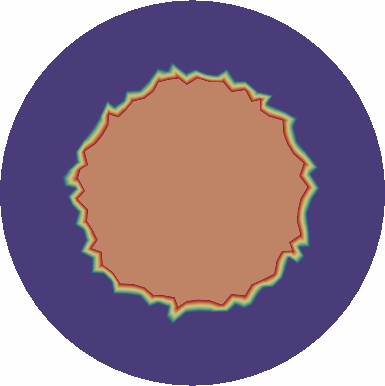}}
		\subfloat[$\omega$]{\includegraphics[width=0.19\linewidth]{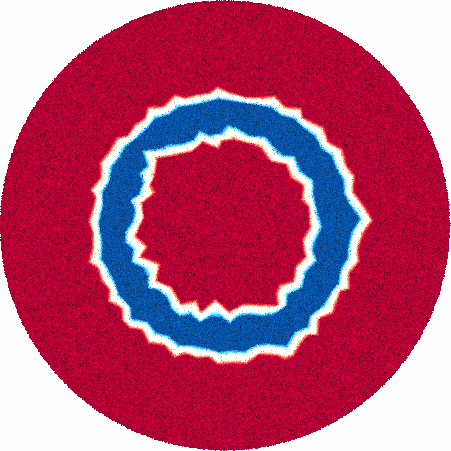}}
		\subfloat[$\alpha$]{\includegraphics[width=0.19\linewidth]{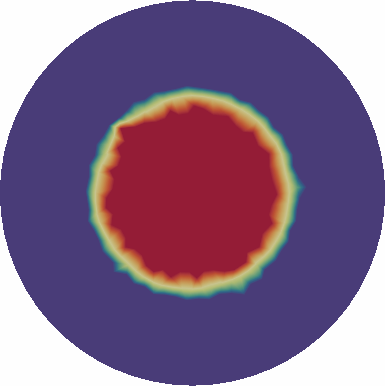}}\\
		{\footnotesize{$\boldsymbol{\ell=0.3 \ }$ \; }} \subfloat[$\omega^0$]{\includegraphics[width=0.19\linewidth]{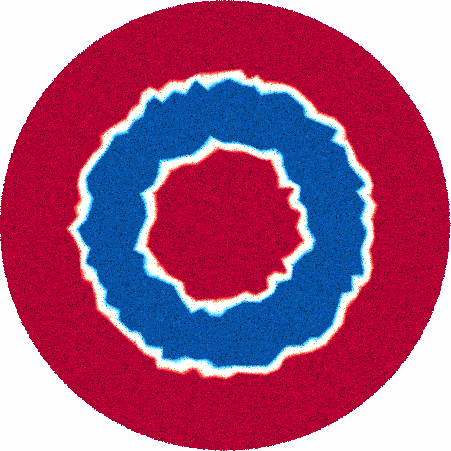}}
		\subfloat[$\alpha^0$]{\includegraphics[width=0.19\linewidth]{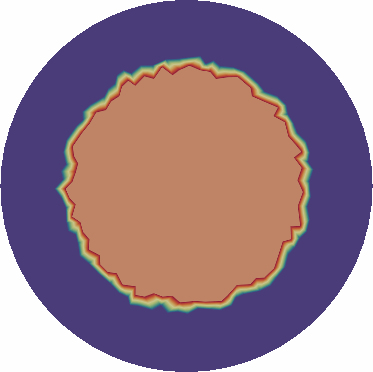}}
		\subfloat[$\omega$]{\includegraphics[width=0.19\linewidth]{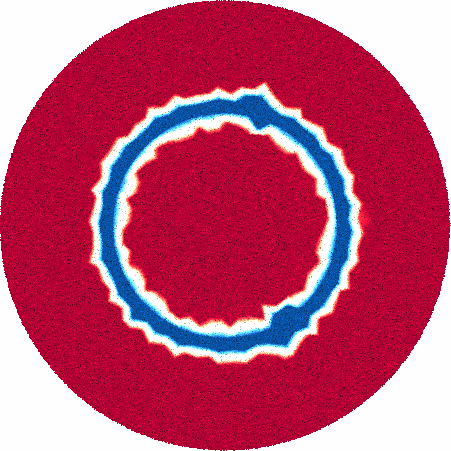}}
		\subfloat[$\alpha$]{\includegraphics[width=0.19\linewidth]{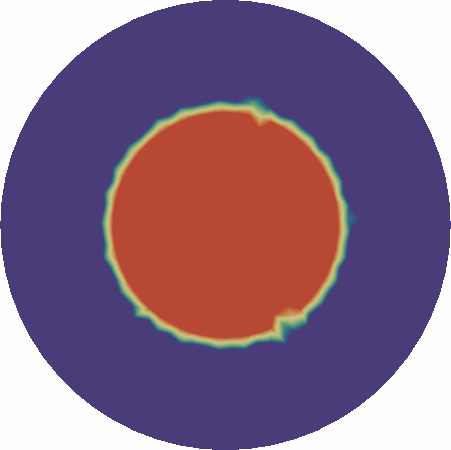}}\\
		{\footnotesize{$\boldsymbol{\ell=0.4 \ }$ \; }} \subfloat[$\omega^0$]{\includegraphics[width=0.19\linewidth]{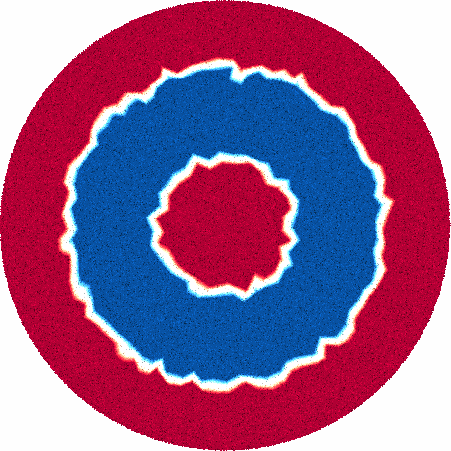}}
		\subfloat[$\alpha^0$]{\includegraphics[width=0.19\linewidth]{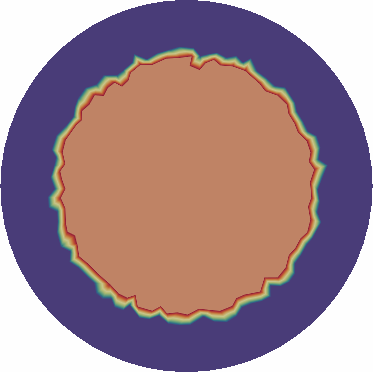}}
		\subfloat[$\omega$]{\includegraphics[width=0.19\linewidth]{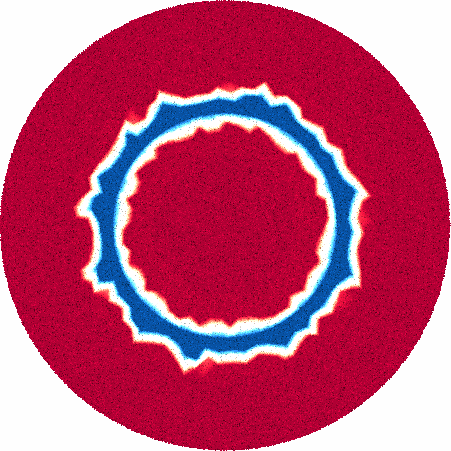}}
		\subfloat[$\alpha$]{\includegraphics[width=0.19\linewidth]{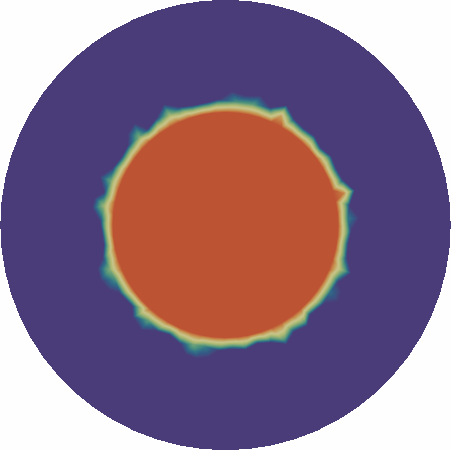}}\\
			{\footnotesize{$\boldsymbol{\omega \equiv 1 \ \ }$ \;  }}
		\subfloat[$\omega^0$]{\includegraphics[width=0.19\linewidth]{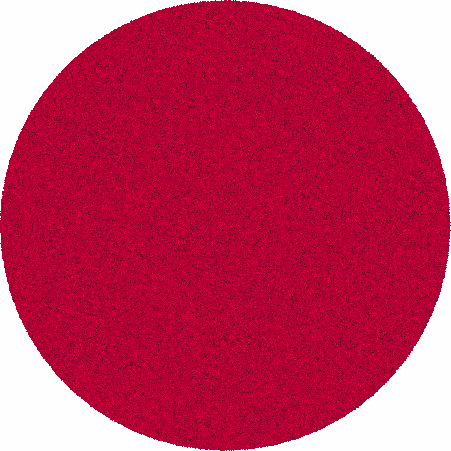}}
 		\subfloat[$\alpha^0$]{\includegraphics[width=0.19\linewidth]{pl02a0.png}}
 		\subfloat[$\omega$]{\includegraphics[width=0.19\linewidth]{Tpw.png}}
 		\subfloat[$\alpha$]{\includegraphics[width=0.19\linewidth]{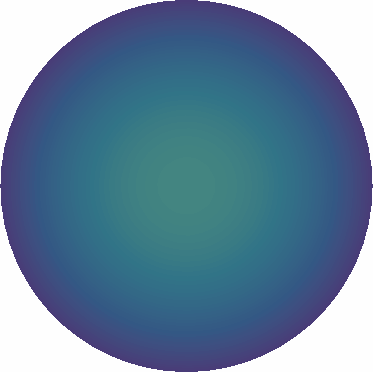}}
		\end{minipage}
\caption{\label{figure:conceverylayer}$BV$ {\it Regularized Inversion} solution pair $(\omega,\alpha)$ at the final iteration for the $\omega^0$ profiles: $(\ell=\lbrace 0.2, 0.3, 0.4\rbrace,\ \omega\equiv1)$ and $\mu=\lbrace 1,0.1,0.1,1 \rbrace$ respectively for the concentric inhomogeneity problem.}
\end{figure}

\begin{figure}[H]\onecolumn
  \centering
 	\hspace{-1.5cm}
	\begin{minipage}{.08\textwidth}
  	\vspace{0.2cm}
		\includegraphics[width=0.71\textwidth,height=0.50\textheight]{colorbar12.png}
 \end{minipage} \hspace{.3cm}
	\centering
\begin{minipage}{.96\textwidth}
\centering
		{\footnotesize{$\boldsymbol{\ell=0.2 \ }$ \; }} \subfloat[$\omega^0$]{\includegraphics[width=0.19\linewidth]{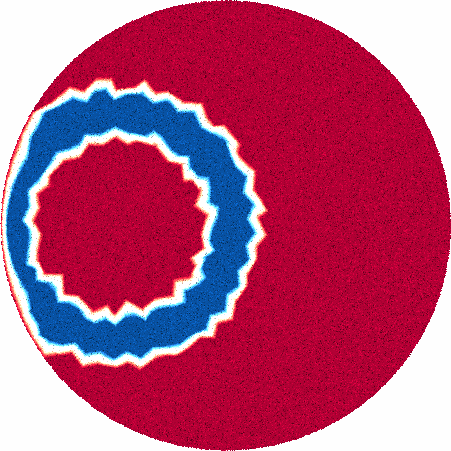}}
		\subfloat[$\alpha^0$]{\includegraphics[width=0.19\linewidth]{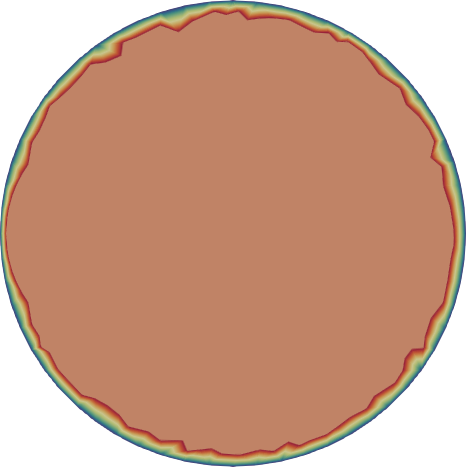}}
		\subfloat[$\omega$]{\includegraphics[width=0.19\linewidth]{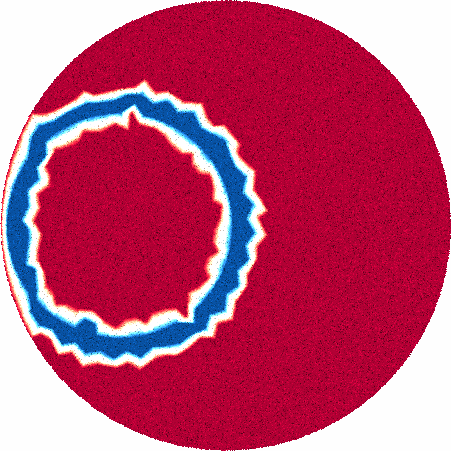}}
		\subfloat[$\alpha$]{\includegraphics[width=0.19\linewidth]{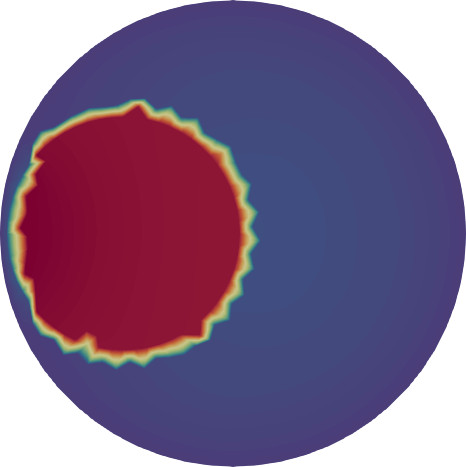}}\\
		{\footnotesize{$\boldsymbol{\ell=0.3 \ }$ \; }} \subfloat[$\omega^0$]{\includegraphics[width=0.19\linewidth]{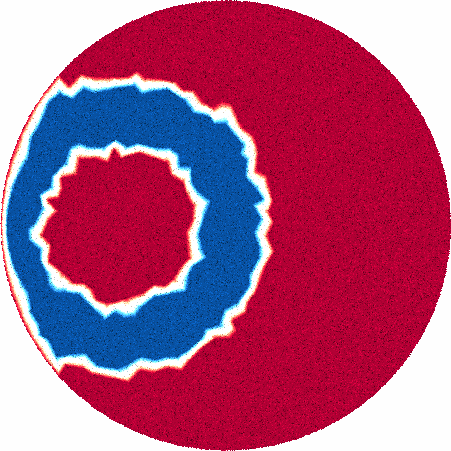}}
		\subfloat[$\alpha^0$]{\includegraphics[width=0.19\linewidth]{sel02a0.png}}
	\subfloat[$\omega$]{\includegraphics[width=0.19\linewidth]{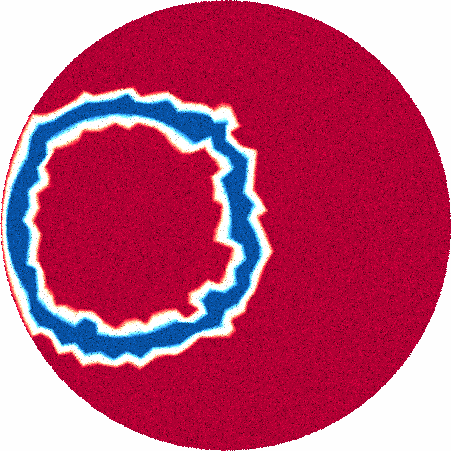}}
	   \subfloat[$\alpha$]{\includegraphics[width=0.19\linewidth]{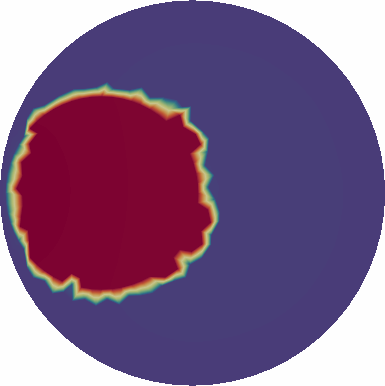}}\\
		{\footnotesize{$\boldsymbol{\ell=0.4 \ }$ \; }} \subfloat[$\omega^0$]{\includegraphics[width=0.19\linewidth]{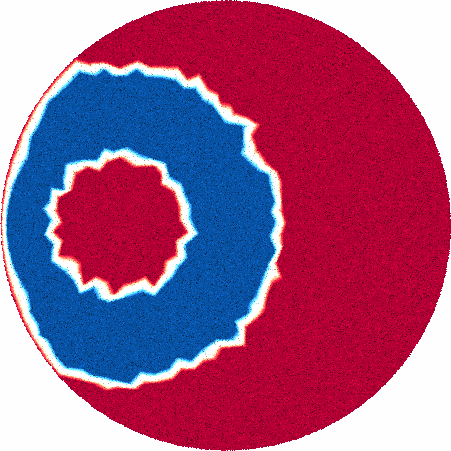}}
		\subfloat[$\alpha^0$]{\includegraphics[width=0.19\linewidth]{sel02a0}}
		\subfloat[$\omega$]{\includegraphics[width=0.19\linewidth]{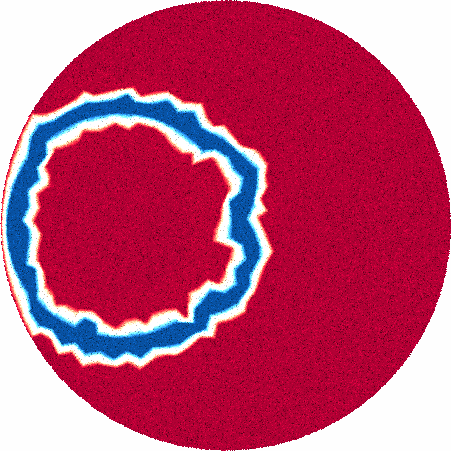}}
		\subfloat[$\alpha$]{\includegraphics[width=0.19\linewidth]{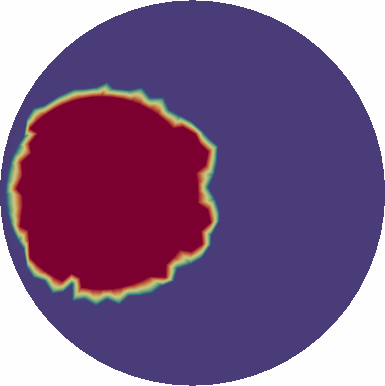}}\\
			{\footnotesize{$\boldsymbol{\omega \equiv 1 \ \ }$ \;  }}
		\subfloat[$\omega^0$]{\includegraphics[width=0.19\linewidth]{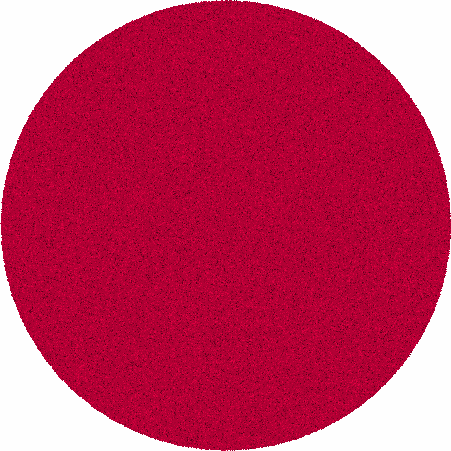}}
		\subfloat[$\alpha^0$]{\includegraphics[width=0.19\linewidth]{sel02a0.png}}
		\subfloat[$\omega$]{\includegraphics[width=0.19\linewidth]{Tsew0.png}}
		\subfloat[$\alpha$]{\includegraphics[width=0.19\linewidth]{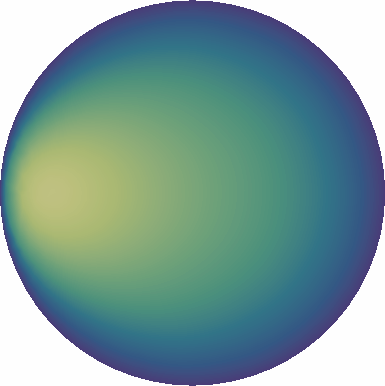}}
		\end{minipage}
\caption{\label{figure:seeverylayer}$BV$ {\it Regularized Inversion} solution pair $(\omega,\alpha)$ at the final iteration for the $\omega^0$ profiles: $(\ell=\lbrace 0.2, 0.3, 0.4\rbrace,\ \omega\equiv1)$ and $\mu=\lbrace 0.1,0.5,0.5,0.1 \rbrace$ respectively for the strong eccentricity problem.}
\end{figure}

\begin{figure}[H]\onecolumn 
  \centering
 	\hspace{-1.5cm}
	\begin{minipage}{.08\textwidth}
  	\vspace{0.2cm}
	\includegraphics[width=0.71\textwidth,height=0.50\textheight]{colorbar12.png}
 \end{minipage} \hspace{.3cm}
	\centering
\begin{minipage}{.96\textwidth}
\centering
		{\footnotesize{$\boldsymbol{\ell=0.2 \ }$ \; }} \subfloat[$\omega^0$]{\includegraphics[width=0.19\linewidth]{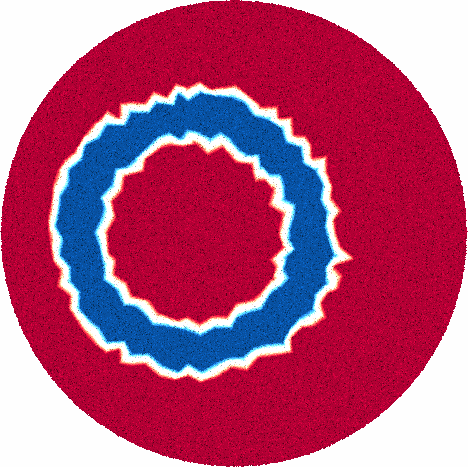}}
		\subfloat[$\alpha^0$]{\includegraphics[width=0.19\linewidth]{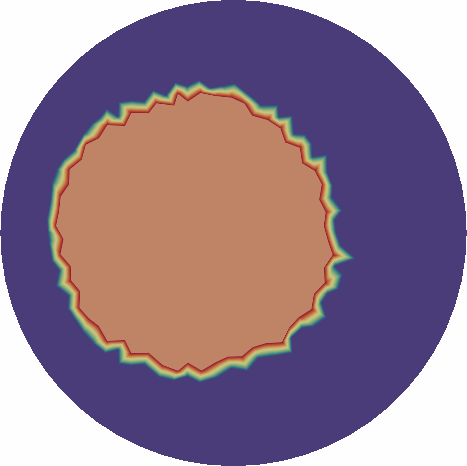}}
		\subfloat[$\omega$]{\includegraphics[width=0.19\linewidth]{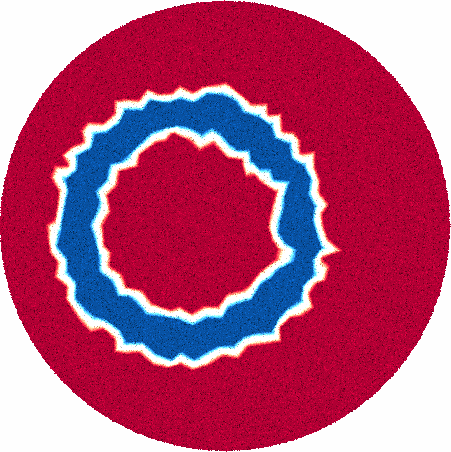}}
		\subfloat[$\alpha$]{\includegraphics[width=0.19\linewidth]{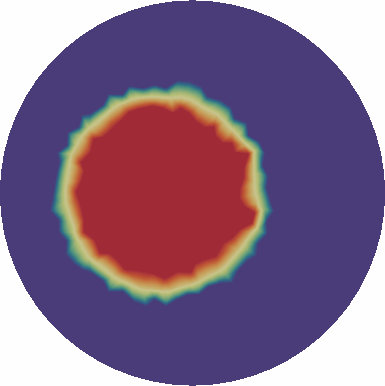}}\\
		{\footnotesize{$\boldsymbol{\ell=0.3 \ }$ \; }} \subfloat[$\omega^0$]{\includegraphics[width=0.19\linewidth]{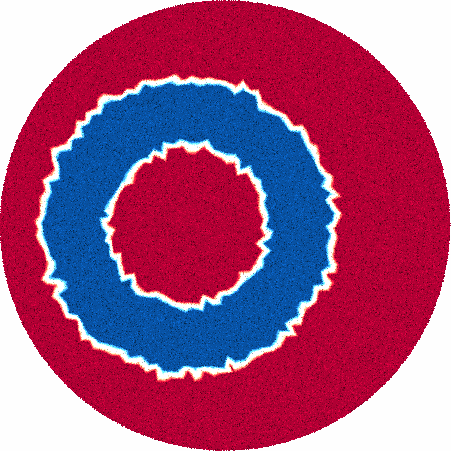}}
		\subfloat[$\alpha^0$]{\includegraphics[width=0.19\linewidth]{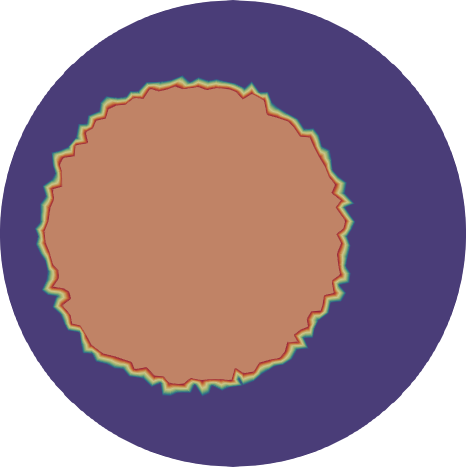}}
		\subfloat[$\omega$]{\includegraphics[width=0.19\linewidth]{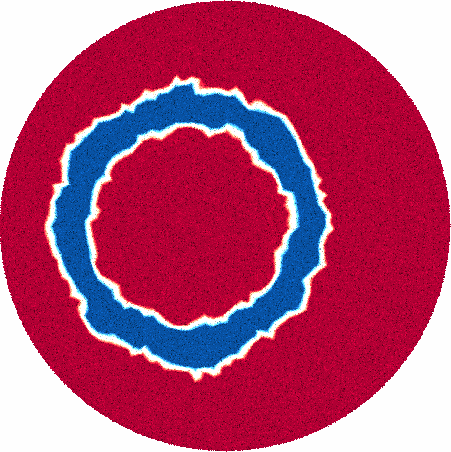}}
		\subfloat[$\alpha$]{\includegraphics[width=0.19\linewidth]{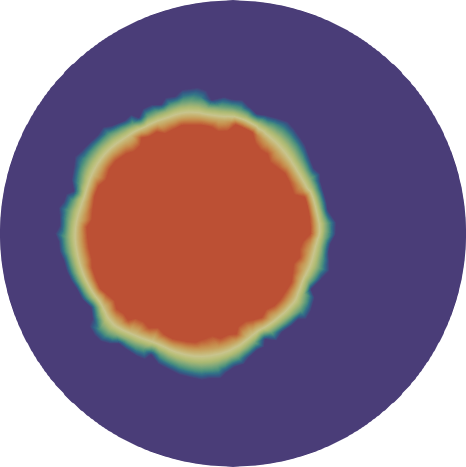}}\\
		{\footnotesize{$\boldsymbol{\ell=0.4 \ }$ \; }} \subfloat[$\omega^0$]{\includegraphics[width=0.19\linewidth]{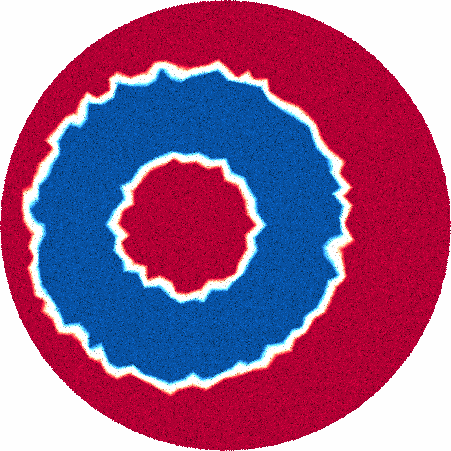}}
		\subfloat[$\alpha^0$]{\includegraphics[width=0.19\linewidth]{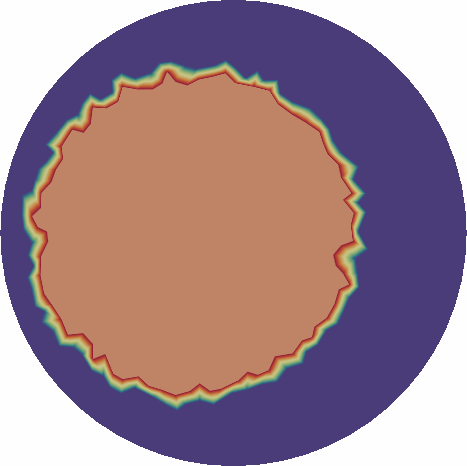}}
		\subfloat[$\omega$]{\includegraphics[width=0.19\linewidth]{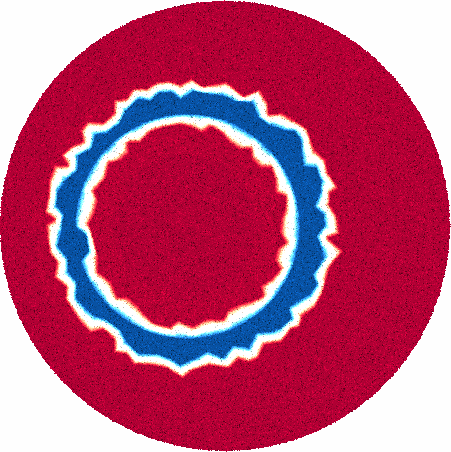}}
		\subfloat[$\alpha$]{\includegraphics[width=0.19\linewidth]{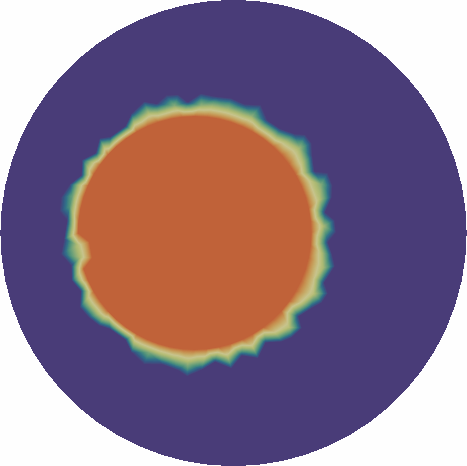}}\\
 			{\footnotesize{$\boldsymbol{\omega \equiv 1 \ \ }$ \;  }}
		\subfloat[$\omega^0$]{\includegraphics[width=0.19\linewidth]{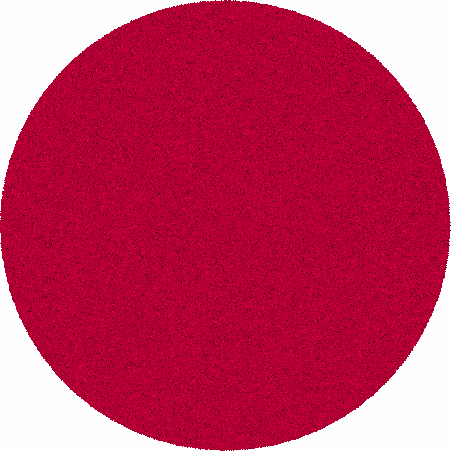}}
		\subfloat[$\alpha^0$]{\includegraphics[width=0.19\linewidth]{mel02a0.png}}
		\subfloat[$\omega$]{\includegraphics[width=0.19\linewidth]{Tmew.png}}
		\subfloat[$\alpha$]{\includegraphics[width=0.19\linewidth]{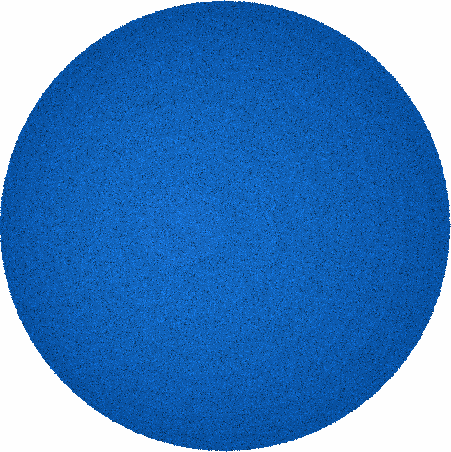}}
		\end{minipage}
\caption{\label{figure:meeverylayer}$BV$ {\it Regularized Inversion} solution pair $(\omega,\alpha)$ at the final iteration for the $\omega^0$ profiles: $(\ell=\lbrace 0.2, 0.3, 0.4\rbrace,\ \omega\equiv1)$ and $\mu=\lbrace 0.1,1,0.1,0.1 \rbrace$ respectively for the mild eccentricity problem.}
\end{figure}

   We emphasize here the sensitivity of the solutions with respect to the two fundamental parameters of the algorithm, $\omega$ and $\mu$. With regard to the first, we introduce the parameter $\ell$ for representing the average width of the ring-shaped $\omega^0$ profile. So, the initial guess for $\omega^0$ is selected by 
\begin{equation}\label{eq:lcontrol}
 \omega^0=\begin{cases}
 0, \qquad \mbox{if}\quad \vert r_i-\rho_{D} \vert \le \frac{\ell}{2} ,\\
 1, \qquad \mbox{otherwise},\end{cases}
\end{equation}
 where $\vert r_i-\rho_{D} \vert$ stands for the Euclidean distance of the $i-$node from $\partial \mathcal{D},\;(\rho_{D}=1)$.
 We consider three different configurations for $\omega^0$ constructed via \ref{eq:lcontrol}, with $\ell=\lbrace 0.2, 0.3, 0.4\rbrace$. Higher values of $\ell$ express higher uncertainty with respect to the location of the inclusion's boundary. The region coloured blue defines the area where $\omega^0(x)=0$, representing so the discontinuity surface (Remark \ref{discandhom}) (see for example the first column of Figure \ref{figure:conceverylayer}, Figure \ref{figure:seeverylayer}, Figure \ref{figure:meeverylayer}). We also consider the case where the initial $\omega$ is a constant function, $\omega^0\equiv 1$.  
	\begin{table}[!t]
\centering
		\caption{\label{table:infomesh}Delaunay Triangulations.} 
	\begin{tabular}{@{}lccc}
		\bottomrule
		 Mesh case & $h$ & Elements & Nodes\\
		\midrule
		 {\it concentric} & $0.27$ & $1850$ & $956$ \\
	 {\it strong eccentric} & $0.24$ & $1834$ & $948$\\
    {\it mild eccentric} & $0.26$ & $1854$ & $958$ \\
    \midrule
{\it concentric} & $0.15$ &  $2024$ & $3926$ \\
 {\it strong eccentric} & $0.14$ & $2036$ & $3950$ \\
 {\it mild eccentric} & $0.15$ & $2036$ & $3950$ \\
 \bottomrule
\end{tabular}	
\end{table}
\newline The regularization parameter $\mu$ represents the influence of the $BV$ structure on the problem (although slightly regularized by the parameter $\epsilon$). It is surprising that only with one measurement, we are able to reconstruct both the geometrical configuration and the exact values of the conductivity profile, even though the initial guess (expressed via the configuration $\omega^0$) protects the minimization process from rambling. Besides, our aim is to examine numerically the robustness of the proposed methodology with respect to the choice of $\omega^0$. \newline
The minimization reveals that when $\mu$ increases, the $BV$ norm of the conductivity dominates and highlights the discontinuities of the conductivity function. It is interesting to mention that when we impose constant $\omega$,  (for example $\omega^0 \equiv 1$), the associated term of the functional becomes $\mu \int_{\Omega} \vert \nabla\alpha\vert^2$ and we recover the commonly employed Tikhonov regularization, forcing the conductivity profiles to be $H^1-$functions, preventing the abruptness of the variations to be revealed across inclusion boundaries. The smoothing of the inclusion boundary, observed for example at the last row (p) of Figure \ref{figure:conceverylayer}, is indicative for the inexpediency of the common Tikhonov regularization $(\omega \equiv 1)$, revealing convergence to false local minima.

 We proceed with the implementation of the algorithm for the numerical solution of the problem (\ref{problemtoreconstruct}) setting the smoothing parameter $\epsilon=0.1$ and the parameter $\kappa=10$. The initial approximation for the conductivity is chosen to be the constant function $\alpha^0=2.5$. Namely, in the concentric (\ref{subsubsec:CONC}) and mild eccentric inhomogeneity (\ref{EC1}), we set 
\begin{equation}\label{eq:a0control}
 \alpha^0 =\begin{cases} 
     1, \qquad &r_i>\rho_D+\frac{\ell}{2},\\
     2.5, \qquad &\mbox{otherwise},     
       \end{cases}
\end{equation}
 (second column of Figure \ref{figure:conceverylayer} and Figure \ref{figure:meeverylayer}), while in the strong eccentricity inhomogeneity (\ref{EC2}) we establish $\alpha_0=2.5$ throughout the domain $\Omega\setminus\Omega_{\delta}$ (second column of Figure \ref{figure:seeverylayer}). The inversion is attempted for each aforementioned value of $\ell$ as well as the case $\omega^0\equiv1$ for $n=10$ iterations. The optimal solutions $\alpha^n$ are depicted for the best - for each case - value of the regularization parameter $\mu$ selected from the value range $\mu=\lbrace 0.1, 0.5, 1, 5 \rbrace$. We clarify here, that the ``best'' choice of the $BV$ regularization term consists in that value of $\mu$, which improves the solution $\alpha^n$ to yield a better result both in the geometrical and physical sense (Table \ref{table:infoP3Conc}, Table \ref{table:infoP3SE}, Table \ref{table:infoP3ME}). We display the solution pair $(\alpha,\omega)$ of the final iteration ($n=10$) for $\ell=\lbrace 0.2,0.3,0.4\rbrace$ and we choose only one value of $\ell$ to present the $(\alpha^n,\omega^n)$ series in each geometry ( Figure \ref{figure:concl02series}, Figure \ref{figure:sel04series}, Figure \ref{figure:mel03series}). In the latter display, we provide also the conductivity contour lines (black lines) along with their values, in order to highlight the uniformity of the physics in the conductivity reconstruction. Finally, we point out that the Tikhonov regularization corresponds to the case $\omega^0\equiv1$ (consequently $\omega\equiv1$). Hereafter, we use the notation $\omega\equiv1$ accounting for the $\omega^0-$ constant profile of the Tikhonov reconstruction. 

\subsubsection{Concentric inhomogeneity}\label{subsubsec:CONC} 
We consider the simple case where the inclusion $D$ is located at $(x_D,y_D)=(0,0)$.
 We additionally consider the synthetic data
\begin{equation}
\begin{aligned}
& f= 1 + \frac{11}{4} \cos(\phi),\quad & 0\leqslant\phi < 2\pi, \\ & g= \frac{13}{8} \cos(\phi),\quad & 0\leqslant\phi < 2\pi,
\end{aligned}\label{pilotfg} 
\end{equation}
produced by the analytical solution
\begin{equation}
 u(\rho, \phi) = \begin{cases}
 1 +\rho \cos (\phi), \quad & 0\leq\rho \le  1, \\
 1 + (\frac{3}{2} \rho  - \frac{1}{2 \rho}) \cos(\phi), \quad & 1 < \rho \le 2.\end{cases}\label{eq:AnalyticalProblem1}
\end{equation}

\begin{figure}[H]\onecolumn
	\centering
  	\hspace{-1.5cm}
	\begin{minipage}{.08\textwidth}
  	\vspace{-1.5cm}
	\includegraphics[width=.45\textwidth,height=0.28\textheight]{colorbar12.png}
 \end{minipage} \hspace{.35cm} 
	\centering
\begin{minipage}{.96\textwidth}
$\alpha^n$\subfloat[] {	\includegraphics[width=0.19\linewidth]{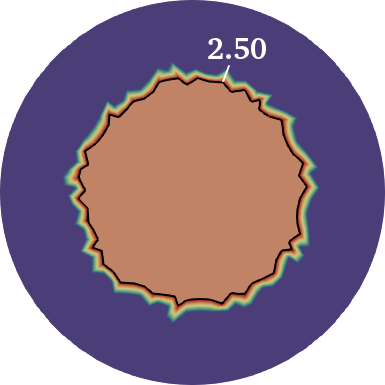}
}
 \subfloat[] {   			\includegraphics[width=0.19\linewidth]{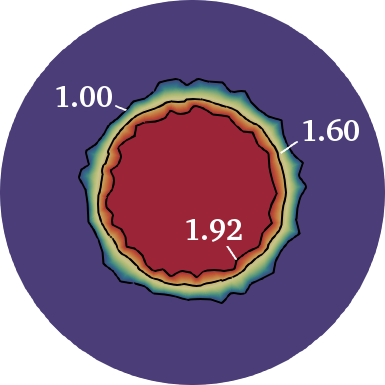}	}             		  	
  \subfloat[] { 			\includegraphics[width=0.19\linewidth]{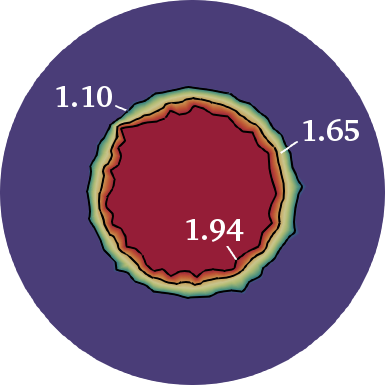}
    		} 
 \subfloat[] {     			\includegraphics[width=0.19\linewidth]{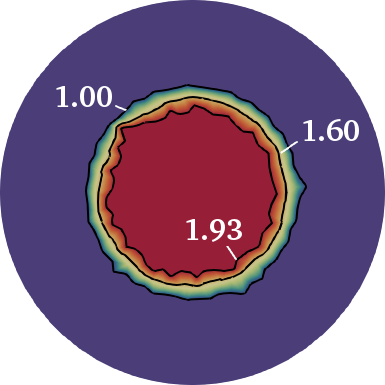}}  
  \subfloat[] {       			\includegraphics[width=0.19\linewidth]{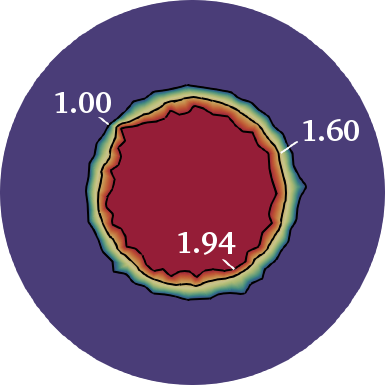}}  \\
 
  $\omega^n$
  \subfloat[$n=0$]  {
  \includegraphics[width=0.19\linewidth]{pl02w0.png}  		} \subfloat[$n=1$] {             			\includegraphics[width=0.19\linewidth]{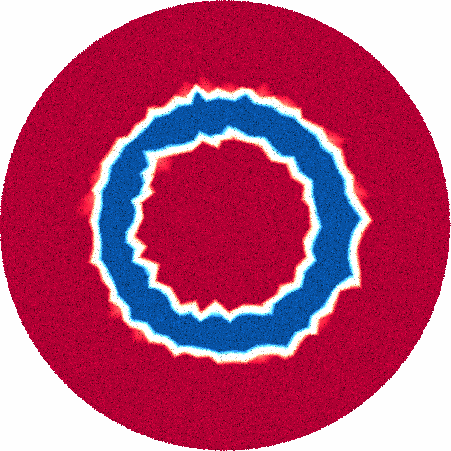}
             		}             		  	
  \subfloat[$n=2$] {
 \includegraphics[width=0.19\linewidth]{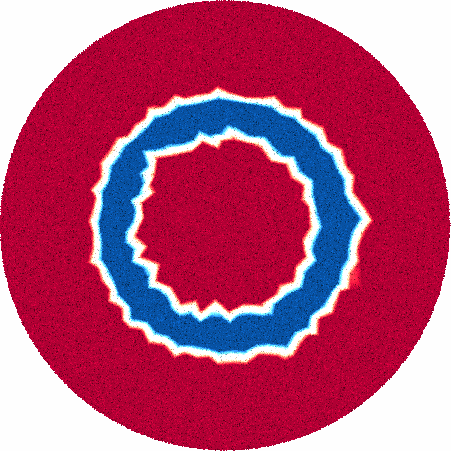}
             		}
  \subfloat[$n=3$] {
  	\includegraphics[width=0.19\linewidth]{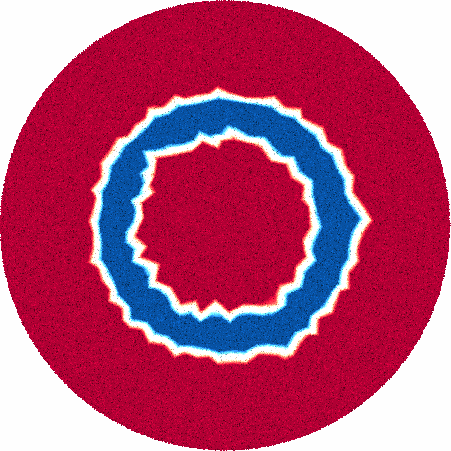}
                 }
\subfloat[$n=6$] {
  	\includegraphics[width=0.19\linewidth]{pl02wend.png}
                 }                 \\
\caption{\label{figure:concl02series}$BV$ {\it Regularized Inversion} $\alpha^n-$ solutions and $\omega^n-$ profiles at the $n-$th iteration for $\ell=0.2$ of the concentric inhomogeneity problem.}
  \end{minipage}      
\end{figure}

\noindent Using the computational process which was described in Section \ref{subsec:Algorithm}, in Figure \ref{figure:conceverylayer} from left to right we display the initial $(\omega^0, \alpha^0)$ and the final corresponding numerical solution $(\omega, \alpha)$ for the $\omega^0-$ configurations: $\ell=\lbrace 0.2,0.3,0.4\rbrace$ and $\omega^0\equiv1$ (last row). The inversion is performed for the different values of $\mu=\lbrace 1, 0.1, 0.1, 1\rbrace$ for $\ell=\lbrace 0.2,0.3,0.4\rbrace,\ \omega^0\equiv 1$ respectively and mesh size $h=0.27$. In Figure \ref{figure:concl02series},  we present extensively the series of the optimal solutions $\alpha^n$ (1st row) corresponding to the profiles $\omega^n$ (2nd row) for the chosen $\omega^0-$ configuration $\ell=0.2$. We summarize the numerical results for all iterations in Table \ref{table:infoP3Conc}. \newline
All tables hereafter show the $\omega^0-$ refinement level $(\ell)$ along with the selected regularization value $(\mu)$, the value of the computed conductivity in the inclusion $D$, $(\alpha_{in})$, and that of the outer domain $\Omega\setminus D$, $(\alpha_{out})$, in every iteration $n$. The variables $(\alpha_{in, out})$ denote the uniform values of the reconstructed conductivity $\alpha$.
We  point out here for readers' convenience, that the color-bar in all figures pertains to the conductivity range ($\alpha\in[1,2]$).

 \subsubsection{Eccentric inhomogeneity}\label{subsubsec:EEC}
 \paragraph{Strong eccentricity problem.}\label{EC2}
The center of an indicative strongly eccentric inclusion is chosen to be at $(x_{D},y_{D})=(\frac{\sqrt{5}-\sqrt{17}}{2},0)$. We assume that synthetic data are available in the form
\begin{equation}
\begin{aligned}
& f=s_f(\alpha_2,\tau_1,\tau_2)\frac{\sin\phi}{\sqrt{17}+4\cos\phi}, \\ 
& g=s_g(\alpha_2,\tau_1,\tau_2)\frac{1}{2}\frac{\sin\phi}{(\sqrt{17}+4\cos\phi)^2},\label{extremeec1}
\end{aligned}
\end{equation}
where 
\begin{equation}
\begin{aligned}
&s_f(\alpha_2,\tau_1,\tau_2)=\frac{1}{2}(1+\alpha_2) e^{-\tau_1}+\frac{1}{2}(1-\alpha_2) e^{-2\tau_2} e^{\tau_1}, \\
&s_g(\alpha_2,\tau_1,\tau_2)= \frac{1}{2}(1+\alpha_2) e^{-\tau_1}-\frac{1}{2}(1-\alpha_2) e^{-2\tau_2} e^{\tau_1}, 
\end{aligned}
\end{equation}
and $\phi \in [0,2\pi)$, generated by the exact solution (obtained by an analytical method in bipolar coordinates)
\begin{align}\label{extremeec2}
& u_{exact}=\nonumber\\ &\begin{cases}
 (\rho^{4}+c_1(\phi)\rho^3+\rho^2c_2(\phi)+c_3(\phi)\rho+16)^{-\frac{1}{2}}e^{-\tau} \rho \sin\phi,\\
\hspace{3cm}\mbox{if}\ \tau > \tau_2,\\
  (\rho^{4}+c_1(\phi)\rho^3 + c_2(\phi)\rho^2+c_3(\phi)\rho +16)^{-\frac{1}{2}}\tilde{c} \rho \sin\phi,
  \\  \hspace{3cm}\mbox{if}\ \tau_1<\tau<\tau_2, 
 \end{cases}
  \end{align}
 with
 \begin{align}
 &a=\frac{1}{2},\ \alpha_2=2,\ \tau_1=\ln[(\sqrt{17}+1)\frac{1}{4}],\nonumber \\  &\tau_2=\ln[(\sqrt{5}+1)\frac{1}{2}],\  c_1(\phi)=2 \sqrt{17}\cos\phi,\nonumber \\ &c_2(\phi)=17+8\cos(2\phi), \ c_3(\phi)=8\sqrt{17}\cos\phi, \nonumber \\   &\tilde{c}=\frac{1}{2}(1+\alpha_2) e^{-\tau}+\frac{1}{2}(1-\alpha_2) e^{-2\tau_2} e^{\tau}, \nonumber\\ &\tau=\frac{1}{2}\ln\Bigg[ \frac{\rho^2+\frac{4a^2}{(1-e^{-2\tau_1})^2}+\frac{4a\rho \cos\phi}{1-e^{-2\tau_1}}}{\rho^2+\frac{4a^2e^{-4\tau_1}}{(1-e^{-2\tau_1})^2}+\frac{4a\rho \cos\phi e^{-2\tau_1}}{1-e^{-2\tau_1}}}\Bigg].
    \end{align}

\begin{figure}[H]\onecolumn
	\centering
  	\hspace{-1.5cm}
	\begin{minipage}{.08\textwidth}
  	\vspace{-1.5cm}
	\includegraphics[width=.45\textwidth,height=0.28\textheight]{colorbar12.png}
 \end{minipage} \hspace{.35cm} 
	\centering
\begin{minipage}{.96\textwidth}
$\alpha^n$\subfloat[] {       
\includegraphics[width=0.19\linewidth]{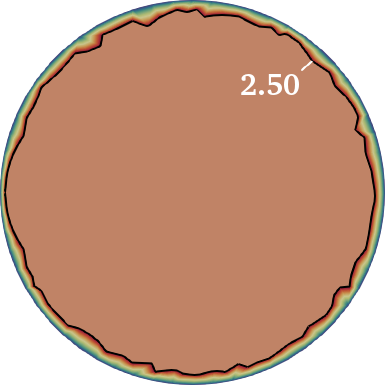}  
}
 \subfloat[] {             			\includegraphics[width=0.19\linewidth]{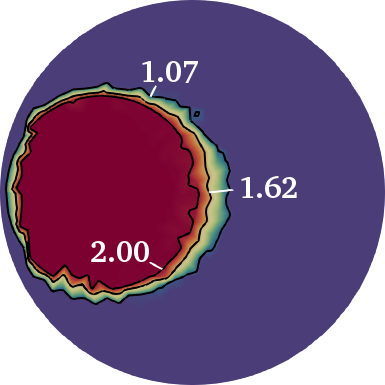}	}             		  	
  \subfloat[] {             			\includegraphics[width=0.19\linewidth]{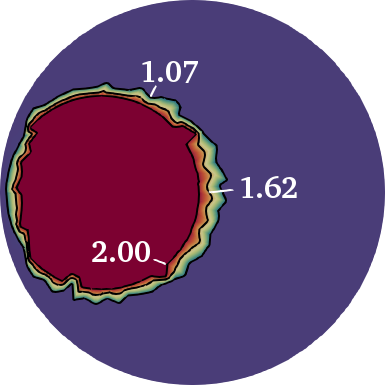}
    		} 
 \subfloat[] {           			\includegraphics[width=0.19\linewidth]{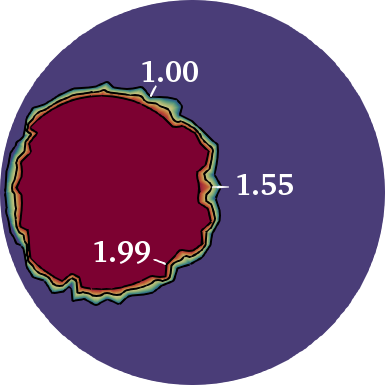}}  
  \subfloat[] {           			\includegraphics[width=0.19\linewidth]{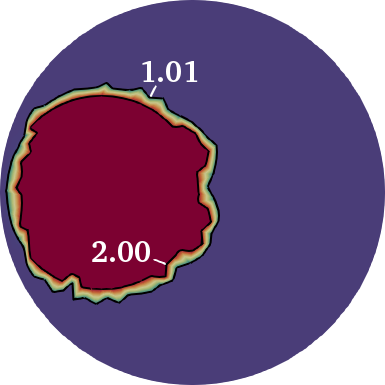}}  \\
 
  $\omega^n$
  \subfloat[$n=0$]  {
  \includegraphics[width=0.19\linewidth]{sel04w0.png}  		} \subfloat[$n=1$] {             			\includegraphics[width=0.19\linewidth]{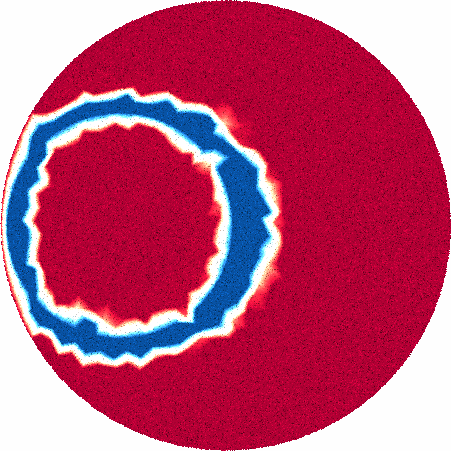}
             		}             		  	
  \subfloat[$n=2$] {
 \includegraphics[width=0.19\linewidth]{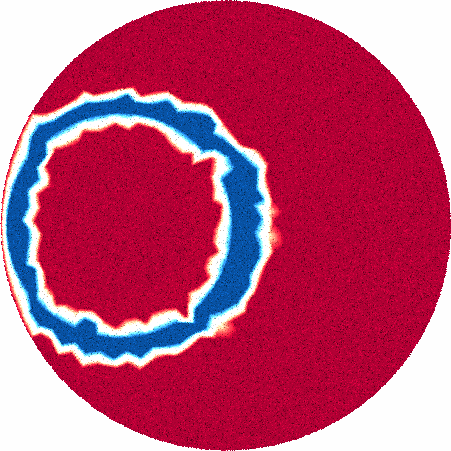}
             		}
  \subfloat[$n=4$] {
  	\includegraphics[width=0.19\linewidth]{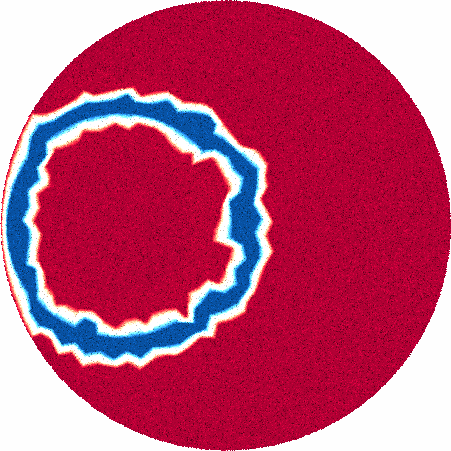}
                 }
\subfloat[$n=8$] {
  	\includegraphics[width=0.19\linewidth]{sel04iter10w.png}
                 }                 \\
	\caption{\label{figure:sel04series} $BV$ {\it Regularized Inversion} $\alpha^n-$ solutions and $\omega^n-$ profiles at the $n-$th iteration for $\ell=0.4$ of the strong eccentricity problem.}
  \end{minipage}      
	\end{figure}

 \noindent In our computational process, we have observed that the regularization values $\mu=\lbrace 0.1,0.5,0.5, 0.1\rbrace$, corresponding to $\ell=\lbrace 0.2,0.3,0.4\rbrace$ and $\omega^0\equiv1$, reveal appropriately the geometrical features of that case. We present the numerical results with respect to these values of $\mu$ and mesh size $h=0.24$. \newline  
Figure \ref{figure:seeverylayer} from left to right displays the chosen initial guess $(\omega^0,\alpha^0)$ along with the numerical solution pair $(\omega,\alpha)$ of the final iteration. Next in Figure \ref{figure:sel04series}, we show the sequence of $(\alpha^n, \omega^n)$ for the most disturbed $\omega^0-$ configuration $\ell=0.4$. We perform the numerical results in Table \ref{table:infoP3SE}. 

\paragraph{Mild eccentricity problem.}\label{EC1}

In this exa\-mple, the coordinates of the center of $D$ are taken as $(x_{D},y_{D})=(-\frac{1}{3},0)$.
We handle the following pair of synthetic data
\begin{equation}
\begin{aligned}
&f=m_f(\alpha_2,\tau_1,\tau_2)\frac{\frac{8}{3}\sqrt{10}\sin\phi}{\frac{28}{3}+4\cos\phi},  \\
&g=m_g(\alpha_2,\tau_1,\tau_2)\frac{320}{9}\frac{\sin\phi}{(\frac{28}{3}+4\cos\phi)^2},
\end{aligned}\label{eq:ecn1}
\end{equation}
where
\begin{align}
&m_f(\alpha_2,\tau_1,\tau_2)= \frac{1}{2}(1+\alpha_2) e^{-\tau_1}+\frac{1}{2}(1-\alpha_2) e^{-2\tau_2} e^{\tau_1},\nonumber \\
&m_g(\alpha_2,\tau_1,\tau_2)=\frac{1}{2}(1+\alpha_2) e^{-\tau_1}-\frac{1}{2}(1-\alpha_2) e^{-2\tau_2} e^{\tau_1},
\end{align} and $\phi \in [0,2\pi)$, stemmed from the exact solution
\begin{align}\label{eq:ecn2}
u_{exact}=\begin{cases}
 e^{-\tau}\frac{\sin\phi \sinh\tau}{\cosh\tau+\cos\phi}, \qquad \mbox{if} \;\tau > \tau_2,\\
 m_u(\alpha_2,\tau_2,\tau)\frac{\sin\phi \sinh\tau}{\cosh\tau+\cos\phi},\\  \hspace{3.1cm}\mbox{if}\;\tau_1<\tau<\tau_2, \end{cases}
\end{align}
with $$m_u(\alpha_2,\tau_2,\tau)=\Big[ \frac{1}{2}(1+\alpha_2) e^{-\tau}+\frac{1}{2}(1-\alpha_2) e^{-2\tau_2} e^{\tau}\Big]$$
and
\begin{align}
&a=\frac{4}{3} \sqrt{10},\ \alpha_2=2,\ \tau_1=\ln[\frac{(2\sqrt{10}+7)}{3}],\nonumber\\ &\tau_2=\ln[\frac{(4\sqrt{10}+13)}{3}],\nonumber\\
&\tau=\frac{1}{2}\ln\Bigg[ \frac{\rho^2+\frac{4a^2}{(1-e^{-2\tau_1})^2}+\frac{4a\rho \cos\phi}{1-e^{-2\tau_1}}}{\rho^2+\frac{4a^2e^{-4\tau_1}}{(1-e^{-2\tau_1})^2}+\frac{4a\rho \cos\phi e^{-2\tau_1}}{1-e^{-2\tau_1}}}\Bigg].
\end{align}

\begin{figure}[H]\onecolumn
	\centering
  	\hspace{-1.5cm}
	\begin{minipage}{.08\textwidth}
  	\vspace{-1.5cm}
	\includegraphics[width=.45\textwidth,height=0.28\textheight]{colorbar12.png}
 \end{minipage} \hspace{.35cm} 
	\centering
\begin{minipage}{.96\textwidth}
$\alpha^n$\subfloat[] {             			\includegraphics[width=0.19\linewidth]{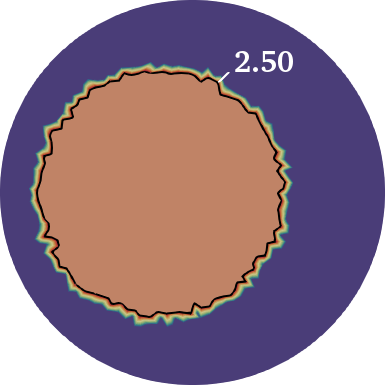}
}
 \subfloat[] {             			\includegraphics[width=0.19\linewidth]{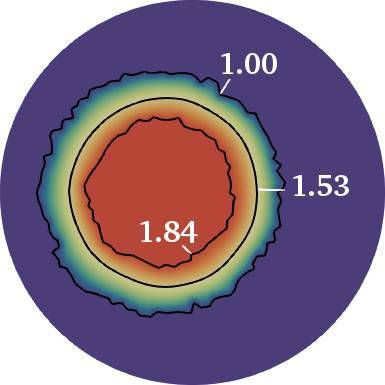}
 }             		  	
  \subfloat[] {             			\includegraphics[width=0.19\linewidth]{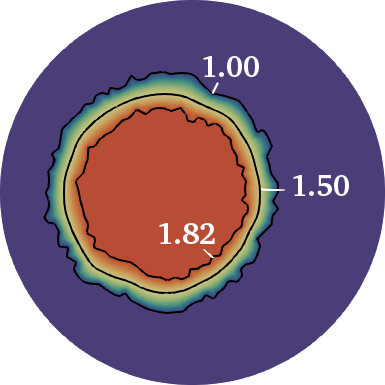}
    		} 
 \subfloat[] {           			\includegraphics[width=0.19\linewidth]{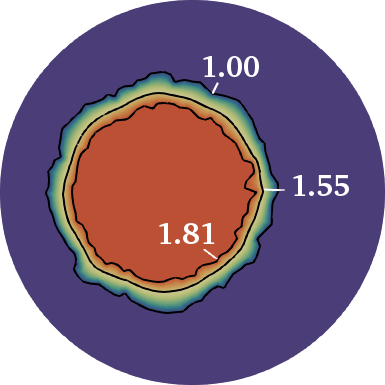}}  
  \subfloat[] {           			\includegraphics[width=0.19\linewidth]{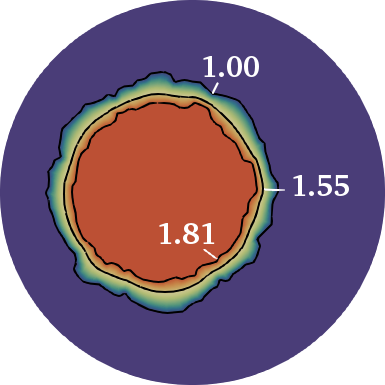}}  \\
 
  $\omega^n$
  \subfloat[$n=0$]  {
  \includegraphics[width=0.19\linewidth]{mel03w0.png}  	
  }
  \subfloat[$n=1$] {             			\includegraphics[width=0.19\linewidth]{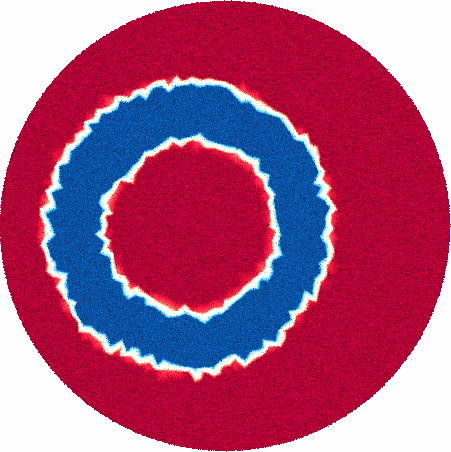}
             		}             		  	
  \subfloat[$n=3$] {
 \includegraphics[width=0.19\linewidth]{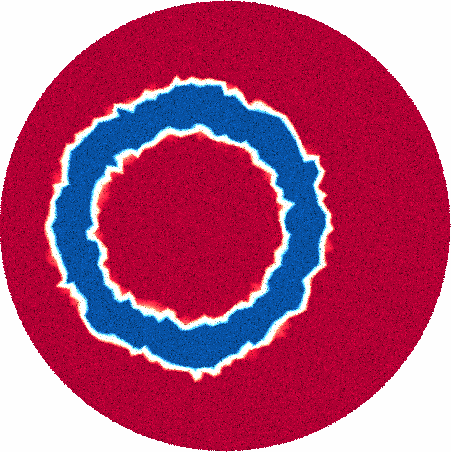}
             		}
  \subfloat[$n=5$] {
  	\includegraphics[width=0.19\linewidth]{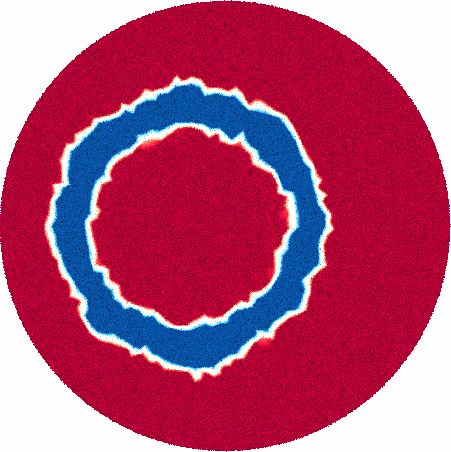}
                 }
\subfloat[$n=10$] {
  	\includegraphics[width=0.19\linewidth]{mel03wend.png}
                 }                 \\
	\caption{\label{figure:mel03series}$BV$ {\it Regularized Inversion} $\alpha^n-$ solutions and $\omega^n-$ profiles at the $n-$th iteration for $\ell=0.3$ of the mild eccentricity problem.}
  \end{minipage}      
	\end{figure}

 \noindent In Figure \ref{figure:meeverylayer} from left to right we plot the initial pair $(\omega^0,\alpha^0)$ and its computed solution pair $(\omega,\alpha)$ on the tenth iteration, for the optimal regularization values - observed in computations - $\mu=\lbrace 0.1,1,0.1,0.1\rbrace$ for $\ell=\lbrace 0.2, 0.3, 0.4\rbrace$ and $\omega^0\equiv1$ apiece, with respect to the mesh size $h=0.26$. Figure \ref{figure:mel03series} presents the $(\omega^n, \alpha^n)-$ reconstruction sequence generated by the value $\ell=0.3$ for the finer mesh size $h=0.15$. The reconstruction parameters are given in Table \ref{table:infoP3ME}.

\subsection{$BV$ Physical reconstruction\label{subsec:physicalrecon}}
\begin{figure}[H]\onecolumn
	\centering
  	\hspace{-1.5cm}
	\begin{minipage}{.07\textwidth}
  	\vspace{-.2cm}
\includegraphics[width=0.54\textwidth,height=0.32\textheight]{colorbar12.png}
 \end{minipage} \hspace{.35cm} 
	\centering
	\begin{minipage}{.97\textwidth}
	 \subfloat[$\alpha^0$]{\includegraphics[width=0.19\linewidth]{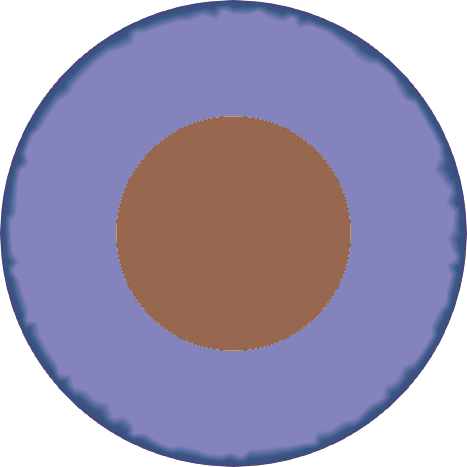}
	 }
		\subfloat[$\omega$]{
		\includegraphics[width=0.19\linewidth]{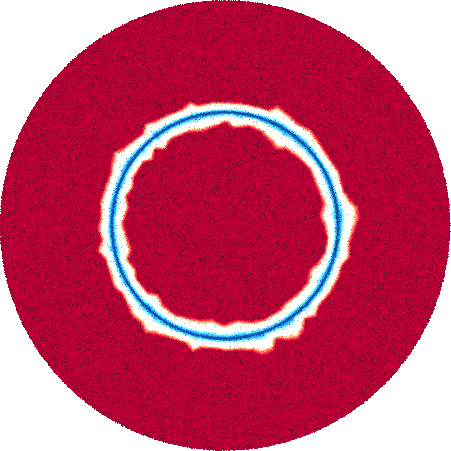}
}
\subfloat[$\alpha$]{
\includegraphics[width=0.19\linewidth]{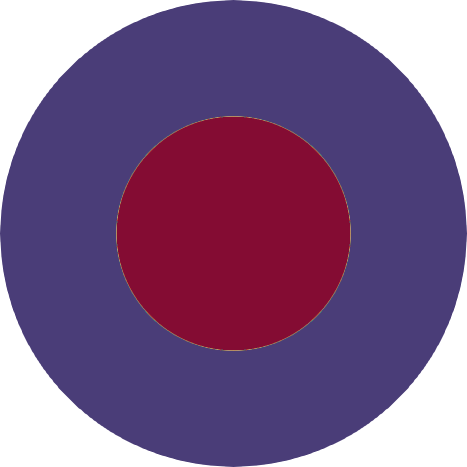}
}\
		\tikz{\draw [dashed] (0,0) -- (0,2.65);}
		\subfloat[$\omega$]{
		\includegraphics[width=0.19\linewidth]{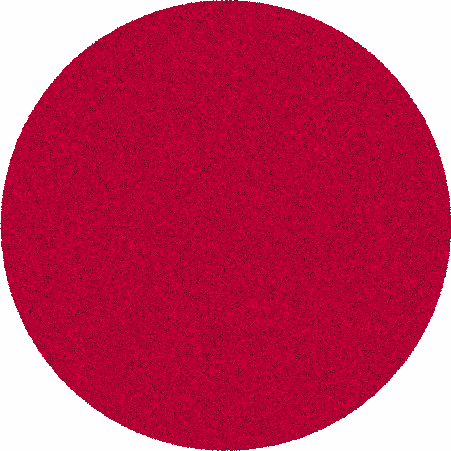}
		}
		\subfloat[$\alpha$]{
		\includegraphics[width=0.19\linewidth]{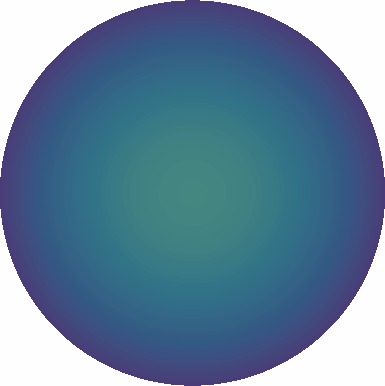}
		}\\
\subfloat[$\alpha^0$]{\includegraphics[width=0.19\linewidth]{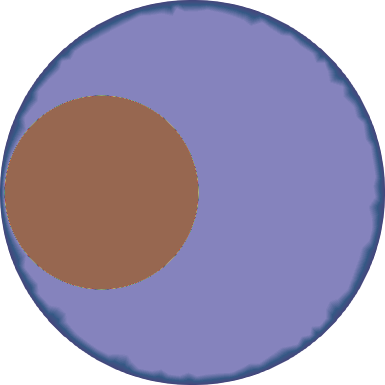} }
		\subfloat[$\omega$]{
		\includegraphics[width=0.19\linewidth]{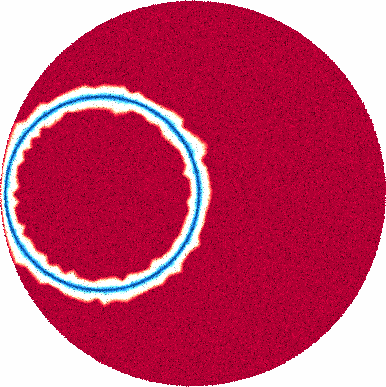}
		}
		\subfloat[$\alpha$]{
		\includegraphics[width=0.19\linewidth]{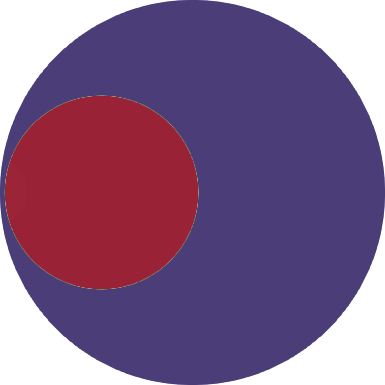}
		}\
		\tikz{\draw [dashed] (0,0) -- (0,2.65);}
		\subfloat[$\omega$]{
		\includegraphics[width=0.19\linewidth]{Tphysicspilotw.png}
		}
		\subfloat[$\alpha$]{
		\includegraphics[width=0.19\linewidth]{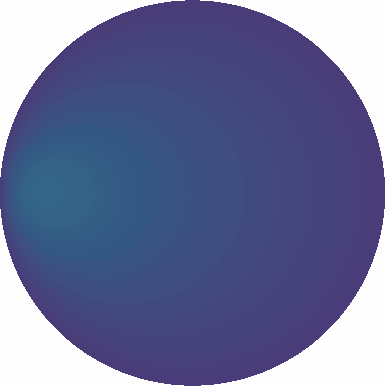}
		}\\
		\subfloat[$\alpha^0$]{\includegraphics[width=0.19\linewidth]{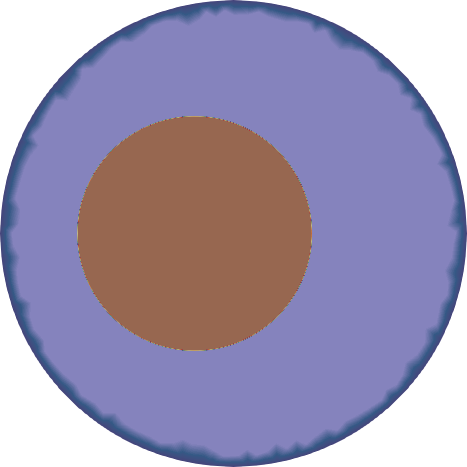}		}
		\subfloat[$\omega$]{
		\includegraphics[width=0.19\linewidth]{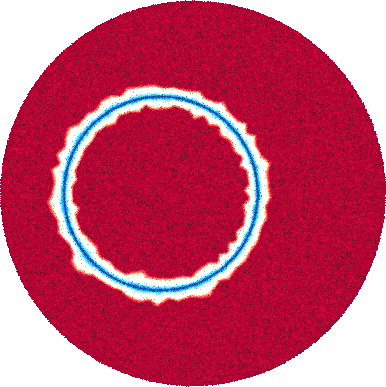}
		}
		\subfloat[$\alpha$]{
		\includegraphics[width=0.19\linewidth]{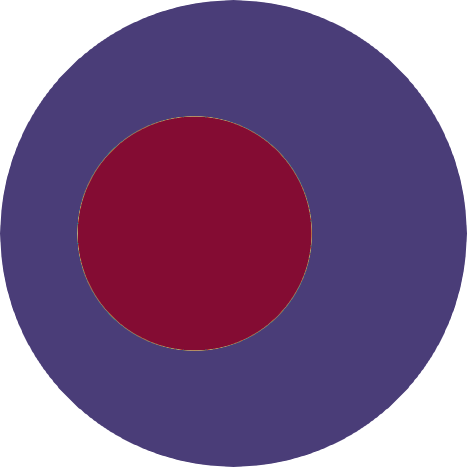}
		}\
		\tikz{\draw [dashed] (0,0) -- (0,2.65);}
		\subfloat[$\omega$]{
		\includegraphics[width=0.19\linewidth]{Tphysicspilotw.png}
		}
		\subfloat[$\alpha$]{
		\includegraphics[width=0.19\linewidth]{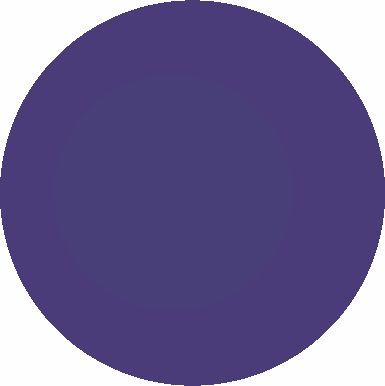}
		}
  \end{minipage}     
  \caption{Single-iteration $BV$ {\it Regularized Inversion} $\alpha-$ solutions  for the $\omega^0$ profiles: ($\ell=0.02,\ \omega\equiv1)$ of every inhomogeneity problem and $\mu=1$. \label{figure:physicalreconalltogether} }
  \end{figure}

In this paragraph, we would like to present the efficiency of the method in case we know a priori the correct geometry and we are aiming at the values of the conductivity function throughout the well defined components of the whole structure. There exist two reasons why we present this analysis. First of all, it is really amazing to remark the efficiency as well as the abrupt convergence of the methodology to the desired exact solution with high accuracy. Secondly, the a priori information of the correct geometry of the inhomogeneity is introduced in the reconstruction technique via the specific form of the field $\omega$ and this implication is strongly related to the efficiency of the numerical treatment of the method, on the basis of the hidden $BV-$structure carried by $\omega$ itself. 
\newline In what follows, we assume that $\omega$ has indeed the expected behavior in the vicinity of the inhomogeneity, as depicted in second column of Figure \ref{figure:physicalreconalltogether}, generated by the value $\ell=0.02$ in formula (\ref{eq:lcontrol}). For all geometries (\ref{subsubsec:CONC}, \ref{EC2}, \ref{EC1}), we set $\epsilon=0.1,\ \kappa=10$, the regularization parameter $\mu=1$ and as initial approximation for the conductivity the three - valued constant function 
\begin{align}
\alpha^0 =\begin{cases} 
    5, \qquad &(x,y) \in D,  \\  0.5, \qquad &(x,y) \in \Omega \setminus (\Omega_{\delta} \cup D),  \\  1, \qquad &(x,y) \in \Omega_{\delta} ,  \end{cases} \label{alpha0physics}
\end{align}
(Figure \ref{figure:physicalreconalltogether}, 1st column). Here, we have selected a very small lower and a very large upper threshold value to show that even under an extended deviation range, the conductivity is recovered immediately and exactly. \newline The $BV$ {\it Regularized Inversion} is performed for one i\-te\-ration with respect to the mesh sizes ($h$) displayed in Table \ref{table:infomeshphys} following. We also consider for comparison the case $\omega^0\equiv 1$, accounting for the Tikhonov reconstru\-ction.

\begin{table}[ht!]
\centering
		\caption{\label{table:infomeshphys}Delaunay Triangulations. } 
		\begin{tabular}{@{}lccc}
		\bottomrule
		 Mesh case & $h$ & Elements & Nodes\\
		\midrule
{\it concentric} & $0.21$ & $2698$ & $1450$ \\
{\it strong eccentric} & $0.21$ & $2738$ & $1470$\\ 
{\it mild eccentric} & $0.20$ & $2746$ & $1474$ \\
	\bottomrule
\end{tabular}	
\end{table}
\noindent
 In Figure \ref{figure:physicalreconalltogether} from left to right we display the initial guess $(\alpha^0)$, the ``right'' $(\omega)$ and the computed numerical solution $(\alpha)$ for the concentric, the strong eccentric and the mild eccentric inhomogeneity. Next to each (dashed line), we present the $\omega\equiv1$ profile along with the solution $(\alpha)$ provided by the Tikhonov reconstru\-ction. We display the numerical results in Table \ref{table:infoPhysical}.

\subsection[Multi - data reconstruction of the concentric problem]{Multi - data reconstruction of the concentric inhomogeneity}\label{subsec:multidatarecon}
\subsubsection{Exact observations.}\label{subsubsec:exactdatamd}
\begin{figure}[H] \onecolumn
	\hspace{-1.5cm}
	\begin{minipage}{.07\textwidth}
		\vspace{-.6cm}
		\includegraphics[width=0.65\textwidth,height=0.45\textheight]{colorbar12.png}
	\end{minipage}\hspace{.8cm}
	\begin{minipage}{.06\textwidth}
		\begin{itemize}\itemsep20em
			\vskip-14.5em
			\item[\footnotesize{$\boldsymbol{\ell=0.2}$}]
			\item[\footnotesize{$\boldsymbol{\ell=0.4}$}]
		\end{itemize}		
	\end{minipage}\hspace{.5cm}
	\centering
	\begin{minipage}{.85\textwidth}
	 \subfloat[$\alpha^0_{N=\lbrace 1,2,5\rbrace}$]{\includegraphics[width=0.19\linewidth]{sel02a0.png}}
		\subfloat[$\alpha_{(N=1)}$]{\includegraphics[width=0.19\linewidth]{pl02aend.png}}
		\subfloat[$\alpha_{(N=2)}$]{\includegraphics[width=0.19\linewidth]{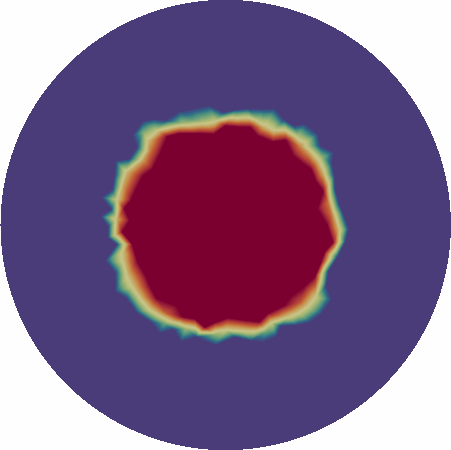}}
		\subfloat[$\alpha_{(N=5)}$]{\includegraphics[width=0.19\linewidth]{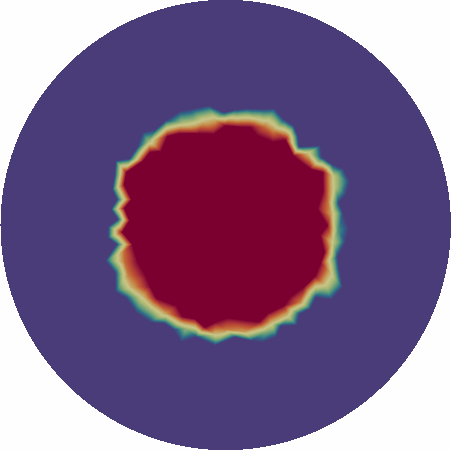}}\\
	   \subfloat[$\omega^0_{N=\lbrace 1,2,5\rbrace}$]{\includegraphics[width=0.19\linewidth]{pl02w0.png}}
		\subfloat[$\omega_{(N=1)}$]{\includegraphics[width=0.19\linewidth]{pl02wend.png}}
		\subfloat[$\omega_{(N=2)}$]{\includegraphics[width=0.19\linewidth]{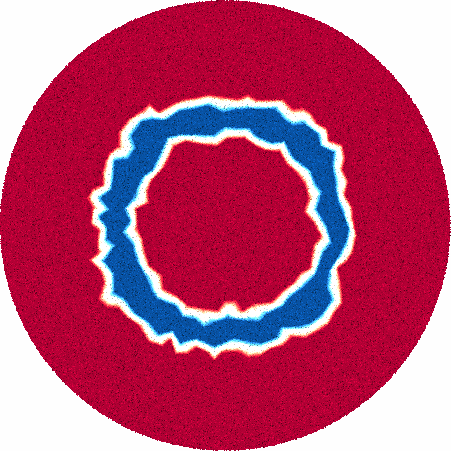}}
		\subfloat[$\omega_{(N=5)}$]{\includegraphics[width=0.19\linewidth]{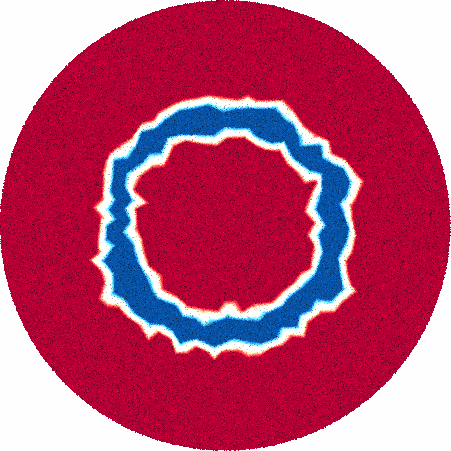}}	\\
		 \subfloat[$\alpha^0_{N=\lbrace 1,2,5\rbrace}$]{\includegraphics[width=0.19\linewidth]{sel02a0.png}}
		\subfloat[$\alpha_{(N=1)}$]{\includegraphics[width=0.19\linewidth]{pl04aend.png}}
		\subfloat[$\alpha_{(N=2)}$]{\includegraphics[width=0.19\linewidth]{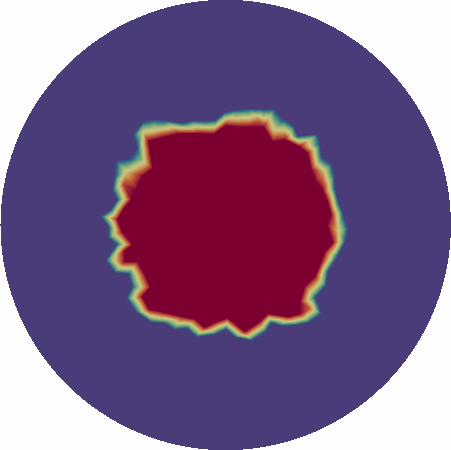}}
		\subfloat[$\alpha_{(N=5)}$]{\includegraphics[width=0.19\linewidth]{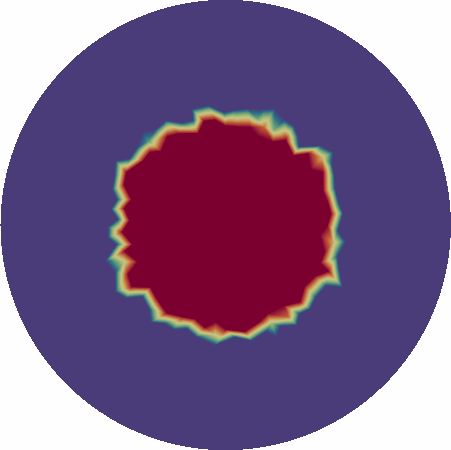}}\\
		\subfloat[$\omega^0_{N=\lbrace 1,2,5\rbrace}$]{\includegraphics[width=0.19\linewidth]{pl04w0.png}}
		\subfloat[$\omega_{(N=1)}$]{\includegraphics[width=0.19\linewidth]{pl04wend.png}}
		\subfloat[$\omega_{(N=2)}$]{\includegraphics[width=0.19\linewidth]{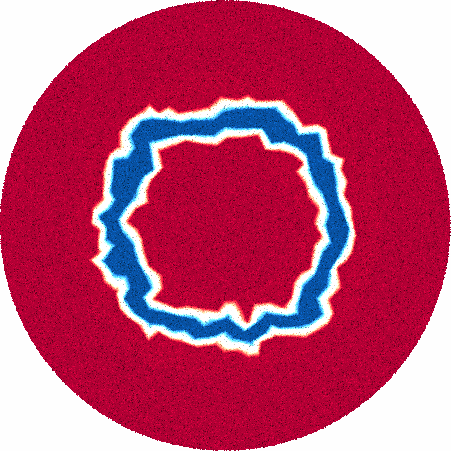}}
		\subfloat[$\omega_{(N=5)}$]{\includegraphics[width=0.19\linewidth]{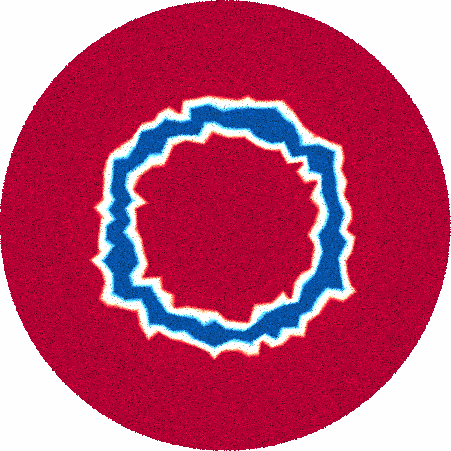}}	
	\end{minipage}
	\caption{$N-$ data $BV$ {\it Regularized Inversion} solution pairs $(\alpha, \omega)$ at the final iteration for $\ell=\lbrace 0.2, 0.4\rbrace$ of the concentric inhomogeneity problem and $\mu=\lbrace 1, 0.1 \rbrace$ respectively.\label{figure:multidataa}}
\end{figure}

In this section, we investigate a simple case of multiple measurements. Thoroughly, we consider a specific family of exact synthetic data for the concentric inhomogeneity problem (\ref{subsubsec:CONC}). More precisely, we assume that we know a sufficient number $m\;(m \in \mathbb{N})$ of boundary data pairs of the form
\begin{equation}
\begin{aligned}
 (f_m,g_m)=\Big(1+\frac{13}{8}\frac{(3\cdot 2^{2m+1}-2)}{m(3\cdot 2^{2m}+1)}\cos(m\phi),\frac{13}{8} \cos(m\phi)\Big),
\end{aligned}\label{multidatafg}
\end{equation}
where $\;\phi \in [0,2\pi),\; m=1,2,...,N$, each-one suppor\-ting the common exact solution
\begin{equation}
u^{(m)}(\rho,\phi)=1+\frac{13}{8} \begin{cases}
 \frac{2^{m+2}}{m(3\cdot 2^{2m}+1)}\rho^m \cos(m\phi),\\ \hspace{3cm} 0\leqslant\rho\leqslant1,\\
\frac{2^{m+1}}{m(3\cdot 2^{2m}+1)}(3\rho^m -\rho^{-m}) \cos(m\phi), \\ \hspace{3cm} 1 < \rho \leqslant 2.
\end{cases}
\end{equation}
We perform the algorithm of Section \ref{subsec:Algorithm} (setting $\lambda=0$) on the problem (\ref{problemtoreconstruct}), where now the functional (\ref{funtoreconstruct}) can be rewritten as
\begin{align}\label{funtoreconstructmanydata}
J_\epsilon(u,w,\alpha,\omega,N)&=\sum_{m=1}^{N}\vert\tilde{E}_m(u,w,\alpha)\vert\nonumber \\ &+\frac{\lambda}{2}\|R\alpha-\alpha_0\|^2_{L^2(\Omega)}\nonumber \\ &+\mu \int_\Omega (\omega |\nabla \alpha{|}^2+\psi_\epsilon(\omega))dx
\end{align}
where
\begin{align}
\tilde{E}_m(u,w,\alpha)&=\frac{\kappa+1}{2}\int_{\Omega}\alpha|\nabla (u+\eta f_m)|^2 dx \nonumber \\ & + (1-\kappa)\left[\frac{1}{2}\int_{\Omega}\alpha |\nabla w|^2 d x  -\left\langle g_m,\gamma w \right\rangle_{\mathcal{A} \times \mathcal{B}}\right]\nonumber\\ & -\kappa\left\langle g_m,f_m \right\rangle_{\mathcal{A} \times \mathcal{B}}. 
\end{align}
The problem parameters $\epsilon,\ \kappa$ and the mesh size $h$ are selected as in (\ref{subsubsec:CONC}), whereas we choose $\alpha_0=2.5$ through the whole domain $\Omega\setminus\Omega_{\delta}$. We repeat the inversion for the $\omega^0 -$ configurations: $\ell=\lbrace 0.2,0.4,0.6\rbrace$ and $\omega^0\equiv1$ with regularization values $\mu=\lbrace1,0.1,1,1 \rbrace$ apiece, for $N=\lbrace 2,5\rbrace$ data pairs of the form (\ref{multidatafg}) and we compare the results with the single data pair case $N=1$ (\ref{subsubsec:CONC}). We notice here that the boundary data (\ref{multidatafg}) for $N=1$ are normalized to coincide with those in (\ref{pilotfg}).\newline
Figure \ref{figure:multidataa}a and \ref{figure:multidataa}b from left to the right displays the initial conductivity $(\alpha^0)$ and the final reconstructed solutions $(\alpha)$ for $N=\lbrace 1,2,5 \rbrace$ data pairs along with their corresponding solution pair $(\omega)$ (beneath each), for the aforementioned $\omega^0 -$ configurations. Table \ref{table:infoMD} presents the numerical results for every selection of the number of the boundary data pairs $(N)$ and every choice of the iteration level $n$ $(n=10)$. \newline
We observe that the use of multiple measurements accelerates the computational process to an ac\-ce\-ptable convergence both in the geometrical and physical sense. It is remarkable that for this mesh case the reconstruction with just $N=2$ data pairs is exact, coinciding with the inversion based on several data pairs $(N>2)$.

\subsubsection{Noisy observations.}
In this example, we assume that noisy observations are available. More precisely, the range of the $\Lambda_{N.t.D}$ o\-pe\-rator is subject to noise additions in the form 
\begin{equation}
\begin{aligned}
(\tilde{f}_m,\tilde{g}_m)=\big(f_m+\vert f_m\vert\cdot \mathcal{R}_{f_m}\cdot \theta, \ g_m\big),
\end{aligned}\label{multidatafgnoise}
\end{equation}
\  where $(f_m,g_m)$ are given in (\ref{multidatafg}), $\mathcal{R}_{f_m}$ are $\partial\Omega^h \times 1 -$ matrices of random numbers on the interval $(-1,1)$, which are generated by the MATLAB function ``rand'' with $\partial\Omega^h$ the number of boundary nodes of the Delaunay triangulation and $\theta$ some noise level.
\newline
The minimization process on (\ref{funtoreconstructmanydata}) will be started with $\lambda=0,\ \epsilon=0.1,\ \kappa=10$ and $\alpha_0=2.5$ as the initial approximation for the con\-ductivity function throu\-ghout the domain $\Omega\setminus\Omega_{\delta}$. We perform $n=10$ ite\-ra\-tions of the {\it $BV$ Regu\-larized Inversion} algorithm for the $\omega^0-$ configurations: $\ell=\lbrace 0.2,0.4,0.6\rbrace$ with respect to the mesh size $h=0.27$. At each $\omega^0-$ refinement level $(\ell)$, we recur the implementation for each value of the noise level $\theta=\lbrace 0.005,0.01,0.05\rbrace$ and e\-ve\-ry choice of the number of the boundary data pairs $N=\lbrace 1,2,5\rbrace$.
\newline
 Figure \ref{figure:noisel02md}a-\ref{figure:noisel06md}c performs from left to the right the reconstructed solution pair $(\alpha,\omega)$ at the final iteration for each number of boundary data pairs $(N)$ as the measurement noise $\theta$ becomes relatively larger. At each $\ell-$ test case, the initial guess $(\alpha^0,\omega^0)$ is depicted on the right-hand side of the color-bar. We start from the finest level $\ell=0.2$. The nume\-rical results are summarized in Table \ref{table:infoNMDl02}a-\ref{table:infoNMDl06}c, where we present the uniform values $(\alpha_{in,out})$ of the computed solution $\alpha$ at every iterate $n$ from the finest $(\ell=0.2)$ to the coarsest refinement level $(\ell=0.6)$.
 \newline In the realm of the noisy observations, the Tikhonov approach $(\omega^0\equiv1)$ totally fails, as expected, and therefore is omitted. 
\newline
To conclude, it is revealing how the reconstruction with multiple measurements forcefully contributes to the geometrical and physical performance of the computed solution even in the presence of both relatively large noise and higher uncertainty with respect to the location of the inclusion's boundary. 
 
\clearpage

\section*{Figures and figure captions}

\captionsetup[figure]{labelformat=cancaptionlabel}

\addtocounter{figure}{-1}
\captionsetup[figure]{labelformat=cancaptionlabel2}
\begin{figure}[H] \onecolumn
	\hspace{-1.5cm}
	\begin{minipage}{.07\textwidth}
		\vspace{-.6cm}
		\includegraphics[width=0.65\textwidth,height=0.45\textheight]{colorbar12.png}
	\end{minipage}\hspace{.8cm}
	\begin{minipage}{.06\textwidth}
		\begin{itemize}\itemsep20em
			\vskip-14.5em
			\item[\footnotesize {$\boldsymbol{\ell=0.6}$}]
			\item[\footnotesize{$\boldsymbol{\omega\equiv1}$}]
		\end{itemize}		
	\end{minipage}\hspace{.5cm}
	\centering
	\begin{minipage}{.85\textwidth}
		\subfloat[$\alpha^0_{N=\lbrace 1,2,5\rbrace}$]{\includegraphics[width=0.19\linewidth]{sel02a0.png}}
		\subfloat[$\alpha_{(N=1)}$]{\includegraphics[width=0.19\linewidth]{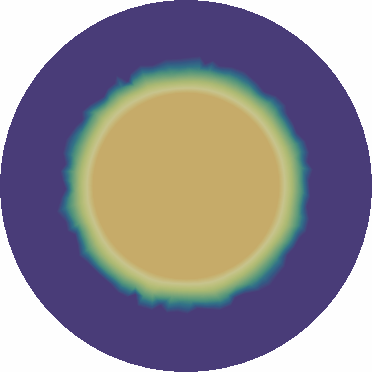}}
		\subfloat[$\alpha_{(N=2)}$]{\includegraphics[width=0.19\linewidth]{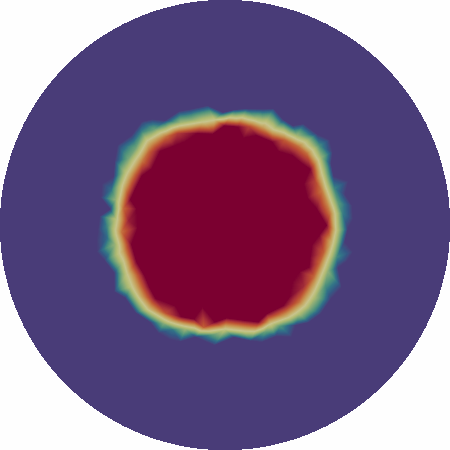}}
		\subfloat[$\alpha_{(N=5)}$]{\includegraphics[width=0.19\linewidth]{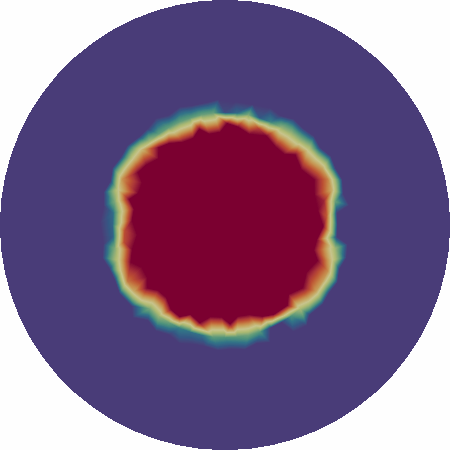}}\\
		\subfloat[$\omega^0_{N=\lbrace 1,2,5\rbrace}$]{\includegraphics[width=0.19\linewidth]{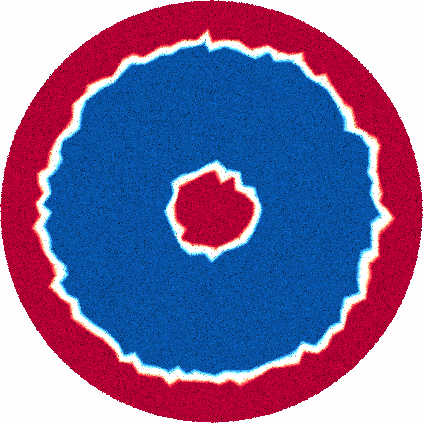}}
		\subfloat[$\omega_{(N=1)}$]{\includegraphics[width=0.19\linewidth]{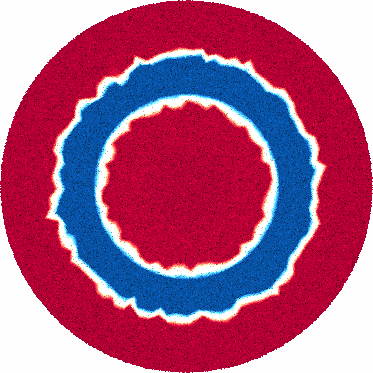}}
		\subfloat[$\omega_{(N=2)}$]{\includegraphics[width=0.19\linewidth]{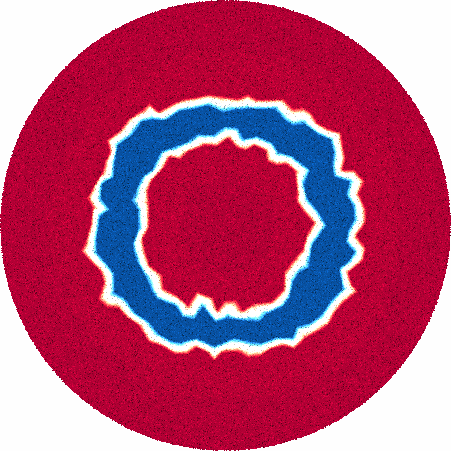}}
		\subfloat[$\omega_{(N=5)}$]{\includegraphics[width=0.19\linewidth]{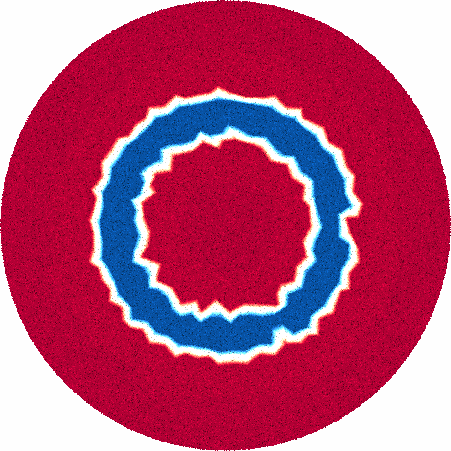}}	\\
		\subfloat[$\alpha^0_{N=\lbrace 1,2,5\rbrace}$]{\includegraphics[width=0.19\linewidth]{sel02a0.png}}
		\subfloat[$\alpha_{(N=1)}$]{\includegraphics[width=0.19\linewidth]{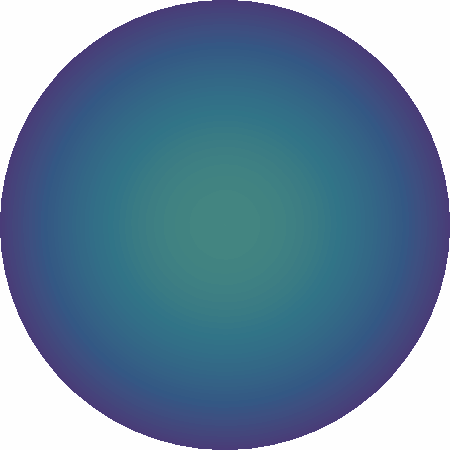}}
		\subfloat[$\alpha_{(N=2)}$]{\includegraphics[width=0.19\linewidth]{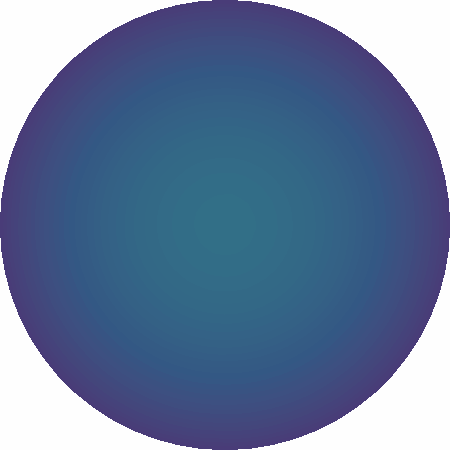}}
		\subfloat[$\alpha_{(N=5)}$]{\includegraphics[width=0.19\linewidth]{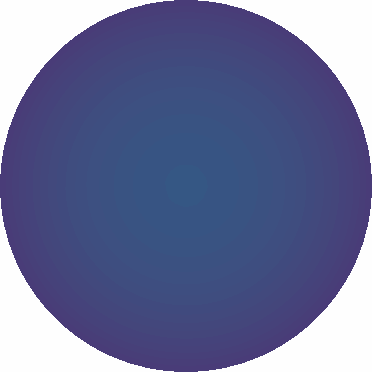}}\\
		\subfloat[$\omega^0_{N=\lbrace 1,2,5\rbrace}$]{\includegraphics[width=0.19\linewidth]{Tpw.png}}
		\subfloat[$\omega_{(N=1)}$]{\includegraphics[width=0.19\linewidth]{Tpw.png}}
		\subfloat[$\omega_{(N=2)}$]{\includegraphics[width=0.19\linewidth]{Tpw.png}}
		\subfloat[$\omega_{(N=5)}$]{\includegraphics[width=0.19\linewidth]{Tpw.png}}
	\end{minipage}
	\caption{$N-$ data $BV$ {\it Regularized Inversion} solution pairs $(\alpha, \omega)$ at the final iteration for $(\ell=0.6, \omega\equiv1)$ of the concentric inhomogeneity problem and $\mu=\lbrace 1, 1 \rbrace$ respectively.\label{figure:multidatab}}
\end{figure}

\captionsetup[figure]{labelformat=cancaptionlabel}
  \begin{figure}[H]\onecolumn
  	\hspace{-1.5cm}
  	\begin{minipage}{.08\textwidth}
  		\includegraphics[width=0.66\textwidth,height=0.47\textheight]{colorbar12.png}
  	\end{minipage}\hspace{.4cm}
  	\begin{minipage}{.6\textwidth}
  		\subfloat[$\alpha^0_{N=\lbrace 1,2,5\rbrace}$]{\includegraphics[width=0.24\linewidth]{sel02a0.png}}\\
  		\subfloat[$\omega^0_{N=\lbrace 1,2,5\rbrace}$]{\includegraphics[width=0.24\linewidth]{pl02w0.png}}
  	\end{minipage}\hspace{-6.3cm}
  	\begin{minipage}{.05\textwidth}
  		\begin{itemize}\itemsep20em
  			\vskip-14.5em
  			\item[\footnotesize $\boldsymbol{\theta=0.005}$]
  			\item[\footnotesize $\boldsymbol{\theta=0.01}$]
  			\item[\footnotesize $\boldsymbol{\theta=0.05}$]
  		\end{itemize}
  	\end{minipage}\hspace{.4cm}
  	\begin{minipage}{.65\textwidth}		
  		\subfloat[$\alpha_{(N=1)}$]{\includegraphics[width=0.27\linewidth]{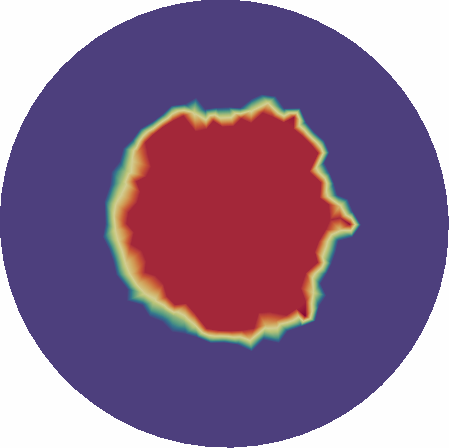}}
  		\subfloat[$\alpha_{(N=2)}$]{\includegraphics[width=0.27\linewidth]{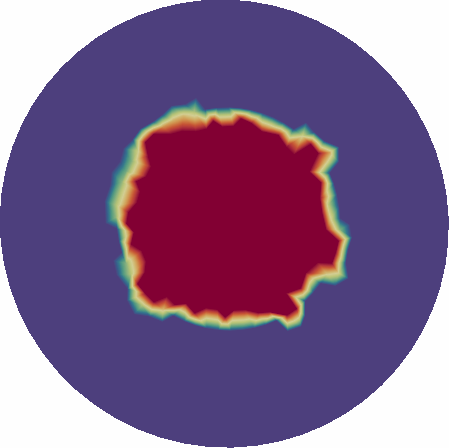}}
  		\subfloat[$\alpha_{(N=5)}$]{\includegraphics[width=0.27\linewidth]{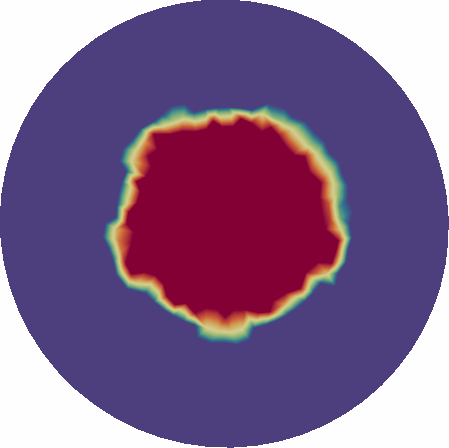}}  \\
  		\subfloat[$\omega_{(N=1)}$]{\includegraphics[width=0.27\linewidth]{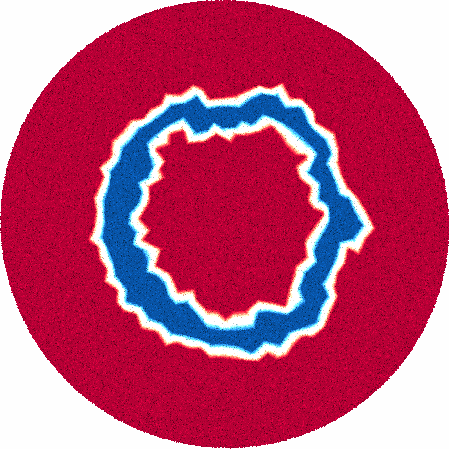}}
  		\subfloat[$\omega_{(N=2)}$]{\includegraphics[width=0.27\linewidth]{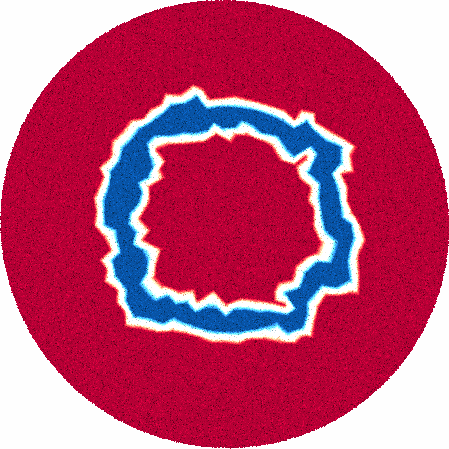}}
  		\subfloat[$\omega_{(N=5)}$]{\includegraphics[width=0.27\linewidth]{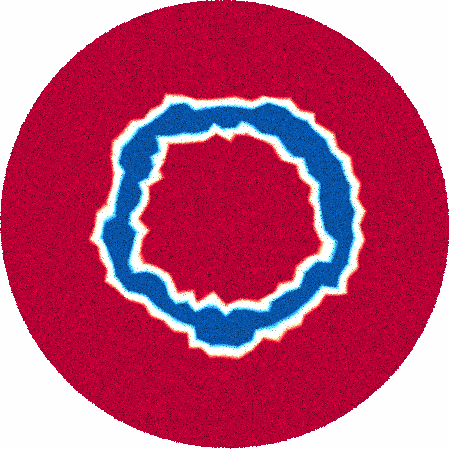}}	\\
  		\subfloat[$\alpha_{(N=1)}$]{\includegraphics[width=0.27\linewidth]{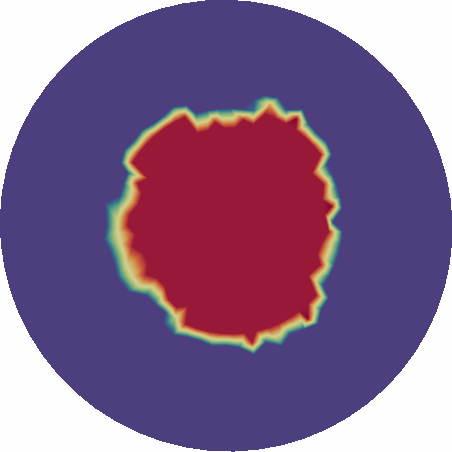}}
  		\subfloat[$\alpha_{(N=2)}$]{\includegraphics[width=0.27\linewidth]{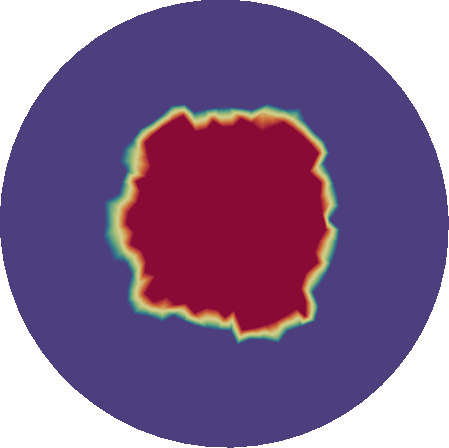}}
  		\subfloat[$\alpha_{(N=5)}$]{\includegraphics[width=0.27\linewidth]{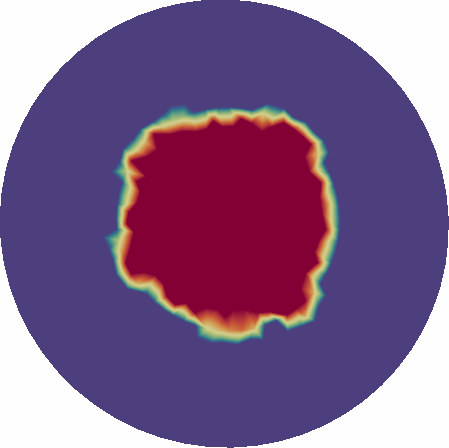}}\\
  		\subfloat[$\omega_{(N=1)}$]{\includegraphics[width=0.27\linewidth]{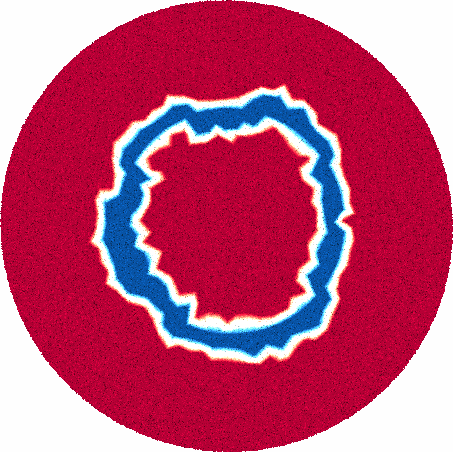}}
  		\subfloat[$\omega_{(N=2)}$]{\includegraphics[width=0.27\linewidth]{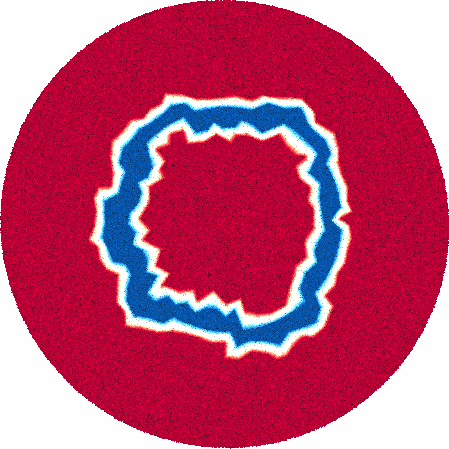}}
  		\subfloat[$\omega_{(N=5)}$]{\includegraphics[width=0.27\linewidth]{aNoise001N5l02}}\\
  		\subfloat[$\alpha_{(N=1)}$]{\includegraphics[width=0.27\linewidth]{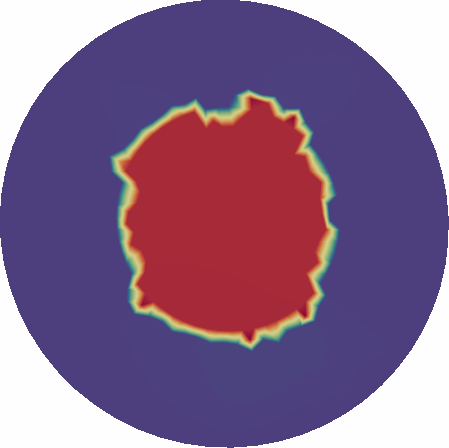}}
  		\subfloat[$\alpha_{(N=2)}$]{\includegraphics[width=0.27\linewidth]{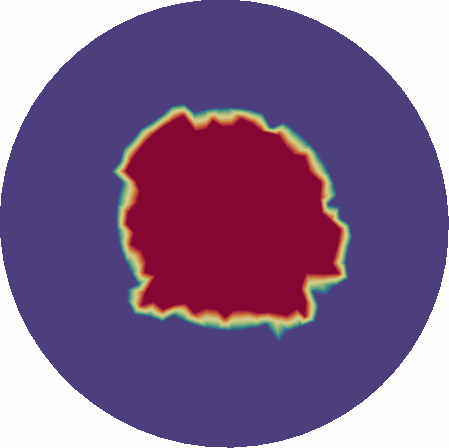}}
  		\subfloat[$\alpha_{(N=5)}$]{\includegraphics[width=0.27\linewidth]{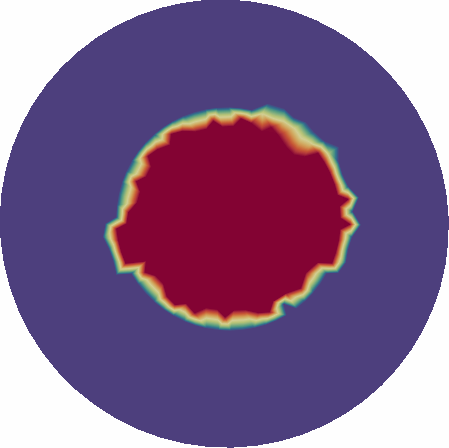}}\\
  		\subfloat[$\omega_{(N=1)}$]{\includegraphics[width=0.27\linewidth]{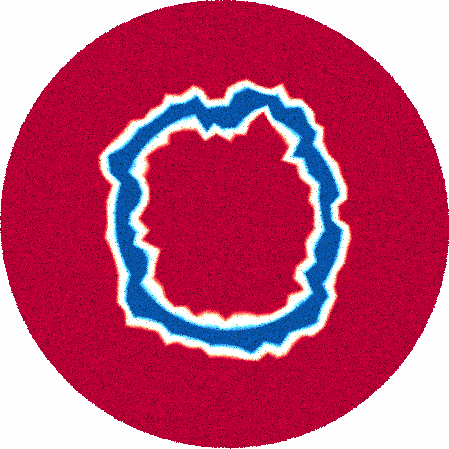}}
  		\subfloat[$\omega_{(N=2)}$]{\includegraphics[width=0.27\linewidth]{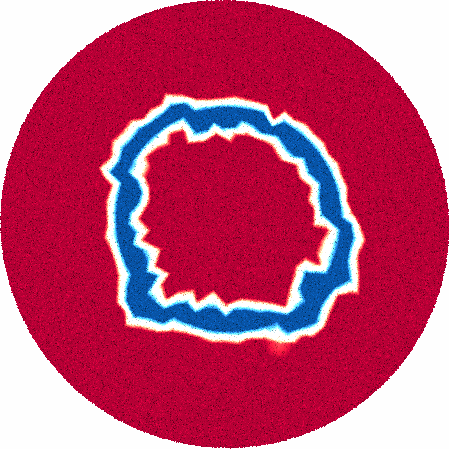}}
  		\subfloat[$\omega_{(N=5)}$]{\includegraphics[width=0.27\linewidth]{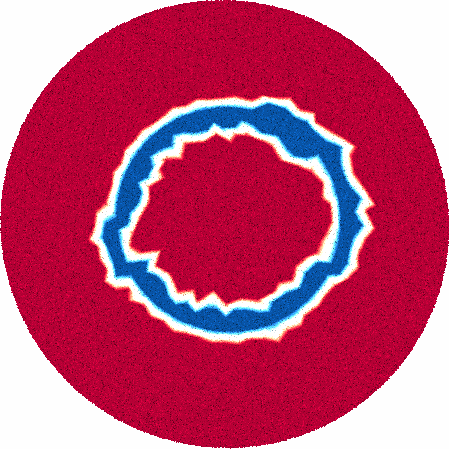}}
  	\end{minipage}
  	\caption{$N-$ data $BV$ {\it Regularized Inversion} solution pairs $(\alpha, \omega)$ at the final iteration for $\ell= 0.2$ of the concentric inhomogeneity problem with measurement noise $\theta=\lbrace 0.5\%, 1\%, 5\% \rbrace$ and $\mu=1$. \label{figure:noisel02md}}
  \end{figure}

\addtocounter{figure}{-1}
\captionsetup[figure]{labelformat=cancaptionlabel2}  
\begin{figure}[H]\onecolumn
	\hspace{-1.5cm}
	\begin{minipage}{.08\textwidth}
		\includegraphics[width=0.66\textwidth,height=0.47\textheight]{colorbar12.png}
	\end{minipage}\hspace{.4cm}
	\begin{minipage}{.6\textwidth}
		\subfloat[$\alpha^0_{N=\lbrace 1,2,5\rbrace}$]{\includegraphics[width=0.24\linewidth]{sel02a0.png}}\\
		\subfloat[$\omega^0_{N=\lbrace 1,2,5\rbrace}$]{\includegraphics[width=0.24\linewidth]{pl04w0.png}}
	\end{minipage}\hspace{-6.3cm}
	\begin{minipage}{.05\textwidth}
		\begin{itemize}\itemsep20em
			\vskip-14.5em
			\item[\footnotesize $\boldsymbol{\theta=0.005}$]
			\item[\footnotesize $\boldsymbol{\theta=0.01}$]
			\item[\footnotesize $\boldsymbol{\theta=0.05}$]
		\end{itemize}
	\end{minipage}\hspace{.4cm}
	\begin{minipage}{.65\textwidth}		
		\subfloat[$\alpha_{(N=1)}$]{\includegraphics[width=0.27\linewidth]{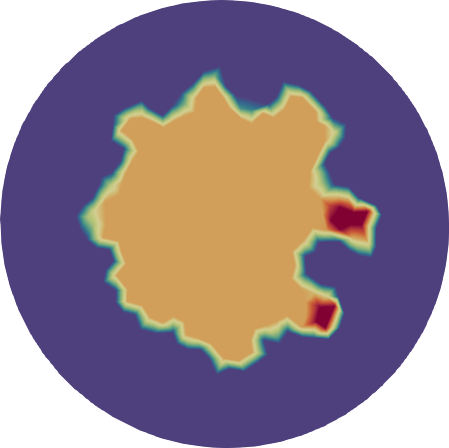}}
		\subfloat[$\alpha_{(N=2)}$]{\includegraphics[width=0.27\linewidth]{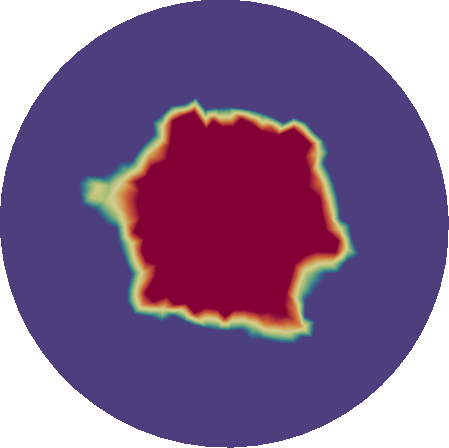}}
		\subfloat[$\alpha_{(N=5)}$]{\includegraphics[width=0.27\linewidth]{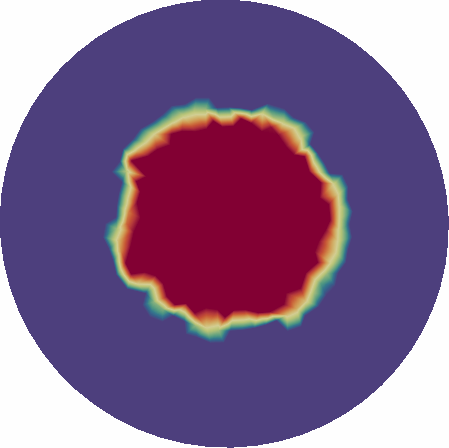}}  \\
		\subfloat[$\omega_{(N=1)}$]{\includegraphics[width=0.27\linewidth]{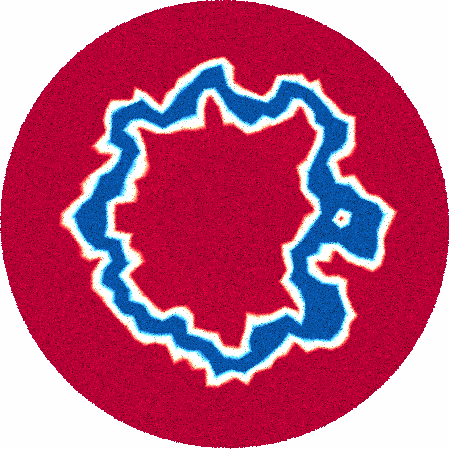}}
		\subfloat[$\omega_{(N=2)}$]{\includegraphics[width=0.27\linewidth]{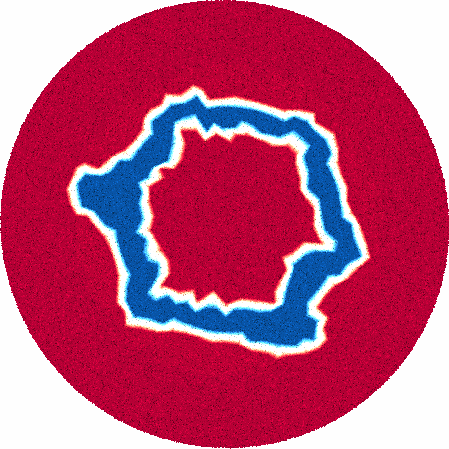}}
		\subfloat[$\omega_{(N=5)}$]{\includegraphics[width=0.27\linewidth]{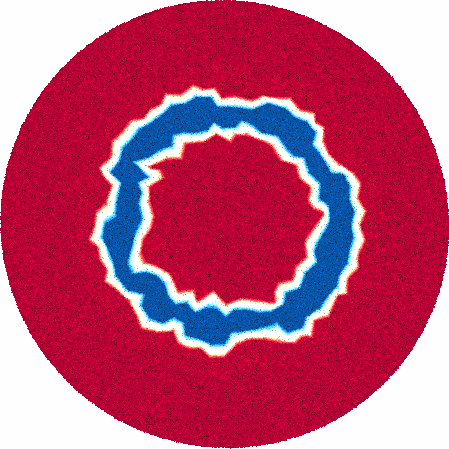}}	\\
		\subfloat[$\alpha_{(N=1)}$]{\includegraphics[width=0.27\linewidth]{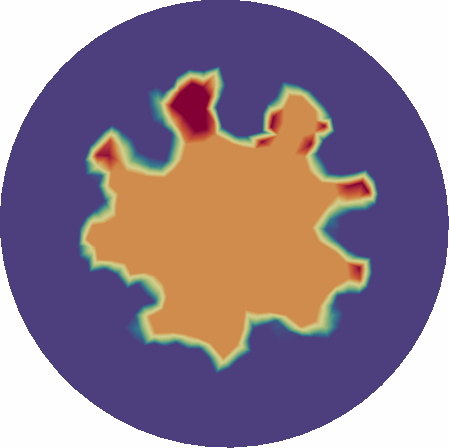}}
		\subfloat[$\alpha_{(N=2)}$]{\includegraphics[width=0.27\linewidth]{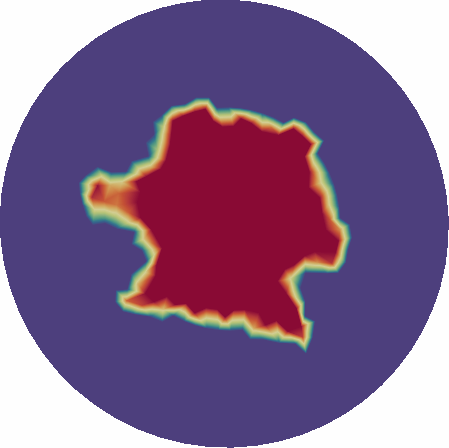}}
		\subfloat[$\alpha_{(N=5)}$]{\includegraphics[width=0.27\linewidth]{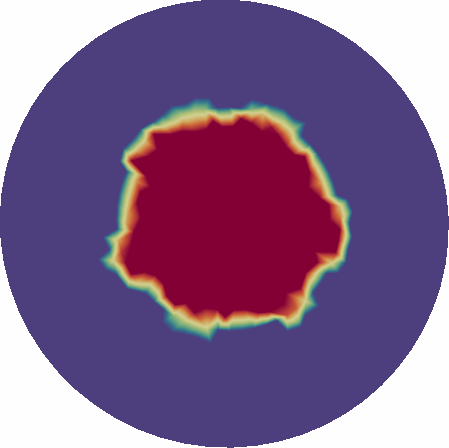}}\\
		\subfloat[$\omega_{(N=1)}$]{\includegraphics[width=0.27\linewidth]{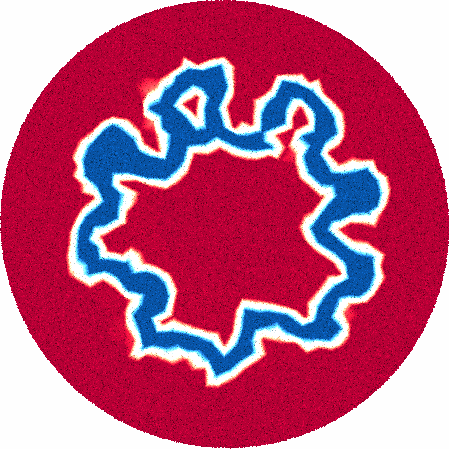}}
		\subfloat[$\omega_{(N=2)}$]{\includegraphics[width=0.27\linewidth]{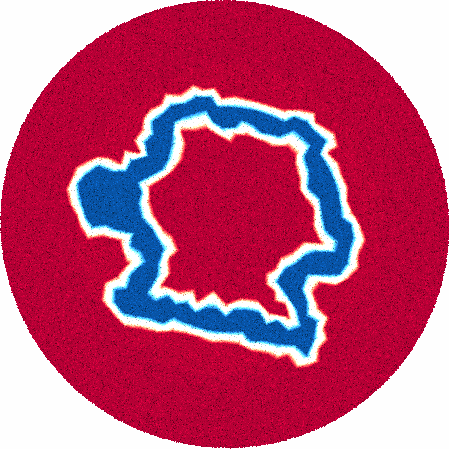}}
		\subfloat[$\omega_{(N=5)}$]{\includegraphics[width=0.27\linewidth]{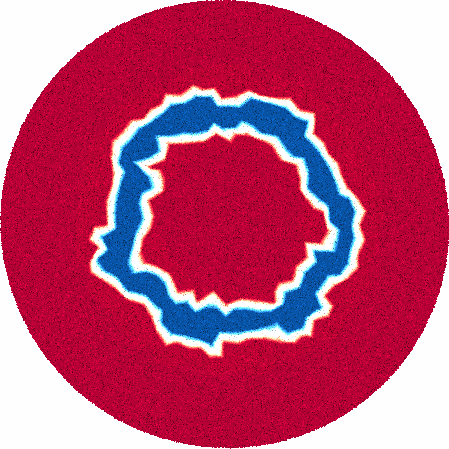}}\\
		\subfloat[$\alpha_{(N=1)}$]{\includegraphics[width=0.27\linewidth]{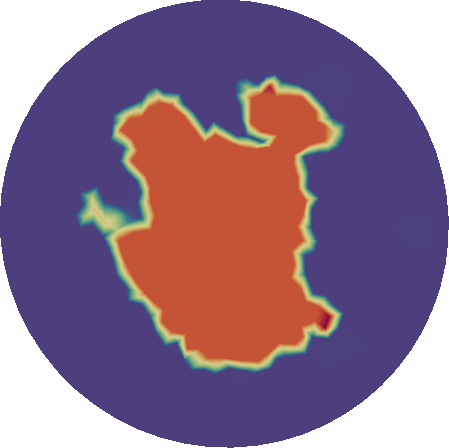}}
		\subfloat[$\alpha_{(N=2)}$]{\includegraphics[width=0.27\linewidth]{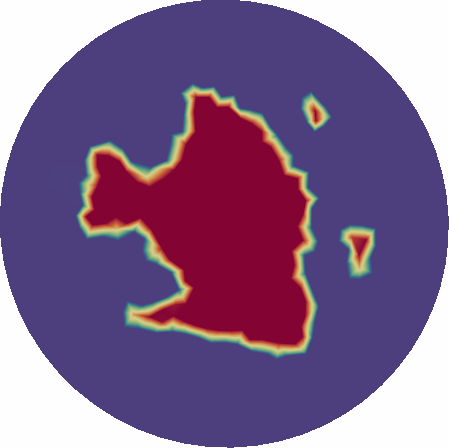}}
		\subfloat[$\alpha_{(N=5)}$]{\includegraphics[width=0.27\linewidth]{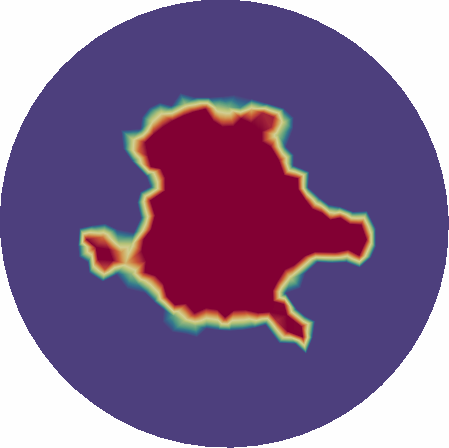}}\\
		\subfloat[$\omega_{(N=1)}$]{\includegraphics[width=0.27\linewidth]{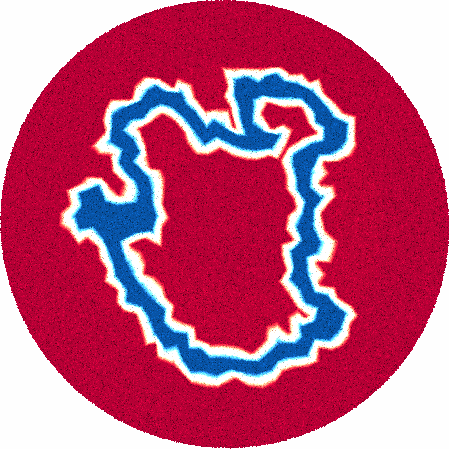}}
		\subfloat[$\omega_{(N=2)}$]{\includegraphics[width=0.27\linewidth]{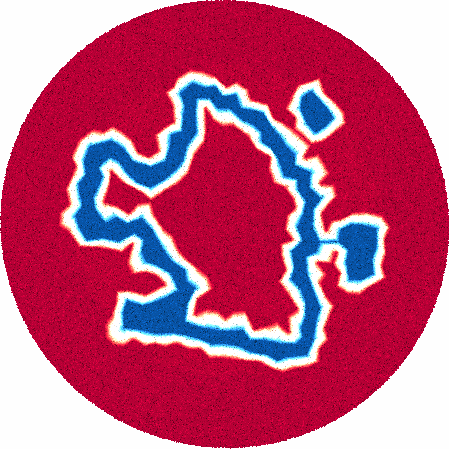}}
		\subfloat[$\omega_{(N=5)}$]{\includegraphics[width=0.27\linewidth]{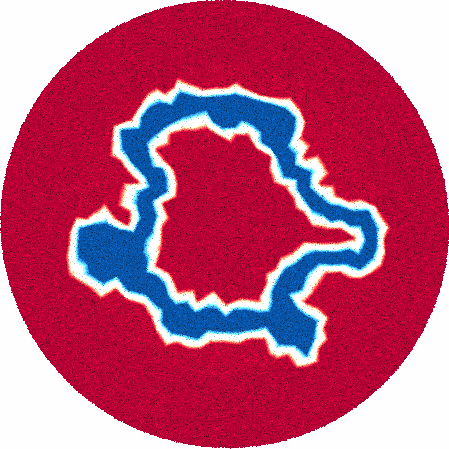}}
	\end{minipage}
	\caption{$N-$ data $BV$ {\it Regularized Inversion} solution pairs $(\alpha, \omega)$ at the final iteration for $\ell= 0.4$ of the concentric inhomogeneity problem with measurement noise $\theta=\lbrace 0.5\%, 1\%, 5\% \rbrace$ and $\mu=1$. \label{figure:noisel04md}}
\end{figure}

\addtocounter{figure}{-1}
\captionsetup[figure]{labelformat=cancaptionlabel3}  
\begin{figure}[H]\onecolumn
	\hspace{-1.5cm}
	\begin{minipage}{.08\textwidth}
		\includegraphics[width=0.66\textwidth,height=0.47\textheight]{colorbar12.png}
	\end{minipage}\hspace{.4cm}
	\begin{minipage}{.6\textwidth}
		\subfloat[$\alpha^0_{N=\lbrace 1,2,5\rbrace}$]{\includegraphics[width=0.24\linewidth]{sel02a0.png}}\\
		\subfloat[$\omega^0_{N=\lbrace 1,2,5\rbrace}$]{\includegraphics[width=0.24\linewidth]{w0l06.png}}
	\end{minipage}\hspace{-6.3cm}
	\begin{minipage}{.05\textwidth}
		\begin{itemize}\itemsep20em
			\vskip-14.5em
			\item[\footnotesize $\boldsymbol{\theta=0.005}$]
			\item[\footnotesize $\boldsymbol{\theta=0.01}$]
			\item[\footnotesize $\boldsymbol{\theta=0.05}$]
		\end{itemize}
	\end{minipage}\hspace{.4cm}
	\begin{minipage}{.65\textwidth}		
		\subfloat[$\alpha_{(N=1)}$]{\includegraphics[width=0.27\linewidth]{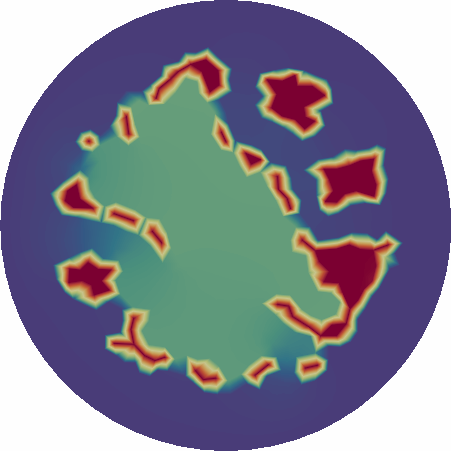}}
		\subfloat[$\alpha_{(N=2)}$]{\includegraphics[width=0.27\linewidth]{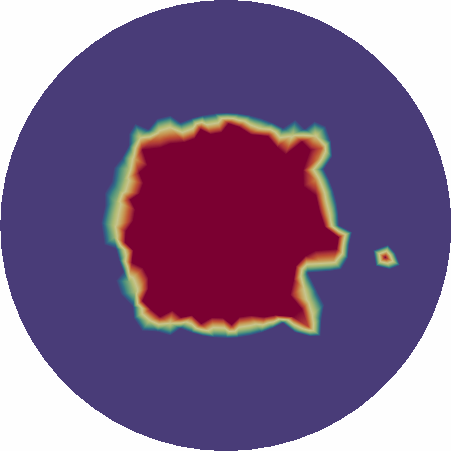}}
		\subfloat[$\alpha_{(N=5)}$]{\includegraphics[width=0.27\linewidth]{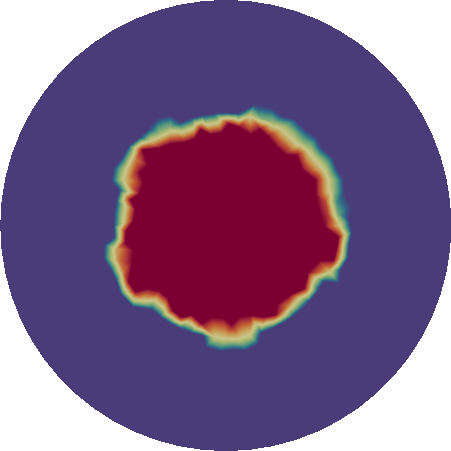}}  \\
		\subfloat[$\omega_{(N=1)}$]{\includegraphics[width=0.27\linewidth]{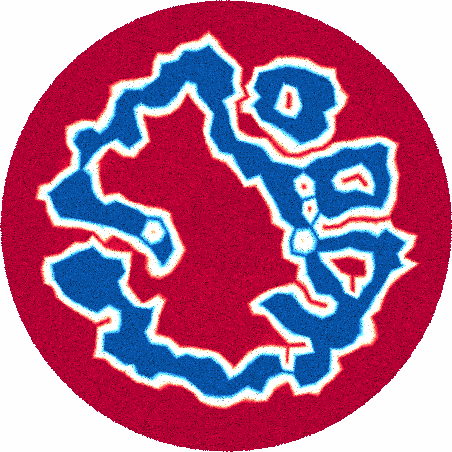}}
		\subfloat[$\omega_{(N=2)}$]{\includegraphics[width=0.27\linewidth]{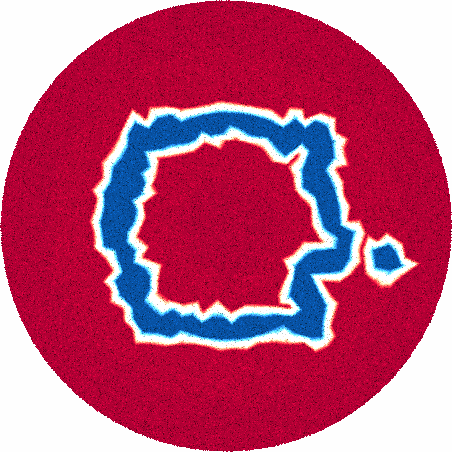}}
		\subfloat[$\omega_{(N=5)}$]{\includegraphics[width=0.27\linewidth]{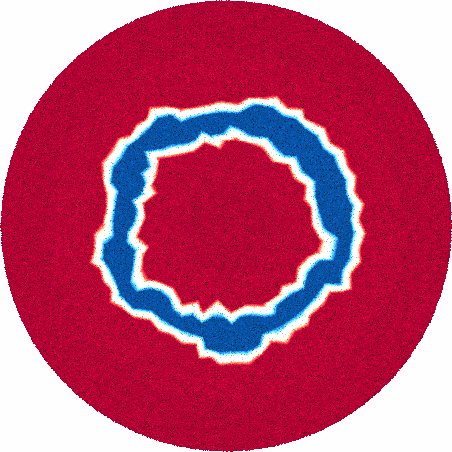}}	\\
		\subfloat[$\alpha_{(N=1)}$]{\includegraphics[width=0.27\linewidth]{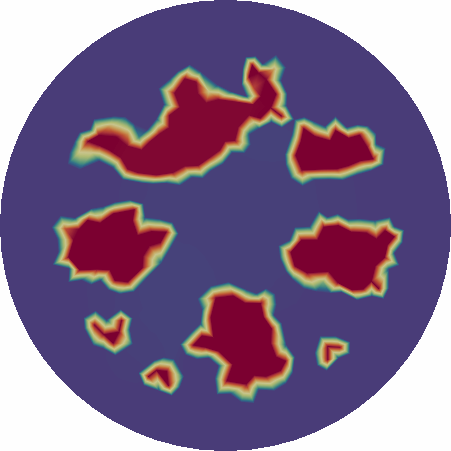}}
		\subfloat[$\alpha_{(N=2)}$]{\includegraphics[width=0.27\linewidth]{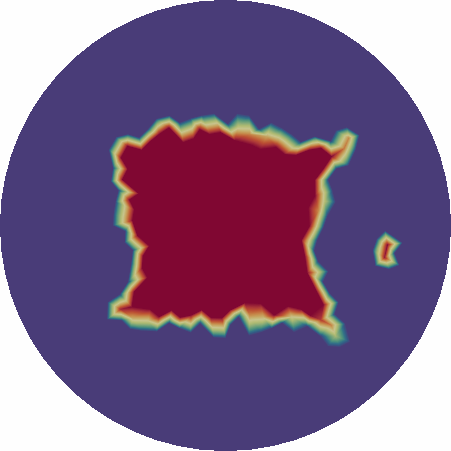}}
		\subfloat[$\alpha_{(N=5)}$]{\includegraphics[width=0.27\linewidth]{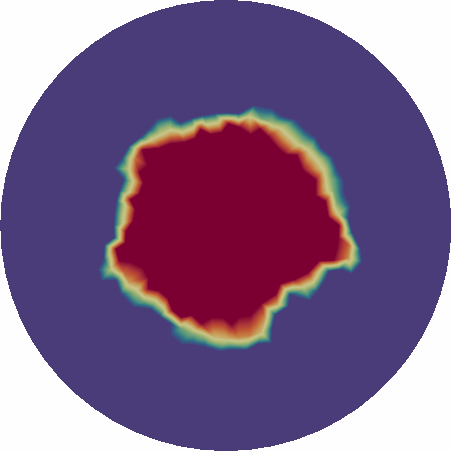}}\\
		\subfloat[$\omega_{(N=1)}$]{\includegraphics[width=0.27\linewidth]{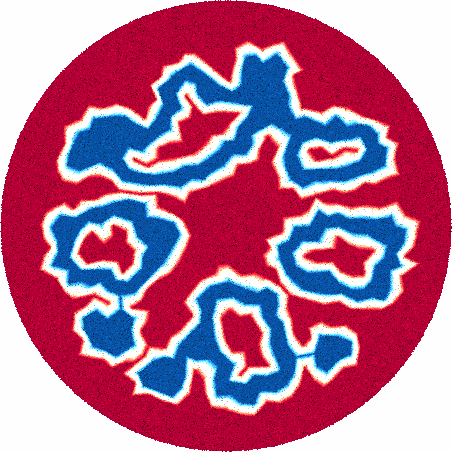}}
		\subfloat[$\omega_{(N=2)}$]{\includegraphics[width=0.27\linewidth]{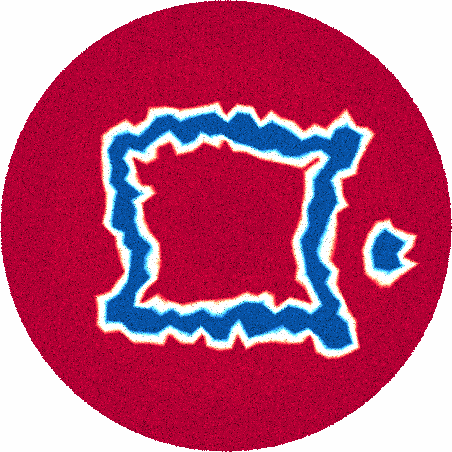}}
		\subfloat[$\omega_{(N=5)}$]{\includegraphics[width=0.27\linewidth]{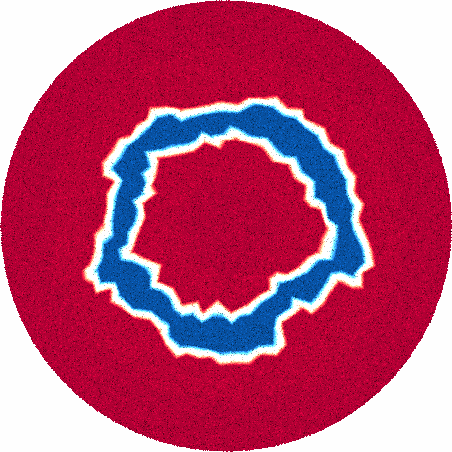}}\\
		\subfloat[$\alpha_{(N=1)}$]{\includegraphics[width=0.27\linewidth]{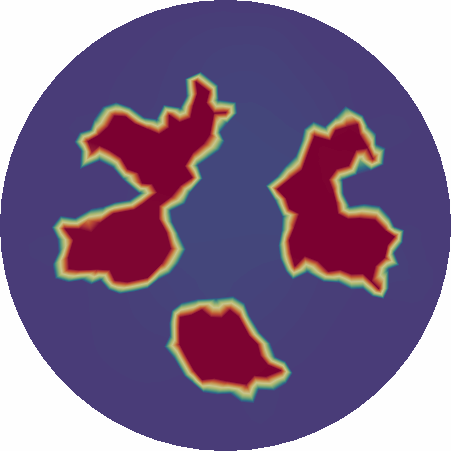}}
		\subfloat[$\alpha_{(N=2)}$]{\includegraphics[width=0.27\linewidth]{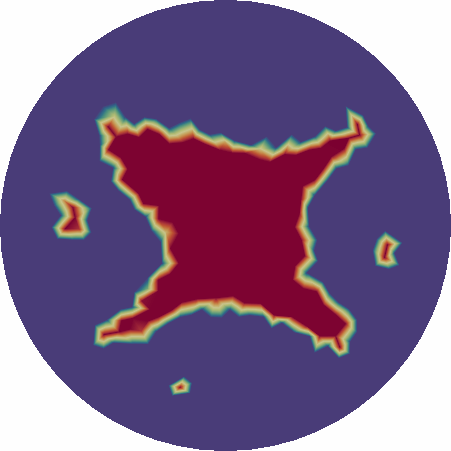}}
		\subfloat[$\alpha_{(N=5)}$]{\includegraphics[width=0.27\linewidth]{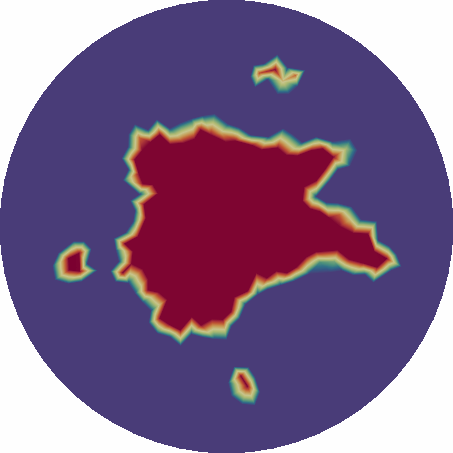}}\\
		\subfloat[$\omega_{(N=1)}$]{\includegraphics[width=0.27\linewidth]{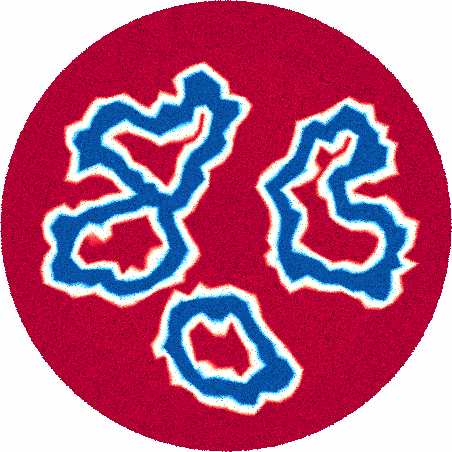}}
		\subfloat[$\omega_{(N=2)}$]{\includegraphics[width=0.27\linewidth]{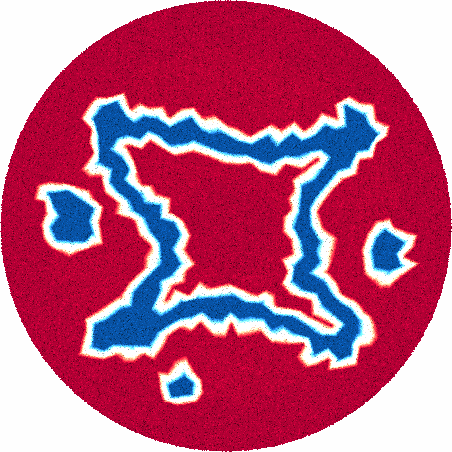}}
		\subfloat[$\omega_{(N=5)}$]{\includegraphics[width=0.27\linewidth]{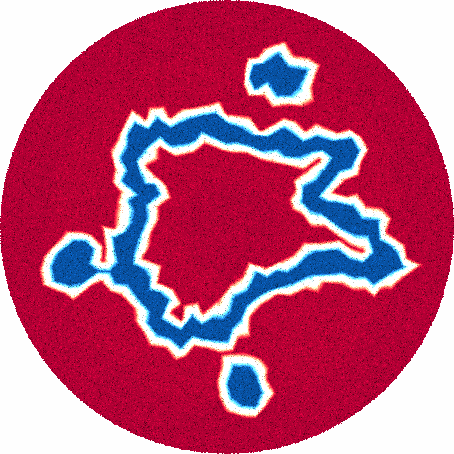}}
	\end{minipage}
	\caption{$N-$ data $BV$ {\it Regularized Inversion} solution pairs $(\alpha, \omega)$ at the final iteration for $\ell= 0.6$ of the concentric inhomogeneity problem with measurement noise $\theta=\lbrace 0.5\%, 1\%, 5\% \rbrace$ and $\mu=1$. \label{figure:noisel06md}}
\end{figure}

 \clearpage
\section*{Tables of results}
\setcounter{table}{2}
\begin{table}[ht]\onecolumn
 \centering
\caption{\label{table:infoP3Conc}Reconstructed values $(\alpha_{in},\alpha_{out})$ of the target conductivity $\alpha$ for the concentric inhomogeneity problem.}
 \begin{tabular}{|p{2cm}|l | c c c c c c c|}
  \bottomrule
\multicolumn{9}{|c|}{Number of iterations $n=10$} \\
\bottomrule
\multirow{2}{4cm}{$\ell=0.2$}&  \multicolumn{1}{l}{$n$} & $1$ & $2$ &$3$ & $\geq 4$ & & &\\
 \cline{2-6}
 \multirow{2}{4cm}{$\mu=1$}&$\alpha_{in}$ & $1.92$ & $1.94$ &$1.93$ &$1.94$ & & & \\
 &$\alpha_{out}$ & $1.00$ & $1.00$ &$1.00$ & $1.00$ & & &   \\
 \midrule
 \multirow{2}{4cm}{$\ell=0.3$}&  \multicolumn{1}{l}{$n$}& $1$ & $2-3$ &$4-5$ & $6$ & $\geq7$ & &\\
 \cline{2-7}
 \multirow{2}{4cm}{$\mu=0.1$}&$\alpha_{in}$ & $1.84$ & $1.86$ &$1.87$ & $1.86$ & $1.83$ & & \\
 &$\alpha_{out}$ & $1.00$ & $1.00$ &$1.00$ & $1.00$ & $1.00$ & & \\
 \midrule
 \multirow{2}{4cm}{$\ell=0.4$}&  \multicolumn{1}{l}{$n$} & $1$ & $2$ &$3-4$ & $5$& $6-7$& $8-9$ &$10$\\
 \cline{2-9}
\multirow{2}{4cm}{$\mu=0.1$}&$\alpha_{in}$ & $1.75$ & $1.76$ &$1.74$ & $1.76$ & $1.78$  & $1.79$ &$1.80$   \\
 &$\alpha_{out}$ & $1.00$ & $1.00$ & $1.00$ & $1.00$  & $1.00$ &$1.00$ &$1.00$  \\
 \midrule
 \multirow{2}{4cm}{$\omega\equiv1$}&  \multicolumn{1}{l}{$n$} & $1$ & $\geq2$ & &  & & & \\
 \cline{2-4}
\multirow{2}{4cm}{$\mu=1$}&$\alpha_{in}$ & $1.12$ & $1.17$ & &  & &  &  \\
&$\alpha_{out}$ & $1.00$ & $1.00$ &  &  &   & &  \\
 \bottomrule
 \end{tabular}
 \end{table}

\begin{table}[H]
\centering
\caption{\label{table:infoP3SE}Reconstructed values $(\alpha_{in},\alpha_{out})$ of the target conductivity $\alpha$ for the strong eccentricity problem.}
  \begin{tabular}{|p{2cm}|l | c c c c c c|}
  \bottomrule
\multicolumn{8}{|c|}{Number of iterations $n=10$} \\
\bottomrule
\multirow{2}{4cm}{$\ell=0.2$}& \multicolumn{1}{l}{$n$} & $\geq 1$ & & &   & & \\
 \cline{2-3}
 \multirow{2}{4cm}{$\mu=0.1$}&$\alpha_{in}$ & $1.99$ &  & &   &  & \\
 & $\alpha_{out}$ & $1.00$ &  & &  &  & \\
 \midrule
 \multirow{2}{4cm}{$\ell=0.3$} & \multicolumn{1}{l}{$n$} & $1$ &  $2-5$ & $ 6$ & $7$  & $8-9$  & $10$ \\
 \cline{2-8}
\multirow{2}{4cm}{$\mu=0.5$}&$\alpha_{in}$ & $2.00$ & $1.99$ & $2.00$ & $1.99$  & $2.00$  & $1.99$ \\
 &$\alpha_{out}$ & $1.00$ & $1.00$ & $1.00$ & $1.00$  & $1.00$  & $1.00$  \\
 \midrule
\multirow{2}{4cm}{$\ell=0.4$} & \multicolumn{1}{l}{$n$} & $1-3$ &  $4$ & $5-6$ & $7-8$ & $\geq 9$ &  \\
 \cline{2-7}
 \multirow{2}{4cm}{$\mu=0.5$}& $\alpha_{in}$ & $2.00$ & $1.99$ & $2.00$ & $1.99$  & $2.00$ & \\
 &$\alpha_{out}$ & $1.00$ & $1.00$ & $1.00$ & $1.00$ & $1.00$  & \\
 \midrule
\multirow{2}{4cm}{$\omega\equiv1$} & \multicolumn{1}{l}{$n$} & $1$ & $2$  & $\geq 3$ &  &  &  \\
 \cline{2-5}
 \multirow{2}{4cm}{$\mu=0.1$}&$\alpha_{in}$ & $1.46$ & $1.45$ & $1.46$  &   &  & \\
 &$\alpha_{out}$ & $1.00$ & $1.00$  & $1.00$  & &  &  \\
 \bottomrule
\end{tabular}
\end{table}

\begin{table}[H]
\centering
\caption{\label{table:infoP3ME}Reconstructed values $(\alpha_{in},\alpha_{out})$ of the target conductivity $\alpha$ for the mild eccentricity problem.}
  \begin{tabular}{|p{2cm}|l | c c c c c c|}
  \bottomrule
\multicolumn{8}{|c|}{Number of iterations $n=10$} \\
\bottomrule
\multirow{2}{4cm}{$\ell=0.2$}&  \multicolumn{1}{l}{$n$} & $1$ & $\geq2$ & &   & & \\
 \cline{2-4}
 \multirow{2}{4cm}{$\mu=0.1$}& $\alpha_{in}$ & $1.92$ & $1.91$ & &   &  &\\
 & $\alpha_{out}$ & $1.00$ & $1.00$ & &  &  & \\
 \midrule
 \multirow{2}{4cm}{$\ell=0.3$} & \multicolumn{1}{l}{$n$} & $1$ &  $2$ & $3-4$ & $5$ & $6-9$ & $10$ \\
 \cline{2-8}
\multirow{2}{4cm}{$\mu=1$}&  $\alpha_{in}$ & $1.84$ & $1.83$ & $1.82$ & $1.81$ & $1.80$ & $1.81$\\
 &  $\alpha_{out}$ & $1.00$ & $1.00$ & $1.00$ & $1.00$  & $1.00$ & $1.00$\\
 \midrule
\multirow{2}{4cm}{$\ell=0.4$} & \multicolumn{1}{l}{$n$} & $1$ &  $2-3$ & $4$ & $5$ & $\geq 6$ &  \\
 \cline{2-7}
 \multirow{2}{4cm}{$\mu=0.1$}& $\alpha_{in}$ & $1.79$ & $1.78$ & $1.77$ & $1.76$  & $1.77$ & \\
 &$\alpha_{out}$ & $1.00$ & $1.00$ & $1.00$ & $1.00$ & $1.00$  & \\
 \midrule
\multirow{2}{4cm}{$\omega\equiv1$} & \multicolumn{1}{l}{$n$} & $1$ &  $\geq 2$ &  &  &  &  \\
 \cline{2-4}
 \multirow{2}{4cm}{$\mu=0.1$}&$\alpha_{in}$ & $1.02$ & $1.04$ &  &   &  &  \\
 & $\alpha_{out}$ & $1.00$ & $1.00$ &  & &  &\\
 \bottomrule
\end{tabular}
\end{table} 

\begin{table}[H]
\centering
\caption{\label{table:infoPhysical}$BV-$ physical reconstructed values $(\alpha_{in},\alpha_{out})$ of the target conductivity $\alpha$ for every inhomogeneity problem.}
 \begin{tabular}{|l |l |c| c| }
  \bottomrule
\multicolumn{4}{|c|}{Number of iterations $n=1$} \\
\bottomrule
\multicolumn{1}{|l}{{\it Concentric }}& \multicolumn{1}{l}{}  & \multicolumn{1}{l}{$\ell=0.02$} & \multicolumn{1}{l|}{$\omega\equiv1$}  \\
\cline{1-4}
  \multirow{2}{4cm}{ $\mu=1$}& $\alpha_{in}$ & $1.98$ & $1.17$ \\
 & $\alpha_{out}$ & $1.00$ & $1.00$    \\
 \midrule
\multicolumn{1}{|l}{{\it Strong eccentric }}& \multicolumn{1}{l}{}  &\multicolumn{1}{l}{$\ell=0.02$} & \multicolumn{1}{l|}{$\omega\equiv1$} \\
\cline{1-4}
  \multirow{2}{4cm}{ $\mu=1$}& $\alpha_{in}$ & $1.93$ & $1.10$     \\
 & $\alpha_{out}$ & $1.00$ & $1.00$ \\
 \midrule
\multicolumn{1}{|l}{{\it Mild eccentric }}& \multicolumn{1}{l}{}  &\multicolumn{1}{l}{$\ell=0.02$} & \multicolumn{1}{l|}{$\omega\equiv1$}   \\
\cline{1-4}
  \multirow{2}{4cm}{ $\mu=1$}& $\alpha_{in}$ & $1.98$ & $1.00$ \\
 &$\alpha_{out}$ & $1.00$ & $1.00$  \\
 \bottomrule
\end{tabular}
\end{table}

\begin{table}[H]
	\centering
	\caption{\label{table:infoMD}Reconstructed values $(\alpha_{in},\alpha_{out})$ of the target conductivity $\alpha$ for the concentric inhomogeneity problem with $N-$ exact measurements.}
	\begin{tabular}{|p{1.5cm}|c| l | c c c c c c|}
		\bottomrule
		\multicolumn{9}{|c|}{Number of iterations $n=10$} \\
		\bottomrule
		\multirow{4}{2cm}{$\ell=0.2$} & & \multicolumn{1}{l}{$n$} & $1$ & $2$ &$3$ &$\geq 4$ & & \\
		\cline{3-7}
		\multirow{4}{2cm}{$\mu=1$}&\multicolumn{1}{l|}{} & $N=1$ & $1.92$ & $1.94$ &$1.93$ &$1.94$   &  & \\
	 	&\multicolumn{1}{l|}{$\alpha_{in}$}& $N=2$ & $2.00$ & $2.00$ &  $2.00$ &  $2.00$ &  &  \\
        & \multicolumn{1}{l|}{} & $N=5$ & $2.00$ & $2.00$ & $1.99$ &  $1.99$  & &   \\
        \cline{3-7}
        &\multicolumn{1}{l|}{$\alpha_{out}$}& $N=1,2,5$  & $1.00$ & $1.00$ & $1.00$ & $1.00$  & & \\
		\midrule
		\multirow{4}{2cm}{$\ell=0.4$} & & \multicolumn{1}{l}{$n$} & $1$ & $2-4$ & $\geq5$ & & &   \\
		\cline{3-6}
		\multirow{4}{2cm}{$\mu=0.1$}&\multicolumn{1}{l|}{} & $N=1$ & $1.75$ &  $1.76$ & $1.79$ & &   &  \\
		&\multicolumn{1}{l|}{$\alpha_{in}$}& $N=2$ & $1.99$ & $2.00$ & $2.00$ &  &  & \\
		& \multicolumn{1}{l|}{} & $N=5$ & $1.99$ & $1.99$ & $2.00$ &  & &  \\
		\cline{3-6}
	 &\multicolumn{1}{l|}{$\alpha_{out}$}& $N=1,2,5$   & $1.00$ & $1.00$ & $1.00$ &  & & \\
		\midrule
		\multirow{4}{2cm}{$\ell=0.6$} & &  \multicolumn{1}{l}{$n$} & $1-2$ & $3-5$ & $\geq 6$ &   & & \\
		\cline{3-6}
		\multirow{4}{2cm}{$\mu=1$}&\multicolumn{1}{l|}{} & $N=1$ & $1.60$ & $1.62$& $1.60$ & &   &  \\
		&\multicolumn{1}{l|}{$\alpha_{in}$}& $N=2$ & $2.00$ & $2.00$ & $2.00$ &  & & \\
       & \multicolumn{1}{l|}{} & $N=5$ & $2.00$ & $2.00$ & $1.99$ &  &  &  \\
		\cline{3-6}
	&\multicolumn{1}{l|}{$\alpha_{out}$}& $N=1,2,5$  & $1.00$ & $1.00$ & $1.00$ &  & & \\
		\midrule
		\multirow{4}{2cm}{$\omega\equiv1$} & & \multicolumn{1}{l}{$n$} & $1$ & $\geq2$ & & & &  \\
		\cline{3-5}
		\multirow{4}{2cm}{$\mu=1$}&\multicolumn{1}{l|}{} & $N=1$ & $1.12$ & $1.17$ & &   &  & \\
		&\multicolumn{1}{l|}{$\alpha_{in}$}& $N=2$ & $1.09$ & $1.12$ & &  & &  \\
		& \multicolumn{1}{l|}{} & $N=5$ & $1.04$ & $1.06$ & &  & &  \\
		\cline{3-5}
		&\multicolumn{1}{l|}{$\alpha_{out}$}& $N=1,2,5$  & $1.00$ & $1.00$ & &  & &\\
		\bottomrule
	\end{tabular}
\end{table} 

\captionsetup[table]{labelformat=cancaptionlabel}
\begin{table}[H]
	\centering
	\caption{\label{table:infoNMDl02}Reconstructed values $(\alpha_{in}, \alpha_{out})$ of the target conductivity $\alpha$ for $\ell=0.2$ of the concentric inhomogeneity problem with $N-$ ``noisy'' measurements.}
	\begin{tabular}{|p{1.5cm}|c| l | c c c c c c|}
		\bottomrule
		\multicolumn{2}{|c}{$\boldsymbol{\ell=0.2}$} &\multicolumn{7}{|c|}{Number of iterations $n=10$} \\
		\bottomrule
		\multirow{4}{2cm}{$\theta=0.005$} & & \multicolumn{1}{l}{$n$} & $1$ & $2$ &$\geq3$ & & & \\
		\cline{3-6}
		\multirow{4}{2cm}{$\mu=1$}&\multicolumn{1}{l|}{} & $N=1$ & $1.91$ & $1.92$ &$1.92$ &   &  & \\
		&\multicolumn{1}{l|}{$\alpha_{in}$}& $N=2$ & $1.98$ & $1.99$ &  $2.00$ &  &  &  \\
		& \multicolumn{1}{l|}{} & $N=5$ & $1.98$ & $2.00$ & $2.00$ &    & & \\
		\cline{3-6} 
		&\multicolumn{1}{l|}{$\alpha_{out}$}& $N=1,2,5$  & $1.00$ & $1.00$ & $1.00$ &   & & \\
		\midrule
		\multirow{4}{2cm}{$\theta=0.01$} & & \multicolumn{1}{l}{$n$} & $1$ & $2$ & $\geq3$ & & &   \\
		\cline{3-6}
		\multirow{4}{2cm}{$\mu=1$}&\multicolumn{1}{l|}{} & $N=1$ & $1.88$ &  $1.90$ & $1.94$ & &   &  \\
		&\multicolumn{1}{l|}{$\alpha_{in}$}& $N=2$ & $1.94$ & $1.96$ & $1.99$ &  &  & \\
		& \multicolumn{1}{l|}{} & $N=5$ & $1.96$ & $1.99$ & $2.00$ &  & &  \\
		\cline{3-6}
		&\multicolumn{1}{l|}{$\alpha_{out}$}& $N=1,2,5$   & $1.00$ & $1.00$ & $1.00$ &  & & \\
		\midrule
		\multirow{4}{2cm}{$\theta=0.05$} & &  \multicolumn{1}{l}{$n$} & $1-2$ & $3$ & $\geq 4$ &   & & \\
		\cline{3-6}
		\multirow{4}{2cm}{$\mu=1$}&\multicolumn{1}{l|}{} & $N=1$ & $1.87$ & $1.88$& $1.90$ & &   &  \\
		&\multicolumn{1}{l|}{$\alpha_{in}$}& $N=2$ & $2.00$ & $1.98$ & $1.99$ &  & & \\
		& \multicolumn{1}{l|}{} & $N=5$ & $2.00$ & $1.99$ & $2.00$ &  &  &  \\
		\cline{3-6}
		&\multicolumn{1}{l|}{$\alpha_{out}$}& $N=1,2,5$  & $1.00$ & $1.00$ & $1.00$ &  & & \\
		\bottomrule
	\end{tabular}
\end{table}

\addtocounter{table}{-1}
\captionsetup[table]{labelformat=cancaptionlabel2}

\begin{table}[H]
	\centering
	\caption{\label{table:infoNMDl04}Reconstructed values $(\alpha_{in}, \alpha_{out})$ of the target conductivity $\alpha$ for $\ell=0.4$ of the concentric inhomogeneity problem with $N-$ ``noisy'' measurements.}
	\begin{tabular}{|p{1.5cm}|c| l | c c c c c c|}
		\bottomrule
		\multicolumn{2}{|c}{$\boldsymbol{\ell=0.4}$} &\multicolumn{7}{|c|}{Number of iterations $n=10$} \\
		\bottomrule
		\multirow{4}{2cm}{$\theta=0.005$} & & \multicolumn{1}{l}{$n$} & $1-2$ & $3$ &$\geq4$ & & & \\
		\cline{3-6}
		\multirow{4}{2cm}{$\mu=1$}&\multicolumn{1}{l|}{} & $N=1$ & $1.60$ & $1.62$ &$1.67$ &   &  & \\
		&\multicolumn{1}{l|}{$\alpha_{in}$}& $N=2$ & $1.99$ & $1.99$ &  $1.99$ &  &  &  \\
		& \multicolumn{1}{l|}{} & $N=5$ & $1.99$ & $1.99$ & $2.00$ &    & & \\
		\cline{3-6} 
		&\multicolumn{1}{l|}{$\alpha_{out}$}& $N=1,2,5$  & $1.00$ & $1.00$ & $1.00$ &   & & \\
		\midrule
		\multirow{4}{2cm}{$\theta=0.01$} & & \multicolumn{1}{l}{$n$} & $1$ & $2-3$ & $\geq4$ & & &   \\
		\cline{3-6}
		\multirow{4}{2cm}{$\mu=1$}&\multicolumn{1}{l|}{} & $N=1$ & $1.60$ &  $1.65$ & $1.70$ & &   &  \\
		&\multicolumn{1}{l|}{$\alpha_{in}$}& $N=2$ & $1.99$ & $1.98$ & $1.98$ &  &  & \\
		& \multicolumn{1}{l|}{} & $N=5$ & $1.99$ & $2.00$ & $1.99$ &  & &  \\
		\cline{3-6}
		&\multicolumn{1}{l|}{$\alpha_{out}$}& $N=1,2,5$   & $1.00$ & $1.00$ & $1.00$ &  & & \\
		\midrule
		\multirow{4}{2cm}{$\theta=0.05$} & &  \multicolumn{1}{l}{$n$} & $1$ & $2$ & $3$  & $\geq 4$  & & \\
		\cline{3-7}
		\multirow{4}{2cm}{$\mu=1$}&\multicolumn{1}{l|}{} & $N=1$ & $1.56$ & $1.71$& $1.76$ & $1.78$ &   &  \\
		&\multicolumn{1}{l|}{$\alpha_{in}$}& $N=2$ & $1.98$ & $1.99$ & $2.00$ & $1.99$  & & \\
		& \multicolumn{1}{l|}{} & $N=5$ & $1.99$ & $2.00$ & $1.99$ & $2.00$  &  &  \\
		\cline{3-7}
		&\multicolumn{1}{l|}{$\alpha_{out}$}& $N=1,2,5$  & $1.00$ & $1.00$ & $1.00$ &  $1.00$  & & \\
		\bottomrule
	\end{tabular}
\end{table} 

\addtocounter{table}{-1}
\captionsetup[table]{labelformat=cancaptionlabel3}
\begin{table}[H]
	\centering
	\caption{\label{table:infoNMDl06}Reconstructed values $(\alpha_{in}, \alpha_{out})$ of the target conductivity $\alpha$ for $\ell=0.6$ of the concentric inhomogeneity problem with $N-$ ``noisy'' measurements.}
	\begin{tabular}{|p{1.5cm}|c| l | c c c c c c|}
		\bottomrule
		\multicolumn{2}{|c}{$\boldsymbol{\ell=0.6}$} &\multicolumn{7}{|c|}{Number of iterations $n=10$} \\
		\bottomrule
		\multirow{4}{2cm}{$\theta=0.005$} & & \multicolumn{1}{l}{$n$} & $1-6$ & $7-8$ &$9-10$ & & & \\
		\cline{3-6}
		\multirow{4}{2cm}{$\mu=1$}&\multicolumn{1}{l|}{} & $N=1$ & $1.67$ & $1.90$ &$1.99$ &   &  & \\
		&\multicolumn{1}{l|}{$\alpha_{in}$}& $N=2$ & $1.97$ & $1.99$ &  $2.00$ &  &  &  \\
		& \multicolumn{1}{l|}{} & $N=5$ & $1.99$ & $2.00$ & $1.99$ &    & & \\
		\cline{3-6} 
		&\multicolumn{1}{l|}{$\alpha_{out}$}& $N=1,2,5$  & $1.00$ & $1.00$ & $1.00$ &   & & \\
		\midrule
		\multirow{4}{2cm}{$\theta=0.01$} & & \multicolumn{1}{l}{$n$} & $1-3$ & $4-5$ & $6-7$ & $8-10$ & &   \\
		\cline{3-7}
		\multirow{4}{2cm}{$\mu=1$}&\multicolumn{1}{l|}{} & $N=1$ & $1.67$ &  $1.94$ & $1.96$ &$2.00$ &   &  \\
		&\multicolumn{1}{l|}{$\alpha_{in}$}& $N=2$ & $1.99$ & $1.98$ & $1.99$ & $1.98$  &  & \\
		& \multicolumn{1}{l|}{} & $N=5$ & $1.99$ & $2.00$ & $2.00$ & $2.00$  & &  \\
		\cline{3-7}
		&\multicolumn{1}{l|}{$\alpha_{out}$}& $N=1,2,5$   & $1.00$ & $1.00$ & $1.00$ & $1.00$  & & \\
		\midrule
		\multirow{4}{2cm}{$\theta=0.05$} & &  \multicolumn{1}{l}{$n$} & $1$ & $2-3$  & $\geq 4$ & & & \\
		\cline{3-6}
		\multirow{4}{2cm}{$\mu=1$}&\multicolumn{1}{l|}{} & $N=1$ & $1.71$ & $1.90$& $1.96$ &  &   &  \\
		&\multicolumn{1}{l|}{$\alpha_{in}$}& $N=2$ & $1.97$ & $1.98$  & $1.99$ & & & \\
		& \multicolumn{1}{l|}{} & $N=5$ & $1.98$ & $2.00$ & $1.99$ & &  &  \\
		\cline{3-6}
		&\multicolumn{1}{l|}{$\alpha_{out}$}& $N=1,2,5$  & $1.00$ & $1.00$ & $1.00$ &   & & \\
		\bottomrule
	\end{tabular}
\end{table}

\clearpage
\twocolumn
\section{Conclusions}
We presented a methodology for the solution of the inverse conductivity problem, honouring the space of functions of bounded variation as the most appropriate space hosting the conductivity profile. We demonstrated through representative examples that the $BV-$ regularization assists the gradient-based interior point algorithm to achieve accurate reconstructions of both the geometry of the inhomogeneity and the values of the conductivity, unlike the commonly employed for this purpose Tikhonov regularization, which enforces $H^1$ regularity. The $BV-$ regularization reconstruction of the geometry seems to be insensitive to the regularization parameter $\mu$, but this is not the same as far as the estimation of the values of conductivity is concerned. Our investigation revealed that accurate reconstructions of both the geometry of the concentric and eccentric inclusions were possible with rather mild assumptions on their locations. In addition, we incorporated indicative cases of multi-data and noisy measurements clarifying the sensitivity trend of the current reconstruction technique. It is of interest to inve\-sti\-gate further the applicability of the proposed metho\-do\-logy for disconnected inclusions and multi-layered materials. This project is under current investigation.

\numberwithin{equation}{section}
\input{appendix.tex}

\section*{Acknowledgements}	
V.Markaki's graduate studies are supported by a grant from the Papakyriakopoulos Foundation.

\end{document}

%% file: figure1smaller.tex

\begin{tikzpicture}[use Hobby shortcut,closed=true]
  \def\OuterCurve{%
    (-3.4,0.5) .. (-3,2.5) .. (-1,3.55) .. (1.5,3.05) .. (4,2.6) ..
    (4.5,0.5) ..(2.5,-1.2).. (0,-.5).. (-2,-0.9)..(-3.4,0.5)
  }
  \coordinate (Left) at (-3.5, 0.5);
  \coordinate (Right) at (5, 0.5);
  \def\RegionB{%
    (0.5,1) to [out=77, in=117] (3,1) to [out=-60, in=-90] (0.5,1)
    %
    %
  }

  \begin{scope}
     \clip[use Hobby shortcut, closed=true] \RegionB;
     \begin{scope}
      \clip \RegionB
        -- (current bounding box.east |- Right)
        -- (current bounding box.south east)
        -- (current bounding box.south west)
        -- (current bounding box.west |- Left)
        -- cycle
      ;
      \fill[blue] (current bounding box.south west)
        rectangle (current bounding box.north east)
      ;
    \end{scope}
    \clip \RegionB
      -- (current bounding box.east |- Right)
      -- (current bounding box.north east)
      -- (current bounding box.north west)
      -- (current bounding box.west |- Left)
      -- cycle
    ;
    \fill[blue] (current bounding box.south west)
      rectangle (current bounding box.north east)
    ;
\end{scope}

  \draw [->] (-2.25,-1.7) -- (1.40,0.55);
  \draw[use Hobby shortcut, closed=true] \OuterCurve;
  \draw \RegionB;
  
  \node[above] at (-2,0.3) {$\Omega$};
  \node[above] at (-2.5,-2.1) {$D$};
  \node[above] at (3.6, 2.9) {$\partial \Omega$};
  \node[above] at (0.3,2) {$\alpha(x)=\tilde{b}(x)\chi_{\Omega \setminus D}(x)+\tilde{c}(x) \chi_{D}(x)$};
\end{tikzpicture}

%% file: appendix.tex
 \appendix
\section{ The framework of functions of bounded variation} \label{ap1}

Let $\Omega$ be an open subset of $\mathbb{R}^N$ and $\mathcal{B}(\Omega)$ its Borel field. $\mathbf{M}(\Omega,\mathbb{R}^N)$ denotes the space of all $\mathbb{R}^N$-valued Borel measures, which is also according to the Riesz theory the dual of the space $C_0(\Omega,\mathbb{R}^N)$ of all continuous functions $\phi$ vanishing at infinity, equipped with the uniform norm $\Vert \phi\Vert_{\infty}=\big(\sum_{i=1}^{N} \sup_{x\in \Omega}\vert \phi_i(x)\vert^2\big)^{1/2}$. We note that $\mathbf{M}(\Omega,\mathbb{R}^N)$ is isomorphic to the product space $\mathbf{M}^N(\Omega)$ and that $$\mu=(\mu_1,...,\mu_N)\in \mathbf{M}(\Omega,\mathbb{R}^N)\Longleftrightarrow$$
$$\mu_i\in C'_0(\Omega),\;i=1,...,N.$$
\begin{definition}
	We say that a function $u:\Omega\rightarrow\mathbb{R}$ is a function of bounded variation iff it belongs to $L^{1}(\Omega)$ and its gradient $Du$ in the distributional sense belongs to $\mathbf{M}(\Omega,\mathbb{R}^N).$ We denote the set of all functions of bounded variation by $BV(\Omega).$ The four following assertions are then equivalent \cite{butazzo}
	
	\renewcommand{\labelenumi}{(\roman{enumi})}
	\begin{enumerate}
		\item $u\in BV(\Omega)$;
		\item $u\in L^1(\Omega)\;\mbox{and}\;\forall i=1,...,N,\;\frac{\partial u}{\partial x_i}\in \mathbf M(\Omega)$;
		\item $u\in L^1(\Omega)\;\mbox{and}\;\Vert Du\Vert:=\sup \lbrace \langle Du,\phi\rangle : \phi\in C_{c}(\Omega,\mathbb{R}^N),\;\Vert \phi\Vert_{\infty}\leqslant 1 \rbrace < +\infty$;
		\item $u\in L^1(\Omega)\;\mbox{and}\;\Vert Du\Vert=\sup \lbrace \int_\Omega u\;div\; \phi\;dx: \phi\in C_c^1(\Omega,\mathbb{R}^N),\;\Vert \phi\Vert_{\infty}\leqslant 1 \rbrace < +\infty,$
	\end{enumerate}
	where $ C_{c}(\Omega,\mathbb{R}^N)$ denotes the space of all continuous functions with compact support in $\Omega$ and the bracket $<,>$ in (iii) is defined by
	$$\langle Du,\phi\rangle:=\sum_{i=1}^N \int_\Omega \phi_i \frac{\partial u}{\partial x_i}.$$
\end{definition}\

\begin{theorem}[Riesz-Alexandroff representation]\label{RR}
	The topological dual of $C_0(\Omega)$ can be isometrically identified with the space of bounded Borel measures. More precisely, to each bounded linear functional $\Phi$ on $C_0(\Omega)$ there is a unique Borel measure $\mu$ on $\Omega$ such that for all $f\in C_0(\Omega),$
	$$\Phi(f)=\int_{\Omega} f(x)d\mu(x).$$
	Moreover, $\Vert \Phi\Vert=\vert \mu \vert(\Omega).$
\end{theorem}
\begin{remark}
	According to the vectorial version of the Riesz-Alexandroff representation theorem (Theorem \ref{RR}), the dual norm $\Vert Du\Vert$ is also the total mass $\vert Du\vert(\Omega)=\int_{\Omega}\vert Du\vert$ of the total variation $\vert Du\vert$ of the measure $Du.$
\end{remark}

\begin{theorem}[Lebesgue decomposition]\label{RD}
	Let $\mu$ a positive bounded measure on $\big(\mathbb{R}^N,\mathcal{B}(\mathbb{R}^N)\big)$ and $\nu$ a vector-valued measure on $\big(\mathbb{R}^N,\mathcal{B}(\mathbb{R}^N)\big).$ Then there exists a unique pair of measures $\nu_{ac}$ and $\nu_{s}$ such that
	$$\nu=\nu_{ac}+\nu_s,\qquad \nu_{ac}\ll \mu,\qquad \nu_s\perp \mu.$$ Moreover,
	$$\frac{d\nu}{d\mu}=\frac{d\nu_{ac}}{d\mu},\qquad \frac{d\nu_s}{d\mu}=0,\;\mu-a.e.\qquad \mbox{and}$$  $$\nu(A)=\int_A \frac{d\nu}{d\mu}d\mu+\nu_s(A) \qquad \forall A\in \mathcal{B}(\mathbb{R}^N),$$ where $\nu_{ac}$ and $\nu_{s}$ are the absolutely continuous part and the singular part of $\nu.$
	\end{theorem}
	
 \noindent If $u$ belongs to $BV(\Omega)$ and in Theorem \ref{RD} we choose $\mu=dx$, the $N$-dimensional Lebesgue measure, and $\nu=Du,$ we get
\begin{eqnarray}\label{con}
Du=\nabla u \ dx\; + \underbrace{D_su}_{\it{singular\ part}},
\end{eqnarray}
where $\nabla u(x)=\displaystyle\frac{d(Du)}{dx}(x)\in L^1(\Omega)$ and $D_s u\perp dx.$ $\nabla u(x)$ is also called the approximate derivative of u. In fact, we can say more for $BV(\Omega)$ functions. The singular part $D_s u$ of $Du$ can be decomposed into a ``jump" part $J_u$ and a ``Cantor" part $C_u$. In order to specify the $J_u$ part, we need first the following definition 
\begin{definition}[notion of approximate limit]
	Let $B(x,r)$ be the ball of center $x$ and radius $r$ and let $u\in BV(\Omega)$. The approximate upper limit $u^{+}(x)$ and the approximate lower limit $u^{-}(x)$ are defined by
	$$u^{+}(x)=inf\Big\lbrace t\in [-\infty,+\infty];\ \lim_{r\to 0}\frac{dx(\lbrace u>t\rbrace\cap B(x,r))}{r^N}=0\Big\rbrace,$$
	$$u^{-}(x)=sup\Big\lbrace t\in [-\infty,+\infty];\ \lim_{r\to 0}\frac{dx(\lbrace u<t\rbrace\cap B(x,r))}{r^N}=0\Big\rbrace.$$
 If $u\in L^{1}(\Omega)$ then
	\begin{eqnarray}\label{lp}
	\lim_{r\to 0}\frac{1}{\vert B(x,r)\vert}\int_{B(x,r)}\vert u(x)-u(y)\vert dy=0\;\;\;a.e.\; x.
	\end{eqnarray}
\end{definition}
A point $x$ for which (\ref{lp}) holds is called a Lebesgue point of u, and we have
\begin{eqnarray}
\lim_{r\to 0}\frac{1}{\vert B(x,r)\vert}\int_{B(x,r)}u(y)dy,
\end{eqnarray}
and $u(x)=u^{+}(x)=u^{-}(x)$. We denote by $S_u$ the jump set, that is, the complement, up to a set of $\mathcal{H}^{N-1}$\footnote{Hausdorff measure} measure zero, of the set of Lebesgue points
$$S_u=\lbrace x\in \Omega;u^{-}(x)<u^{+}(x)\rbrace.$$
Then $S_u$ is countably rectifiable, and for $\mathcal{H}^{N-1}$- a.e.\;$x\in \Omega,$ we can define a normal $n_u(x).$\hfill \break
Therefore, the measure $Du$ admits the following representation \cite{aubert}
\begin{eqnarray}\label{open}
Du=\nabla u dx+\underbrace{(u^{+}-u^{-})n_u\mathcal{H}^{N-1}_{\vert_{S_u}}}_{\it{jump\ part}}\;+ \underbrace{C_u.}_{\it{Cantor\ part}}
\end{eqnarray}
From (\ref{open}) we can deduce the total variation of $Du:$
\begin{align}
 \vert Du\vert (\Omega)&=\int_\Omega \vert Du\vert\nonumber\\ &=
\int_\Omega \vert \nabla u\vert(x) dx  \nonumber \\
& +\int_{S_u} \vert u^{+}-u^{-}\vert d\mathcal{H}^{N-1}+\int_{\Omega-S_u} \vert C_u\vert.
\end{align}

From (\ref{con}) we conclude that $W^{1,1}(\Omega)$ is a subspace of the vectorial space $BV(\Omega)$ and $u\in W^{1,1}(\Omega)$ iff $Du=\nabla u \mathcal{L}^N \lfloor \Omega.$ The space $BV(\Omega)$ is equipped with the following norm, which extends the classical norm in $W^{1,1}(\Omega)$:
$$\Vert u \Vert _{BV(\Omega)}:=\vert u\vert_{L^1(\Omega)} + \Vert Du \Vert.$$
Equipped with its norm, $BV(\Omega)$ is a Banach space.\

We will define two weak convergence processes in $BV(\Omega)$. The first provides compactness of bounded sequences and the second is an intermediate convergence between the weak and the strong convergence associated with the norm.
\begin{definition}
	A sequence $(u_n)_{n\in\mathbb{N}}$ in $BV(\Omega)$ weakly (\footnotemark or weakly-$\star$) converges to some $u$ in $BV(\Omega),$ and we write $u_n \rightharpoonup u$ (or $u_n \overset{\star}{\rightharpoonup} u$) iff the following convergences hold:\footnotetext{In $BV(\Omega),$ ones can refer to the weak convergence by denoting $w$ or $w^{*}$ equivalently}
	\begin{equation}
	\begin{cases}
	u_n\rightarrow u\;\mbox{strongly \ in}\;L^1(\Omega);\\
	Du_n \rightharpoonup Du\;\mbox{weakly \ in}\;\mathbf{M}(\Omega,\mathbb{R}^N).\nonumber
	\end{cases}
	\end{equation}
\end{definition}
In the proposition below we establish a compactness result related to this convergence, together with the lower semicontinuity of the total mass.
\begin{proposition}
	Let $(u_n)_{n\in\mathbb{N}}$ be a sequence in $BV(\Omega)$ strongly converging to some $u$ in $L^{1}(\Omega)$ and satisfying ${\sup_{n\in\mathbb{N}}}\int_{\Omega}\vert Du_n\vert <+\infty.$ Then
	\renewcommand{\labelenumi}{(\roman{enumi})}
	\begin{enumerate}
		\item $u\in BV(\Omega)$ and $\int_{\Omega}\vert Du\vert\leqslant \underset{n\to+\infty}{\underline{\lim}}\int_{\Omega}\vert Du_n\vert;$
		\item $u_n$ weakly converges to $u$ in $BV(\Omega).$
	\end{enumerate}
\end{proposition}

\begin{definition}
	Let $(u_n)_{n\in\mathbb{N}}$  be a sequence in $BV(\Omega)$ and $u\in BV(\Omega).$ We say that $u_n$ converges to $u$ in the sense of the intermediate convergence iff
	\begin{equation}
	\begin{cases}
	u_n\rightarrow u\;\mbox{strongly \ in}\;L^1(\Omega);\nonumber\\
	\int_\Omega\vert Du_n \vert \rightarrow \int_\Omega\vert Du\vert.
	\end{cases}
		\end{equation}\nonumber
\end{definition}
Finally, we state some fundamental theorems \cite{butazzo}
\begin{theorem}
	The space $C^{\infty}(\Omega)\cap BV(\Omega)$ is dense in $BV(\Omega)$ equipped with the intermediate convergence. Consequently $C^{\infty}(\overline{\Omega})$ is also dense in $BV(\Omega)$ for the intermediate convergence.
\end{theorem}

\begin{theorem}
	Let $\Omega$ be a 1-regular open bounded subset of $\mathbb{R}^N.$ For all $p,\;1\leqslant p \leqslant \frac{N}{N-1},$ the embedding
	$$BV(\Omega)\hookrightarrow L^p(\Omega)$$ is continuous. More precisely, there exists a constant $C$ which depends only on $\Omega, p,$ and $N,$ such that for all $u$ in $BV(\Omega),$
	$$ \Big( \int_\Omega \vert u\vert^p dx\Big)^{\frac{1}{p}}\leqslant C\vert u\vert_{BV(\Omega)}.$$
\end{theorem}

\begin{theorem}
	Let $\Omega$ be a 1-regular open bounded subset of $\mathbb{R}^N.$ Then for all $p,\;1\leqslant p < \frac{N}{N-1},$ the embedding $$BV(\Omega)\hookrightarrow L^p(\Omega)$$ is compact.
\end{theorem}

\begin{theorem}
	Let $\Omega$ be a domain of $\mathbb{R}^N$ with a Lipschitz boundary $\Gamma$. There exists a linear continuous map $\gamma_0$ from $BV(\Omega)$ onto $L_{\mathcal{H}^{N-1}}^{1}(\Gamma)$ satisfying
	\renewcommand{\labelenumi}{(\roman{enumi})}
	\begin{enumerate}
		\item for all $u$ in $C(\overline{\Omega})\cap BV(\Omega), \gamma_0(u)=u\lfloor\Gamma;$
		\item the generalized Green's formula holds: $\forall \phi \in C^{1}(\overline{\Omega},\mathbb{R}^N),$
		$$\int_{\Omega}\phi Du=-\int_{\Omega}u\; div\; \phi\; dx + \int_{\Gamma}\gamma_0(u)\phi .\nu\;d\mathcal{H}^{N-1},$$
		where $\nu(x)$ is the outer unit normal at $\mathcal{H}^{N-1}$ almost all $x$ in $\Gamma.$
	\end{enumerate}
\end{theorem}

\section{Results concerning the homogenization process} \label{ap2} 
\bfseries\large{Proof of Lemma \ref{prop1}} \\
\normalfont\normalsize \par
Given that $u_n^{\alpha_n} \underset{n \rightarrow \infty}{-\!\!\!\rightharpoonup} v$ in $H^1_0(\Omega)$, there exists a vector function $\sigma$ such that, up to a subsequence still denoted by $n$, $\alpha_n(x) \nabla u_n^{\alpha_n} \underset{n \rightarrow \infty}{-\!\!\!\rightharpoonup} \sigma $ (\mbox{weakly \ in} ${\left( L^2(\Omega)\right)}^N$). For $\lambda \in \mathbb{R}^N$, we set $w(x)=\lambda \cdot x \phi(x)$, where $\phi \in C_0^\infty(\Omega)$ and $g(x)=-\nabla \cdot (\alpha \nabla w(x))$. Then $w_n$ is defined to be the unique solution of the problem
\begin{align}
&-\nabla \cdot (\alpha_n(x) \nabla w_n(x))=g(x) \qquad && \mbox{in} \  \Omega, \\
&w_n(x)=0 \qquad &&\mbox{on} \ \partial \Omega.
\end{align}
Since $\alpha_n(x) \mathbb{I}_{N \times N}\underset{H}{\longrightarrow} \alpha(x) \mathbb{I}_{N \times N}$, we have that $w_n \underset{H_0^1(\Omega)}{-\!\!\!\rightharpoonup} w$ and $\alpha_n \nabla w_n \underset{{(L_2(\Omega))}^N}{-\!\!\!\rightharpoonup} \alpha \nabla w$ (We just mention that $w$ satisfies the limit problem $-\nabla \cdot (\alpha \nabla w)=g$, $\gamma w=0$ in a trivial manner). By coercivity of $\alpha_n$, we have that 
$(\alpha_n(x)\nabla u_n^{\alpha_n}(x)-\alpha_n(x)\nabla w_n(x)) \cdot (\nabla u_n^{\alpha_n}(x)- \nabla w_n(x)) \geq 0$ a.e. in $\Omega$. We select an arbitrary open set $\omega$, proper subset of $\Omega$ and restrict ourselves to $\phi \in C_0^\infty(\Omega)$, with $\phi{|}_\omega=1$. Then $\nabla w (x)=\lambda, \ x\in \omega$. By construction of the involved fields, the assumptions of the div-curl Lemma 1.3.1 of \cite{allaire} are fulfilled and so we can pass to the limit in the last inequality, obtaining that $(\sigma-\alpha \lambda) \cdot (\nabla v-\lambda) \geq 0$ a.e. in $\omega$.  We select a point $x_0 \in \omega$, where strict inequality holds and choose $\lambda=\nabla v (x_0) +t \mu$, where $t>0$ and $\mu \in {\mathbb{R}}^{N}$. Dividing by $t$, taking the limit $t \rightarrow 0$ and exploiting the arbitrariness of $\mu$, we obtain $\sigma(x_0)=\alpha(x_0) \nabla v(x_0)$, a.e. in $\omega$, which has been considered as an arbitrary subset of $\Omega$, invoked just to avoid vicinity with the boundary $\partial \Omega$. Since the weak limit $\sigma$ is determined uniquely and independently of the initial subsequence, we infer that the entire sequence $\alpha_n(x) \nabla u_n^{\alpha_n}$ weakly converges to $\sigma$ in $(L^2(\omega))^N$. Since $\omega$ is any subset of $\Omega$, this means that $\alpha_n(x) \nabla u_n^{\alpha_n} \underset{n \rightarrow \infty}{-\!\!\!\rightharpoonup} \alpha \nabla v$ \mbox{ weakly \ in} ${\left(L^2(\Omega)\right)}^N$.     \\ \par
\bfseries\large{Proof of Lemma \ref{prop2}} \\
\normalfont\normalsize \par
The variational form of the problems I and (\ref{ouf1}-\ref{ouf2}) can be exploited to lead easily to the following relation
\begin{align}
&\int_{\Omega}\alpha_n \nabla (u_n^{\alpha_n}-{\hat{u}}_n^{\alpha_n}) \cdot \nabla z =\nonumber\\ &\left\langle h_n-h, z \right\rangle_{H^{-1}(\Omega) \times H^{1}_0(\Omega)},\ \forall z \in  H^{1}_0(\Omega). 
\end{align}
Taking $z=u_n^{\alpha_n}-{\hat{u}}_n^{\alpha_n}$ and exploiting the coerciveness of the coefficient $\alpha_n$, we obtain
\begin{eqnarray}
b\|u_n^{\alpha_n}-{\hat{u}}_n^{\alpha_n}{\|}_{H^1_0(\Omega)} \leq \|h_n-h{\|}_{H^{-1}(\Omega)}.  \label{ordinaz}
\end{eqnarray}
Due to H-convergence, the sequence $\hat{u}_n^{\alpha_n}$ has as weak limit (in $H_0^1(\Omega)$) the function $u^{\alpha}$. Thanks to the strong convergence $h_n \rightarrow h$, we infer from (\ref{ordinaz}) that the sequences $u_n^{\alpha_n}$ and $\hat{u}_n^{\alpha_n}$ merit the same weak limit ($u^\alpha$) in $H_0^1(\Omega)$. Using again the strong convergence $h_n \rightarrow h$ (in $H^{-1}(\Omega)$), we find that \cite{brezis}
\begin{eqnarray}
\left\langle h_n, {u}_n^{\alpha_n} \right\rangle_{H^{-1}(\Omega) \times H^{1}_0(\Omega)} \underset{n \rightarrow \infty}{-\!\!\!\rightarrow} \left\langle h, u^{\alpha} \right\rangle_{H^{-1}(\Omega) \times H^{1}_0(\Omega)}. \label{patienz}
\end{eqnarray}
Taking advantage of the variational form of the involved problems, (\ref{patienz}) is another expression of the following relation 
\begin{eqnarray}
\int_{\Omega}\alpha_n |\nabla u_n^{\alpha_n}{|}^2 \underset{n \rightarrow \infty}{-\!\!\!\rightarrow} \int_{\Omega}\alpha |\nabla u^{\alpha}{|}^2 . \label{lud}
\end{eqnarray} \\ \par
\bfseries\large{Proof of Proposition \ref{prop3}} \\
\normalfont\normalsize \par
As already stated, when $f \in H^{\frac{1}{2}}(\partial \Omega)$, and on the basis of boundedness of the conductivity coefficients, it is easily proved that $h_n = \nabla \cdot(\alpha_n(x) \nabla (\eta f)(x) ) \in H^{-1}(\Omega)$. In addition, it is a straightforward matter to show that $h_n \rightarrow h=\nabla \cdot(\alpha(x) \nabla (\eta f)(x) )$ \mbox{strongly\ in} $H^{-1}(\Omega)$, due mainly to the $L^1(\Omega)$ - convergence of $\alpha_n$. \\
The variational form of Problem I yields
\begin{align}
&\int_{\Omega}\alpha_n \nabla u_n^{\alpha_n} \cdot \nabla z =\nonumber\\ &\left\langle h_n, z \right\rangle_{H^{-1}(\Omega) \times H^{1}_0(\Omega)},\ \forall z \in  H^{1}_0(\Omega). \label{eks}
\end{align}
Taking $z=u_n^{\alpha_n}$ and exploiting the lower boundedness of $\alpha_n$ and the uniform boundedness of $\|h_n\|_{H^{-1}}$, we easily show that $\|u_n^{\alpha_n}\|_{H^1_0} \leq C$. Then, we extract a subsequence, still denoted $u_n^{\alpha_n}$, such that $u_n^{\alpha_n} \rightharpoonup v$, \mbox{weakly\ in} $H^1_0(\Omega)$, for some element $v \in H^1_0(\Omega)$. According to Lemma \ref{prop1}, it also holds that $\alpha_n(x) \nabla {u}_n^{\alpha_n}(x)  \rightharpoonup  \alpha(x) \nabla v(x)\ \mbox{weakly\ in}\ {(L^{2}(\Omega))}^N$. Taking then the limit ($n \rightarrow \infty$) to the expression (\ref{eks}), we obtain
\begin{eqnarray}
\int_{\Omega}\alpha \nabla v \cdot \nabla z = \left\langle h, z \right\rangle_{H^{-1}(\Omega) \times H^{1}_0(\Omega)}, \ \forall z \in  H^{1}_0(\Omega). \label{eks2}
\end{eqnarray}
But the variational form of the limiting problem (\ref{oriako1})-(\ref{oriako2}) reads
\begin{eqnarray}
\int_{\Omega}\alpha \nabla u^\alpha \cdot \nabla z = \left\langle h, z \right\rangle_{H^{-1}(\Omega) \times H^{1}_0(\Omega)}, \ \forall z \in  H^{1}_0(\Omega). \label{eks3}
\end{eqnarray}
Substracting (\ref{eks2}) and (\ref{eks3}), putting $z= v-u^\alpha$ and exploiting the positiveness of $\alpha$, we obtain the desired result $v=u^\alpha$. Every weakly converging subsequence of the family $u_n^{\alpha_n}$ converges necessarily to $u^\alpha$ and then the whole sequence makes so. Thus, we have already shown that 
\begin{align}
{u}_n^{\alpha_n}(x)  &\rightharpoonup   {u}^{\alpha}(x) \ \mbox{weakly \ in} \ H^{1}_0(\Omega) \label{1965} \\
\alpha_n(x) \nabla {u}_n^{\alpha_n}(x)  &\rightharpoonup  \alpha(x) \nabla {u}^{\alpha}(x)\ \mbox{weakly \ in} \ {(L^{2}(\Omega))}^N. \label{1966}
\end{align}
Due to Lemma \ref{prop2}, we have the validity of (\ref{lud})
\begin{eqnarray*}
	\int_{\Omega}\alpha_n |\nabla u_n^{\alpha_n}{|}^2 \underset{n \rightarrow \infty}{-\!\!\!\rightarrow} \int_{\Omega}\alpha |\nabla u^{\alpha}{|}^2 .
\end{eqnarray*}
For surface data $f \in H^{\frac{1}{2}}(\partial \Omega)$ we have that $\eta f \in H^{1}(\Omega)$ and that $\nabla (\eta f) \in {(L_2(\Omega))}^N$. Then according to (\ref{1966}), we obtain
\begin{eqnarray}
\int_{\Omega}\alpha_n \nabla u_n^{\alpha_n} \cdot \nabla (\eta f) \rightarrow \int_{\Omega}\alpha \nabla u^{\alpha} \cdot \nabla (\eta f). \label{aux1}
\end{eqnarray}  
The convergence $\alpha_n \underset{L^{\infty}-w^{\ast}}{-\!\!\!\rightharpoonup} \alpha$ and the relation ${|\nabla (\eta f)|}^2 \in L^1(\Omega)$ imply that 
\begin{eqnarray}
\int_{\Omega}\alpha_n {|\nabla (\eta f)|}^2  \rightarrow \int_{\Omega}\alpha {|\nabla (\eta f)|}^2. \label{aux2}
\end{eqnarray}
Using (\ref{lud}, \ref{aux1}, \ref{aux2}), we find that 
\begin{align}
&\int_{\Omega}\alpha_n|\nabla (u_n^{\alpha_n}+\eta f)|^2 = \nonumber \\ &\int_{\Omega}\alpha_n|\nabla u_n^{\alpha_n}|^2 +\int_{\Omega}\alpha_n|\nabla (\eta f)|^2+2\int_{\Omega}\alpha_n \nabla u_n^{\alpha_n} \cdot \nabla (\eta f) \nonumber \\
 &\underset{n \rightarrow \infty}{-\!\!\!\rightarrow} \int_{\Omega}\alpha|\nabla u^{\alpha}|^2 +\int_{\Omega}\alpha|\nabla (\eta f)|^2+2\int_{\Omega}\alpha \nabla u^{\alpha} \cdot \nabla (\eta f)\nonumber \\ &= \int_{\Omega}\alpha|\nabla (u^{\alpha}+\eta f)|^2.
\end{align}
\\
\bfseries\large{Proof of Proposition \ref{prop4}} \\
\normalfont\normalsize \par
We make the decomposition
\begin{eqnarray}
w_n^{\alpha_n}=v_n^{\alpha_n}+\eta \Lambda_{\alpha_n}^{-1}g, 
\end{eqnarray}
where we evoke the family of Neumann-to-Dirichlet operators $\Lambda_{\alpha_n}^{-1}$ (or $\Lambda_{N.t.D}$), which are well defined, uniformly bounded operators (over the integer $n$), since the members of the sequence of conductivities $\alpha_n$ are uniformly bounded below, given that $\alpha_n, \alpha \in L^\infty(\Omega,[b,c])$.  
It is straightforward to show that the field $v_n^{\alpha_n}$ satisfies the boundary value problem
\begin{eqnarray}
-\nabla \cdot (\alpha_n(x) \nabla v_n^{\alpha_n}(x))=\tilde{h}_n(x) \qquad &\mbox{in} \  \Omega, \label{prig1}\\
 v_n^{\alpha_n}(x)=0 \qquad &\mbox{on} \ \partial \Omega,
\end{eqnarray}
where $\tilde{h}_n=\nabla \cdot (\alpha_n \nabla (\eta \Lambda_{\alpha_n}^{-1} g))$ $\in H^{-1}(\Omega)$. \\
For readers' convenience, we briefly clarify here the following notations (see extensively in \ref{icp}): $\Lambda_{\alpha}=\Lambda_{D.t.N} \ \mbox{with}\ \Lambda_{\alpha}: \mathcal{B}\left(:=H^{\frac{1}{2}}(\partial \Omega)/\mathbb{R}\right) \rightarrow {\mathcal{B}}^{\ast} \ \mbox{and}\ \Lambda^{-1}_{\alpha}= \Lambda_{N.t.D} \ \mbox{with}\ \Lambda^{-1}_{\alpha}: {\mathcal{B}}^{\ast} \rightarrow \mathcal{B}.$ \newline At this point, we exploit the main result of \cite{Ruiz}, according to which when a very weak assumption is made for the conductivity profiles $\alpha_n$ then the strong convergence ${\|\Lambda_{\alpha_n} - \Lambda_\alpha\|}_{{\mathcal{B}}} \underset{n \rightarrow \infty}{-\!\!\!\rightarrow} 0$ is valid. Actually, the mentioned weak assumption requires a specific convergence of the profiles $\alpha_n$ to $\alpha$ in a very thin region around the boundary $\partial \Omega$. In this work, this assumption is a fortiriori trivially satisfied, given that $\alpha_n, \alpha$ belong to $L_\delta^\infty(\Omega,\left[ b,c \right])$ (Section \ref{func}) and then all the profiles totally agree on a very small region $\Omega_{\delta}=\left\{x \in \Omega : \mbox{dist}(x, \partial \Omega) < \delta \right\}$, with $\delta << 1$. Given that $\Lambda_{\alpha_n}^{-1} - \Lambda_\alpha^{-1}=\Lambda_{\alpha_n}^{-1}(\Lambda_\alpha - \Lambda_{\alpha_n})\Lambda_{\alpha}^{-1}$, it is easily inferred - thanks to the uniform boundedness of $\Lambda_{\alpha_n}^{-1}$, $\Lambda_{\alpha}^{-1}$ - that the strong convergence ${\|\Lambda_{\alpha_n}^{-1} - \Lambda_\alpha^{-1}\|}_{{\mathcal{B}}^{\ast}} \underset{n \rightarrow \infty}{-\!\!\!\rightarrow} 0$ is also valid.  \par
The last result allows us to take the limit in $\tilde{h}_n$ to prove that $\tilde{h}_n\underset{n \rightarrow \infty}{-\!\!\!\rightarrow} \tilde{h}:=\nabla \cdot (\alpha \nabla (\eta \Lambda_{\alpha}^{-1} g))$ $\in H^{-1}(\Omega)$ \mbox{strongly in} $H^{-1}(\Omega)$. \par Indeed, given that $\Lambda_{\alpha_n}^{-1}g \underset{n\rightarrow \infty}{-\!\!\!\rightarrow} \Lambda_{\alpha}^{-1}g$ (\mbox{strongly in} ${\mathcal{B}}^{\ast}$), we infer that $\nabla (\eta \Lambda_{\alpha_n}^{-1} g) \cdot \nabla z$ $\underset{L^1(\Omega)}{-\!\!\!\rightarrow}$ $\nabla (\eta \Lambda_{\alpha}^{-1} g) \cdot \nabla z$, for every $z \in H^1_0(\Omega)$. We already know that $\alpha_n \underset{L^{\infty}-w^{\ast}}{-\!\!\!\rightarrow} \alpha$ and then
\begin{align}
&\left\langle \tilde{h}_n, z \right\rangle_{H^{-1}(\Omega) \times H^{1}_0(\Omega)} = \int_{\Omega} \alpha_n \nabla (\eta \Lambda_{\alpha_n}^{-1} g) \cdot \nabla z \nonumber \\ &\underset{n\rightarrow \infty}{-\!\!\!\rightarrow} \int_{\Omega} \alpha \nabla (\eta \Lambda_{\alpha}^{-1} g) \cdot \nabla z =\left\langle \tilde{h}, z \right\rangle_{H^{-1}(\Omega) \times H^{1}_0(\Omega)}.\nonumber \\
\end{align}
\par According to Lemma \ref{prop2} and Proposition \ref{prop3}, we obtain that 
\begin{align}
&v_n^{\alpha_n} \underset{n \rightarrow \infty}{-\!\!\!\rightharpoonup} v^{\alpha} \ \mbox{weakly \ in}\ H^1_0(\Omega), \\
&\alpha_n \nabla v_n^{\alpha_n} \underset{n \rightarrow \infty}{-\!\!\!\rightharpoonup} \alpha \nabla v^{\alpha}  \ \mbox{weakly \ in}\ {(L_2(\Omega))}^N,  \\
&\int_{\Omega} \alpha_n {|\nabla v_n^{\alpha_n}|}^2 \underset{n \rightarrow \infty}{-\!\!\!\rightarrow} \int_{\Omega} \alpha {|\nabla v^{\alpha}|}^2, \label{finaal1}
\end{align}
where $v^{\alpha}$ is the solution of the limiting problem
\begin{align}
&-\nabla \cdot (\alpha(x) \nabla v^{\alpha}(x))=\tilde{h}(x) \qquad && \mbox{ in} \  \Omega, \\
&v^{\alpha}(x)=0 \qquad &&\mbox{on} \ \partial \Omega.
\end{align}
\noindent In addition
\begin{eqnarray}
\int_{\Omega} \alpha_n {|\nabla (\eta \Lambda^{-1}_{\alpha_n} g)|}^2 \underset{n\rightarrow \infty}{-\!\!\!\rightarrow} \int_{\Omega} \alpha {|\nabla (\eta \Lambda^{-1}_{\alpha} g)|}^2, \label{finaal2}
\end{eqnarray} 
due again to the limits $\alpha_n \underset{L^{\infty}-w^{\ast}}{-\!\!\!\rightarrow} \alpha$ and ${|\nabla (\eta \Lambda^{-1}_{\alpha_n} g)|}^2 \underset{L^1(\Omega) - \mbox{strongly}}{-\!\!\!\rightarrow} {|\nabla (\eta \Lambda^{-1}_{\alpha} g)|}^2$.
\newline Finally, 
\begin{eqnarray}
\int_{\Omega} \alpha_n \nabla v_n^{\alpha_n} \cdot \nabla (\eta \Lambda^{-1}_{\alpha_n} g) \underset{n\rightarrow \infty}{-\!\!\!\rightarrow} \int_{\Omega} \alpha \nabla v^{\alpha} \cdot \nabla (\eta \Lambda^{-1}_{\alpha} g), \label{finaal3}
\end{eqnarray} 
since $\alpha_n \nabla v_n^{\alpha_n} \underset{n \rightarrow \infty}{-\!\!\!\rightharpoonup} \alpha \nabla v^{\alpha}, \ \mbox{weakly \ in}\ {(L_2(\Omega))}^N$ and $ \nabla (\eta \Lambda^{-1}_{\alpha_n} g) \underset{n \rightarrow \infty}{-\!\!\!\rightarrow}  \nabla (\eta \Lambda^{-1}_{\alpha} g)$,  $\mbox{strongly \ in}\ {(L_2(\Omega))}^N$. 
\newline Consequently, the sequence $w_n^{\alpha_n}$ converges (\mbox{weakly \ in} $H^1_0(\Omega)$) to $w^{\alpha}=v^{\alpha}+\eta \Lambda_{\alpha}^{-1}g$, which readily solves the problem (\ref{firstfinal})-(\ref{secondfinal}) and furthermore, on the basis of limit processes (\ref{finaal1}), (\ref{finaal2}) and (\ref{finaal3}), it comes out that  
\begin{align}
 &\int_{\Omega} \alpha_n {|\nabla w_n^{\alpha_n}|}^2 =\int_{\Omega} \alpha_n {|\nabla (\eta \Lambda^{-1}_{\alpha_n} g)|}^2\nonumber \\ &+ 2 \int_{\Omega} \alpha_n \nabla v_n^{\alpha_n} \cdot \nabla (\eta \Lambda^{-1}_{\alpha_n} g)+\int_{\Omega} \alpha_n {|\nabla v_n^{\alpha_n}|}^2 \nonumber \\ &\underset{n\rightarrow \infty}{-\!\!\!\rightarrow} \int_{\Omega} \alpha {|\nabla (\eta \Lambda^{-1}_{\alpha} g)|}^2+2 \int_{\Omega} \alpha \nabla v^{\alpha} \cdot \nabla (\eta \Lambda^{-1}_{\alpha} g)\nonumber \\ &+ \int_{\Omega} \alpha {|\nabla v^{\alpha}|}^2 = \int_{\Omega} \alpha {|\nabla w^{\alpha}|}^2.
\end{align}

\newpage